\documentclass[final]{siamltex}

\usepackage{epsfig,amssymb,latexsym}
\usepackage{amsfonts,psfrag,amsmath,bbm,color}
\usepackage{cancel}
\usepackage{comment}
\usepackage{mathrsfs}
\usepackage{graphicx}
\usepackage{textcomp}
\usepackage{multirow}
\usepackage{multicol}
\usepackage{enumerate}
\usepackage{cancel}
\usepackage{algpseudocode}
\usepackage{caption}  
\usepackage{subcaption}
\usepackage{url}
\usepackage{rotating}
\usepackage{bigints}
\usepackage{hyperref}
\usepackage{xcolor}
\usepackage{ulem}

\usepackage{placeins}

\def\grad{{\nabla}}
\usepackage[ruled,vlined]{algorithm2e}
\usepackage{geometry}
\vbadness=\maxdimen

\newcommand{\footremember}[2]{%
	\footnote{#2}
	\newcounter{#1}
	\setcounter{#1}{\value{footnote}}%
}

\usepackage{tikz}
\usetikzlibrary{shapes.geometric, arrows}
\tikzstyle{startstop} = [rectangle, rounded corners, minimum width=1cm, minimum height=1cm,text centered, draw=black]
\tikzstyle{io} = [trapezium, trapezium left angle=70, trapezium right angle=110, minimum width=1cm, minimum height=1cm, text centered, 
draw=black, fill=blue!30]
\tikzstyle{method} = [rectangle, rounded corners, minimum width=1cm, minimum height =1cm, text centered, draw=black]
\tikzstyle{process} = [rectangle, minimum width=1cm, minimum height=1cm, text centered, draw=black]
\tikzstyle{decision} = [diamond, minimum width=0.5cm, minimum height=0.5cm, text centered, draw=black, fill=green!30]
\tikzstyle{arrow} = [thick,->,>=stealth]

\usepackage{amscd}

\graphicspath{{./figs/}}

\newcommand{\lp}{\left(}
\newcommand{\rp}{\right)}

\newcommand{\lab}{<\hspace{-1mm}}
\newcommand{\rab}{\hspace{-1mm}>}

\newtheorem{remark}{Remark}[section]
\newtheorem{assumption}{Assumption}[section]

\def\PP{{{\rm l}\kern - .15em {\rm P} }}
\def\PN2{{\PP_{N}-\PP_{N-2}}}

\newcommand{\E}{\mathbbm{E}}

\newcommand{\cD}{\mathcal{D}}

\newcommand{\bphi}{\boldsymbol{\varphi}}

\newcommand{\bif}{\textbf{\textit{f}}}
\newcommand{\ba}{\boldsymbol{a}}

\newcommand{\bb}{\boldsymbol{b}}
\newcommand{\bB}{\boldsymbol{B}}

\newcommand{\bchi}{\pmb{\chi}}

\newcommand{\be}{\boldsymbol{e}}
\newcommand{\bep}{\boldsymbol{\epsilon}}

\newcommand{\bg}{\textbf{\textit{g}}}

\newcommand{\bH}{\boldsymbol{H}}
\newcommand{\bs}{{\bf s}}
\newcommand{\bl}{\boldsymbol{l}}
\newcommand{\bL}{\boldsymbol{L}}

\newcommand{\hq}{{\hat{q}}}

\newcommand{\bR}{\boldsymbol{R}}

\newcommand{\bu}{\boldsymbol{u}}

\newcommand{\bv}{\boldsymbol{v}}
\newcommand{\btv}{\tilde{\boldsymbol{v}}}
\newcommand{\bhu}{\hat{\boldsymbol{u}}}
\newcommand{\bhv}{\hat{\boldsymbol{v}}}
\newcommand{\bnh}{\hat{\textbf{\textit{n}}}}
\newcommand{\hlam}{\hat{\lambda}}
\newcommand{\hr}{\hat{r}}
\newcommand{\bV}{\boldsymbol{V}}

\newcommand{\bw}{\boldsymbol{w}}
\newcommand{\bW}{\boldsymbol{W}}
\newcommand{\bhw}{\hat{\boldsymbol{w}}}
\newcommand{\btw}{\tilde{\boldsymbol{w}}}

\newcommand{\bx}{\boldsymbol{x}}
\newcommand{\bX}{\boldsymbol{X}}

\newcommand{\bY}{\boldsymbol{Y}}
\newcommand{\by}{\boldsymbol{y}}
\newcommand{\bz}{\boldsymbol{z}}
\newcommand{\bhz}{\hat{\boldsymbol{z}}}
\newcommand{\bGamma}{\boldsymbol{\Gamma}}
\newcommand{\bxi}{\boldsymbol{\xi}}

\newcommand{\deleted}[1]{{}}

\begin{document}
	\title{A penalty-projection based efficient and accurate Stochastic Collocation Method for magnetohydrodynamic flows}

	\author{
		Muhammad Mohebujjaman\footremember{mit}{D\MakeLowercase{epartment of} M\MakeLowercase{athematics}, U\MakeLowercase{niversity of} 
			A\MakeLowercase{labama at} B\MakeLowercase{irmingham}, AL 35294, USA; T\MakeLowercase{his author's work was Partially supported by the} N\MakeLowercase{ational} 
			S\MakeLowercase{cience} F\MakeLowercase{oundation grant} DMS-2213274, \MakeLowercase{and} T\MakeLowercase{exas} A\&M 
			I\MakeLowercase{nternational} U\MakeLowercase{niversity.}}\footnote{C\MakeLowercase{orrespondence: mmohebuj@uab.edu}}%
		\and  Julian Miranda \footremember{tamiu}{D\MakeLowercase{epartment of} M\MakeLowercase{athematics and} P\MakeLowercase{hysics}, 
			T\MakeLowercase{exas} A\&M I\MakeLowercase{nternational} U\MakeLowercase{niversity}, TX 78041, USA; T\MakeLowercase{his author's work was partially supported by the} 
			N\MakeLowercase{ational} S\MakeLowercase{cience} F\MakeLowercase{oundation grant} DMS-2213274.}
		\and Md. Abdullah Al Mahbub\footremember{comilla}{D\MakeLowercase{epartment of} M\MakeLowercase{athematics}, C\MakeLowercase{omilla} 
			U\MakeLowercase{niversity}, C\MakeLowercase{umilla} 3506, B\MakeLowercase{angladesh;}}
		\and Mengying Xiao\footremember{UWF}{D\MakeLowercase{epartment of} M\MakeLowercase{athematics and} S\MakeLowercase{tatistics}, 
			U\MakeLowercase{niversity of} W\MakeLowercase{est} F\MakeLowercase{lorida}, P\MakeLowercase{ensacosa}, FL 32514, USA.}
	}
	
	\maketitle
	
	\begin{abstract}
		We propose, analyze, and test a penalty projection-based efficient and accurate algorithm for the Uncertainty Quantification (UQ) of the time-dependent Magnetohydrodynamic (MHD) flow problems in convection-dominated regimes. The algorithm uses the Els\"asser variables formulation and discrete Hodge decomposition to decouple the stochastic MHD system into four sub-problems (at each time-step for each realization) which are much easier to solve than solving the coupled saddle point problems. Each of the sub-problems is designed in a sophisticated way so that at each time-step the system matrix remains the same for all the realizations but with different right-hand-side vectors which allows saving a huge amount of computer memory and computational time. Moreover, the scheme is equipped with ensemble eddy-viscosity and grad-div stabilization terms. The stability of the algorithm is proven rigorously. We prove that the proposed scheme converges to an equivalent non-projection-based coupled MHD scheme for large grad-div stabilization parameter values. We examine how Stochastic Collocation Methods (SCMs) can be combined with the proposed penalty projection UQ algorithm. Finally, a series of numerical experiments are given which verify the predicted convergence rates, show the algorithm's performance on benchmark channel flow over a rectangular step, and a regularized lid-driven cavity problem with high random Reynolds number and magnetic Reynolds number.
	\end{abstract}

	
	
	\pagestyle{myheadings}
	\thispagestyle{plain}

	\markboth{\MakeUppercase{A penalty-projection efficient algorithm for stochastic MHD flows}}{\MakeUppercase{ M. Mohebujjaman, J. Miranda, M. A. 
			A. Mahbub, and M. Xiao}}
	
	\section{Introduction}
	Numerical simulation of MHD flow has been explored by many scientists \cite{B03,cibik2020analysis,D01,erkmen2020second,hossain2023mhd,LL60,li2018decoupled,MR17,wang2024convergence,yang2018efficient,zhang2018second} for the last couple of decades. However, their promise to reduce the computational cost and high accuracy for complex and larger MHD problems still remains an open question. The situation becomes even worse for the more realistic MHD flows, where convective-dominated flows interact with magnetic fields and model parameters involve random noises which introduce aleatoric uncertainty into the system and play a key role in determining the characteristic of the final solutions. We will use rigorous mathematics to develop novel computational frameworks to reduce the immense computational complexity involved in commonly used algorithms for Stochastic MHD (SMHD) flow problems.
	
	Let $\cD\subset \mathbb{R}^d\ (d=2,3)$ be a convex polygonal or polyhedral physical domain with boundary $\partial\cD$. A complete probability space is denoted by $(\Omega,\mathcal{F},P)$ with $\Omega$ the set of outcomes, $\mathcal{F}\subset 2^\Omega$ the $\sigma$-algebra of events, and $P:\mathcal{F}\rightarrow [0,1]$ represents a probability measure. We consider the time-dependent, dimensionless, viscoresistive, incompressible  SMHD flow problems for homogeneous Newtonian fluids which are governed by the following non-linear coupled stochastic PDEs \cite{AKMR15,B03,HMR17,LL60, Mohebujjaman2022High,MR17,mohebujjaman2022efficient,trenchea2014unconditional}:\vspace{-2mm}
	\begin{eqnarray}
		\bu_{t}+\bu\cdot\nabla \bu-s\bB\cdot\nabla \bB-\nabla\cdot\left(\nu(\bx,\omega)\nabla \bu\right)+\nabla p &= & \bif(t,\bx,\omega), 
		\hspace{2mm}\text{in}\hspace{2mm} (0,T]\times\cD,\label{momentum}\\
		\bB_{t}+\bu\cdot\nabla \bB-\bB\cdot\nabla \bu-\nabla\cdot\left(\nu_{m}(\bx,\omega)\nabla \bB\right)+\nabla\lambda &= & \nabla\times 
		\bg(t,\bx,\omega),\hspace{2mm}\text{in}\hspace{2mm} (0,T]\times\cD,\label{maxwell}\\
		\nabla\cdot \bu & =& 0, \hspace{2mm}\text{in}\hspace{2mm} (0,T]\times\cD, \\
		\nabla\cdot \bB &=& 0, \hspace{2mm}\text{in}\hspace{2mm} (0,T]\times\cD,\\ 
		\bu(\bx,0,\omega)& =& \bu^0(\bx),\hspace{2mm}\text{in}\hspace{2mm}\cD,\\
		\bB(\bx,0,\omega)& =& \bB^0(\bx),\hspace{2mm}\text{in}\hspace{2mm}\cD,\label{magnet-initial}
	\end{eqnarray}together with appropriate boundary conditions. Where $T>0$ represents the simulation end time and $\bx$ is the spatial variable. The viscosity $\nu(\bx,\omega)$ and magnetic diffusivity $\nu_m(\bx,\omega)$ are modeled as random fields with $\omega\in\Omega$. Here, the unknown quantities are the velocity field $\bu$ and the magnetic flux density $\bB$ which map as  $\bu,\bB:\Lambda\rightarrow\mathbb{R}^d$, and the modified pressure $p:\Lambda\rightarrow\mathbb{R}$, where $\Lambda:=(0,T]\times\cD\times\Omega$, and $d\in\{2,3\}$. The artificial magnetic pressure $\lambda:\Lambda\rightarrow\{0\}$, but $\lambda\neq 0$ in the discrete case. The external forces are represented by $\bif$, and $\nabla\times\bg$ in the momentum equation \eqref{momentum}, and induction equation \eqref{maxwell}, respectively. The coupling parameter is denoted by $s> 0$ which is the coefficient of the Lorentz force into the momentum equation \eqref{momentum}; If $s=0$, the fluid flow is not influenced by the magnetic field. The recent study shows \cite{aggul2023artificial, AKMR15, erkmen2020second,erkmen2019note,HMR17,li2018partitioned,Mohebujjaman2022High,MR17,mohebujjaman2022efficient,trenchea2014unconditional,wilson2015high} instead of solving \eqref{momentum}-\eqref{magnet-initial} (in terms of the original variables) together with appropriate boundary conditions, a change of variable called the Els\"asser variables formulation, allows to propose stable decoupled algorithms. This breakthrough idea was first presented by C. Trenchea in \cite{trenchea2014unconditional}.
	
	Defining $\bv:=\bu+\sqrt{s} \bB$, $\bw:=\bu-\sqrt{s}\bB$, $\bif_{1}:=\bif+\sqrt{s}\nabla\times \bg$, $\bif_{2}:=\bif-\sqrt{s}\nabla\times \bg$, 
	$q:=p+\sqrt{s}\lambda$, and $r:=p-\sqrt{s}\lambda$ produces the Els{\"{a}}sser variables formulation of the above stochastic MHD 
	system:\vspace{-2mm}
	\begin{align}
		\bv_{t}+\bw\cdot\nabla \bv-\nabla\cdot\left[\frac{\nu(\bx,\omega)+\nu_{m}(\bx,\omega)}{2}\nabla 
		\bv\right]-\nabla\cdot\left[\frac{\nu(\bx,\omega)-\nu_{m}(\bx,\omega)}{2}\nabla \bw\right]+\nabla q&=\bif_{1},\label{els1}\\
		\bw_{t}+\bv\cdot\nabla \bw-\nabla\cdot\left[\frac{\nu(\bx,\omega)+\nu_{m}(\bx,\omega)}{2}\nabla 
		\bw\right]-\nabla\cdot\left[\frac{\nu(\bx,\omega)-\nu_{m}(\bx,\omega)}{2}\nabla \bv\right]+\nabla r&=\bif_{2},\label{els2}\\
		\nabla\cdot \bv=\nabla\cdot \bw&=0,\label{els3}
	\end{align}
	together with the initial and boundary conditions. The $L^2(\cD)$ inner product is denoted by $(\cdot,\cdot)$. Defining the function spaces for 
	velocity and magnetic flux density as $\bX:=\bH_0^1(\cD)$, for the pressure and the magnetic pressure as $Q:=L_0^2(\cD)$, and the stochastic 
	space as $\bW:=\bL_P^2(\Omega)$, we get the weak formulation of \eqref{els1}-\eqref{els3} as\vspace{-4mm}
	\begin{align}
		\mathbb{E}[(\bv_{t},\bchi)]&+\mathbb{E}\left[\left(\bw\cdot\nabla 
		\bv,\bchi\right)\right]+\mathbb{E}\left[\frac{\nu(\bx,\omega)+\nu_{m}(\bx,\omega)}{2}\left(\nabla 
		\bv,\nabla\bchi\right)\right]\nonumber\\&+\mathbb{E}\left[\frac{\nu(\bx,\omega)-\nu_{m}(\bx,\omega)}{2}\left(\nabla 
		\bw,\nabla\bchi\right)\right]-\mathbb{E}\left( 
		q,\nabla\cdot\bchi\right)=\mathbb{E}\left[\left(\bif_{1},\bchi\right)\right],\hspace{3mm}\forall\bchi\in \bX\otimes \bW,\label{wels1}\\
		\mathbb{E}\left[\left(\bw_{t},\bl\right)\right]&+\mathbb{E}\left[(\bv\cdot\nabla 
		\bw,\bl)\right]+\mathbb{E}\left[\frac{\nu(\bx,\omega)+\nu_{m}(\bx,\omega)}{2}\left(\nabla 
		\bw,\nabla\bl\right)\right]\nonumber\\&+\mathbb{E}\left[\frac{\nu(\bx,\omega)-\nu_{m}(\bx,\omega)}{2}\left(\nabla 
		\bv,\nabla\bl\right)\right]-\mathbb{E}\left[\left( 
		r,\nabla\cdot\bl\right)\right]=\mathbb{E}\left[\left(\bif_{2},\bl\right)\right],\hspace{5mm}\forall\bl\in \bX\otimes \bW,\label{wels2}\\
		&\;\;\;\;\mathbb{E}[\left(\nabla\cdot \bv,\zeta\right)]=\mathbb{E}[\left(\nabla\cdot \bw,\eta\right)]=0,\hspace{31mm}\;\qquad\qquad\forall 
		\zeta,\eta\in Q\otimes \bW.\label{wels3}
	\end{align}
	In UQ, it is common to assume that the randomness is approximated by a finite number of random variables \cite{BNT07,gunzburger2019evolve}.
	To use SCMs, we consider $\by=(y_1,\cdots,y_N)\in\bGamma\subset\mathbb{R}^N$ be a finite $N\in\mathbb{N}$
	dimensional vector with joint probability density function $\rho(\by)$ in some parameter space $\bGamma=\prod\limits_{l=1}^N\Gamma_l$. Then, the 
	random fields $\nu(\bx,\omega)$, and $\nu_m(\bx,\omega)$ can be approximated in terms of the random variable as $\nu(\bx,\by)$, and 
	$\nu_m(\bx,\by)$, respectively. We define the space of square-integrable functions on $\bGamma$ subject to the weight $\rho(\by)$ as 
	$\bY:=\bL_\rho^2(\bGamma)$, and consider the following weak formulation:
	Find $\bv,\bw \in \bX \otimes \bY$ and $q,r \in Q \otimes \bY$ which, for almost all $t\in(0,T]$, satisfy
	\begin{align}
		\int_{\Gamma}  (\bv_{t},\bchi)\rho(\by)d\by &+ \int_{\Gamma} (\bw\cdot\nabla \bv,\bchi)\rho(\by)d\by + \int_{\Gamma}  
		\left[\frac{\nu(\bx,\by)+\nu_{m}(\bx,\by)}{2}(\nabla \bv,\nabla \bchi)\right]\rho(\by)d\by \nonumber\\&+\int_{\Gamma} 
		\left[\frac{\nu(\bx,\by)-\nu_{m}(\bx,\by)}{2}\left(\nabla \bw,\nabla\bchi\right)\right]\rho(\by)d\by  -\int_{\Gamma}  (q ,\nabla\cdot 
		\bchi)\rho(\by)d\by 
		\nonumber\\ &= \int_{\Gamma}  (\bif_1,\bchi)\rho(\by)d\by, \qquad\qquad\qquad\qquad\qquad\qquad\forall \bchi\in \bX \otimes \bY,\label{weak_formulation_final-start}\\
		\int_{\Gamma}  \left(\bw_{t},\bl\right)\rho(\by)d\by&+\int_{\Gamma} (\bv\cdot\nabla \bw,\bl)\rho(\by)d\by+\int_{\Gamma} 
		\left[\frac{\nu(\bx,\by)+\nu_{m}(\bx,\by)}{2}\left(\nabla \bw,\nabla\bl\right)\right]\rho(\by)d\by\nonumber\\&+\int_{\Gamma} 
		\left[\frac{\nu(\bx,\by)-\nu_{m}(\bx,\by)}{2}\left(\nabla \bv,\nabla\bl\right)\right]\rho(\by)d\by-\int_{\Gamma} \left( 
		r,\nabla\cdot\bl\right)\rho(\by)d\by\nonumber\\&=\int_{\Gamma} 
		\left(\bif_{2},\bl\right)\rho(\by)d\by,\qquad\qquad\qquad\qquad\qquad\qquad\forall\bl\in \bX\otimes \bY,\\
		\int_{\Gamma}  (\nabla\cdot \bv,\zeta) \rho(\by)d\by&=\int_{\Gamma}  (\nabla\cdot \bw,\eta) \rho(\by)d\by    =0, \qquad \qquad \qquad 
		\qquad \;\;\forall \zeta,\eta\in Q \otimes \bY.\label{weak_formulation_final-end}
	\end{align}
	To have an effective and efficient penalty-projection scheme for UQ, we assume affine dependence of the random variables for the viscosity and 
	magnetic diffusivity as below:\vspace{-4mm}
	\begin{align*}
		\nu(\bx,\by)=\nu_0(\bx)+\sum_{l=1}^N\nu_l(\bx)y_l,~~\text{and}~~\nu_m(\bx,\by)=\nu_{m,0}(\bx)+\sum_{l=1}^N\nu_{m,l}(\bx)y_l.
	\end{align*}
	In order to have a robust high fidelity solution of \eqref{els1}-\eqref{els3}, which is often essential for many surrogate models \cite{chen2004subpattern,mohebujjaman2019physically,proctor2016dynamic,XMRI17}, one of the major hurdles is to realize the model over an ensemble 
	of flow parameters with high spatial resolution. In this case, a popular approach is to use SCMs 
	\cite{BNT07,CCMFT2013,gunzburger2019evolve,gunzburger2014stochastic,lee2018stochastic,FTW2008,stoyanov2013hierarchy,stoyanov2016dynamically}, which requires fewer realizations compare to other sampling-based UQ methods, such as the Monte Carlo method which comes with high computational complexity \cite{johnson2021coupled}. SCMs use global polynomial approximation and are independent of the PDE solvers, thus compatible with combining with any legacy code. 
	
	In this paper, we present, analyze, and test, a novel, efficient, and accurate SCMs based Stabilized Penalty Projection algorithm for the UQ of SMHD (SCM-SPP-SMHD) flow problems, which is the conjunction of the SCMs and an SPP - Finite Element Method (SPP-FEM). The SPP-FEM is based on 
	Els{\"{a}}sser formulation which provides a stable decoupling of the SMHD system into two Oseen type sub-problems (velocity-pressure and magnetic 
	field-magnetic pressure types saddle-point sub-problems). These two sub-problems can be solved simultaneously if the resources are available. For 
	each of these sub-problems, we employ ``Projection Methods'' which utilize a discrete Hodge decomposition at each time-step (which was proposed by Chorin and Temam \cite{chorin1968numerical,temam1969approximation} in the early 1960s) together with recent stabilization techniques \cite{AKMR15,linke2017connection}. This allows us to solve the difficult $2\times 2$ block saddle-point sub-problems into two easier linear solves (a $1\times 1$, and a $2\times 2$ block systems), particularly in 3D problems with high Reynolds number and high magnetic Reynolds numbers. The use of a large grad-div stabilization parameter helps us to improve the accuracy of the penalty-projection algorithm which has examined on problems of Navier-Stokes (N-S) flow \cite{linke2017connection}, and in fluid-fluid interaction \cite{aggul2022grad}. Also, each of the sub-problems contains an ensemble eddy viscosity term which is taken from the idea of turbulence modeling techniques to reduce the numerical 
	instability particularly in 3D problems \cite{MR17,mohebujjaman2022efficient}.
	
	Moreover, for each of the linear solves, the scheme is designed in an elegant way so that at each time-step, each realization shares a common 
	system matrix but a different right-hand-side vector, following the breakthrough idea by N. Jiang and W. Layton given in \cite{JL14}. This allows us to use orders of magnitude shorter than 
	the system matrix assembly time (which is assumed to be the most time-consuming step in the finite element assembly process). Moreover, for 
	problems for which a direct solver is appropriate, the LU decomposition or its variants of the system matrix is needed to compute only once per 
	time-step. For large size and complex problems, Krylov subspace methods are appropriate \cite{mohebujjaman2023scalability} for which a single 
	preconditioner is needed to be built for each sub-problem per time-step. Further, the advantage of the block linear solvers can be taken with a 
	single system matrix for multiple right-hand-side vectors at each time-step. This elegant feature leads to saving a huge computational cost and 
	memory for the UQ of complex dynamical systems and is successfully implemented in the surface data assimilation \cite{fujita2007surface}, 
	turbulence modeling \cite{JKL15}, porous media flow \cite{jiang2021artificial}, Boussinesq \cite{jiang2023unconditionally}, weather forecasting 
	\cite{L05,LP08}, spectral methods \cite{LK10}, sensitivity analyses \cite{MX06}, MHD \cite{Mohebujjaman2022High,MR17,mohebujjaman2022efficient},   
	N-S simulations \cite{J15,J16,JL15,NTW16},   and hydrology \cite{GG11}.
	
	The proposed SCM-SPP-SMHD algorithm consists of grad-div stabilization terms with coefficient parameter $\gamma$. Large $\gamma$ values help to 
	achieve optimal temporal accuracy, reducing penalty-projection splitting errors \cite{AKMR15,linke2017connection}. Finally, using straightforward 
	transformations, we get back the solution in terms of the original variables. Thus, we consider a uniform time-step size $\Delta t$ and let 
	$t_n=n\Delta t$ for $n=0, 1, \cdots$., (suppress the spatial discretization momentarily), then computing the $N_{sc}$ (number of stochastic 
	collocation points) solutions independently, takes the following form: For $j$=1,2,...,$N_{sc}$, \\ 
	Step 1: Compute $\bhv_j^{n+1}$:
	\begin{align}
		\frac{\bhv_j^{n+1}}{\Delta t}+&<\bhw>^n\cdot\nabla \bhv_j^{n+1}-\nabla\cdot\left(\frac{\Bar{\nu}+\Bar{\nu}_{m}}{2}\nabla 
		\bhv_j^{n+1}\right)-\nabla\left(\gamma\nabla\cdot\bhv_{j}^{n+1}\right)-\nabla\cdot\left(2\nu_T(\bhw^{'},t^n)\nabla \bhv_j^{n+1}\right)  
		\nonumber\\&= \bif_{1,j}(t^{n+1})+\frac{\btv_j^n}{\Delta t}-\bhw_j^{'n}\cdot\nabla \bhv_j^n+\nabla\cdot\left(\frac{\nu_j-\nu_{m,j}}{2}\nabla 
		\bhw_j^n\right)+\nabla\cdot\left(\frac{\nu_j^{'}+\nu_{m,j}^{'}}{2}\nabla \bhv_j^{n}\right).\label{scheme1}
	\end{align}
	Step 2: Compute $\btv_j^{n+1}$, and $\hq_j^{n+1}$:
	\begin{align} \frac{\btv_{j}^{n+1}}{\Delta t}+\nabla\hq_{j}^{n+1} =&\frac{\bhv_{j}^{n+1}}{\Delta 
			t},\label{weaknewnew2}\\\nabla\cdot\btv_{j}^{n+1}=&0.\label{weaknewnew3}
	\end{align}
	Step 3: Compute $\bhw_j^{n+1}$:
	\begin{align}
		\frac{\bhw_j^{n+1}}{\Delta t}+&<\bhv>^n\cdot\nabla \bhw_j^{n+1}-\nabla\cdot\left(\frac{\Bar{\nu}+\Bar{\nu}_{m}}{2}\nabla 
		\bhw_j^{n+1}\right)-\nabla\left(\gamma\nabla\cdot\bhw_{j}^{n+1}\right)-\nabla\cdot\left(2\nu_T(\bhv^{'},t^n)\nabla \bhw_j^{n+1}\right)   
		\nonumber\\&= \bif_{2,j}(t^{n+1})+\frac{\btw_j^n}{\Delta t}-\bhv_j^{'n}\cdot\nabla \bhw_j^n+\nabla\cdot\left(\frac{\nu_j-\nu_{m,j}}{2}\nabla 	\bhv_j^n\right)+\nabla\cdot\left(\frac{\nu_j^{'}+\nu_{m,j}^{'}}{2}\nabla \bhw_j^{n}\right).\label{scheme2}
	\end{align}
	Step 4: Compute $\btw_j^{n+1}$, and $\hr_j^{n+1}$:
	\begin{align} \frac{\btw_{j}^{n+1}}{\Delta t}+\nabla \hr_{j}^{n+1} =&\frac{\bhw_{j}^{n+1}}{\Delta 
			t},\label{weaknewnew5}\\\nabla\cdot\btw_{j}^{n+1}=&0.\label{weaknewnew4}
	\end{align}
	The ensemble mean and fluctuation about the mean are defined as follows:
	\begin{align*}
		<\bhu>^n:&=\frac{1}{N_{sc}}\sum\limits_{j=1}^{N_{sc}}\bhu_j^n, \hspace{2mm} \bhu_j^{'n}:=\bhu_j^n-<\bhu>^n,\\ 
		\Bar{\nu}:=\frac{1}{N_{sc}}\sum\limits_{j=1}^{N_{sc}}\nu_j,\hspace{2mm}\nu_j^{'}:&=\nu_j-\Bar{\nu},\hspace{2mm} 
		\Bar{\nu}_m:=\frac{1}{N_{sc}}\sum\limits_{j=1}^{N_{sc}}\nu_{m,j},\hspace{2mm}\nu_{m,j}^{'}:=\nu_{m,j}-\Bar{\nu}_m.
	\end{align*}
	The eddy viscosity term, which is of $O(\Delta t)$ accurate, is defined using mixing length phenomenology following \cite{JL15}, and is given by 
	\begin{align}
		\nu_T(\bz^{'},t^n):=\mu\Delta t(l_{\bz}^{n})^2,\;\text{and}\;(l_{\bz}^{n})^2:=\sum\limits_{j=1}^{N_{sc}}|\bz_j^{'n}|^2,\label{eddy-viscosity-term}
	\end{align}where $|\cdot|$ denotes length of a vector. The eddy viscosity term helps the scheme to provide stability for flows that are not resolved on particular meshes. The solutions $\bhv$, and $\bhw$ in Step 1, and Step 3, respectively, do not satisfy the divergence-free conditions, whereas the solution $\btv$, and $\btw$ in Step 2, and Step 3, respectively, do not satisfy the boundary conditions. Step 1 has only unknown $\bhv_j^{n+1}$, and (the finite element variational formation) provides a $1\times 1$ block linear system which is not dependent on the index $j$, thus for all $N_{sc}$ realizations, the system matrix remains the same but the right-hand-side vector varies. Therefore, at each time-step, we need to solve a linear system of equations of the form $A[\bx_1|\bx_2|\cdots|\bx_{N_{sc}}]=[\bb_1|\bb_2|\cdots|\bb_{N_{sc}}],$ where $A$ is a sparse coefficient matrix, $\bx_j$, and $\bb_j$ are the solution, and right-hand-side vector for the $j$-th realization, respectively. This feature prevails in all other three steps which makes the algorithm efficient in saving a huge time in global system matrix assembly and in saving a massive computer memory. The size of the system matrix $A$ is much smaller than the size of the corresponding system matrix (which is the $2\times 2$ block linear system) that arises in the saddle-point sub-problems. Step 2 requires a linear solve for its $2\times 2$ block system, which is a symmetric positive definite matrix (since there is no non-linear term present), thus the advantage of using block conjugate-gradient method \cite{o1980block} can be taken. Therefore, by employing the penalty-projection splitting, we replace the difficult linear solve for each of the saddle point sub-systems (corresponding to the Oseen-type sub-problems that arise after decoupling the system given in \eqref{els1}-\eqref{els3}, see \cite{mohebujjaman2022efficient}) into two easier linear solves at each time-step. In Step 1, and 
	Step 3, the system matrices are nonsymmetric (due to the presence of non-linear terms), and we can take advantage of the block GMRES \cite{simoncini1996hybrid} solver for the nonsymmetric system with multiple right-hand-side vectors. On the other hand, for problems in which a direct solver is more appropriate, its decomposition needs to be done only once per time-step and can be reused for all the realizations.
	
	Therefore, at each time-step, the feature of having a common system matrix in each of Steps 1-4, and the penalty-projection splitting make the 
	algorithm efficient in saving a huge computational time and computer memory. Moreover, having a grad-div stabilization term in Step 1, and in Step 3, helps the scheme in achieving temporal accuracy similar to the non-splitting (velocity-pressure coupled type) algorithms for the coefficient $\gamma\rightarrow\infty$.
	
	Using finite element spatial discretization, in this paper, we investigate the novel SCM-SPP-SMHD ensemble scheme in a fully discrete setting. The efficient SCM-SPP-SMHD scheme is proved to be stable and convergent without any time-step restriction but takes care of uncertainties in all model data. To the best of our knowledge, SCM-SPP-SMHD is new for the UQ of MHD flow problems.
	
	The rest of the paper is organized as follows: We provide necessary notations and mathematical preliminaries in Section \ref{notation-prelims} to 
	follow a smooth analysis. As a benchmark algorithm for \eqref{els1}-\eqref{els3}, we consider a first-order backward-Euler time-stepping fully discrete Coupled (where the velocity and magnetic fields like variables are decoupled but the velocity with the pressure, and magnetic field with the magnetic pressure are not) algorithm given in \cite{mohebujjaman2022efficient} for Stochastic MHD (Coupled-SMHD) in Section \ref{fully-discrete-scheme}. We also discuss the stability, convergence, regularity assumptions, and small data assumptions of the Coupled-SMHD scheme for a fair comparison with the SCM-SPP-SMHD scheme. In Section \ref{SPP-FEM-scheme}, we present the SCM-SPP-SMHD scheme and describe the additional functional space we need for further analysis. We also state and prove the stability and convergence theorem of the SCM-SPP-SMHD scheme in Section \ref{SPP-FEM-scheme}. A brief description of SCMs is given in Section \ref{scm}. To support the theoretical analysis, we compute the convergence rates varying $\gamma$, time-step size, and mesh width, and finally implement the scheme in benchmark channel flow past a rectangular step and a regularized lid-driven cavity problems with space-dependent variable high random Reynolds number and variable high random magnetic Reynolds number in Section \ref{numerical-exp}. Finally, conclusions and future research avenues are given in Section \ref{conclusion-future-works}.

	\section{\large Notation and preliminaries}\label{notation-prelims}
	
	The usual $L^2(\cD)$ norm and inner product are denoted by $\|.\|$ and $(.,.)$, respectively. Similarly, the $L^p(\cD)$ norms and the Sobolev 
	$W_p^k(\cD)$ norms are $\|.\|_{L^p}$ and $\|.\|_{W_p^k}$, respectively for $k\in\mathbb{N},\hspace{1mm}1\le p\le \infty$. Sobolev space 
	$W_2^k(\cD)$ is represented by $H^k(\cD)$ with norm $\|.\|_k$. The vector-valued spaces are $$\bL^p(\cD)=(L^p(\cD))^d, 
	\hspace{1mm}\text{and}\hspace{1mm}\bH^k(\cD)=(H^k(\cD))^d.$$
	For $\bX$ being a normed function space in $\cD$, $L^p(0,T;\bX)$ is the space of all functions defined on $(0,T]\times\cD$ for which the 
	following norm
	\begin{align*}
		\|\bu\|_{L^p(0,T;\bX)}=\lp\int_0^T\|\bu\|_{\bX}^pdt\rp^\frac{1}{p},\hspace{2mm}p\in[1,\infty)
	\end{align*}
	is finite. For $p=\infty$, the usual modification is used in the definition of this space. The natural function spaces for our problem are
	\begin{align*}
		\bX:&=\bH_0^1(\cD)=\{\bv\in \bL^p(\cD) :\nabla \bv\in L^2(\cD)^{d\times d}, \bv=0 \hspace{2mm} \mbox{on}\hspace{2mm}   \partial \cD\},\\
		\tilde{\bY}:&=\{\bv\in\bH^1(\cD):\bv\cdot\bnh=0\;\;\text{on}\;\partial\cD\},\\
		Q:&=L_0^2(\cD)=\{ q\in L^2(\cD): \int_\cD q\hspace{1mm}d\bx=0\},
	\end{align*}
	where $\bnh$ denotes the outward unit normal vector normal to the boundary $\partial\cD$. Recall the Poincar\'e inequality holds in $X$: There 
	exists $C$ depending only on $\cD$ satisfying for all $\bphi\in X$,
	\[
	\| \bphi \| \le C \| \nabla \bphi \|.
	\]
	The divergence free velocity space is given by
	$$\bV:=\{\bv\in \bX:(\nabla\cdot \bv, q)=0, \forall q\in Q\}.$$
	We define the trilinear form $b:\bX\times \bX\times \bX\rightarrow \mathbb{R}$ by
	\[
	b(\bu,\bv,\bw):=(\bu\cdot\nabla \bv,\bw), 
	\]
	and recall from \cite{GR86} that $b(\bu,\bv,\bv)=0$ if $\bu\in \bV$, and 
	\begin{align}
		|b(\bu,\bv,\bw)|\leq C(\cD)\|\nabla \bu\|\|\nabla \bv\|\|\nabla \bw\|,\hspace{2mm}\mbox{for any}\hspace{2mm}\bu,\bv,\bw\in 
		\bX.\label{nonlinearbound}
	\end{align}
	
	The conforming finite element spaces are denoted by $\bX_h\subset \bX$, and  $Q_h\subset Q$, and we assume a regular triangulation $\tau_h(\cD)$, 
	where $h$ is the maximum triangle diameter. We assume that $(\bX_h,Q_h)$ satisfies the usual discrete inf-sup condition
	\begin{eqnarray}
		\inf_{q_h\in Q_h}\sup_{\bv_h\in \bX_h}\frac{(q_h,\grad\cdot \bv_h)}{\|q_h\|\|\grad \bv_h\|}\geq\beta>0,\label{infsup}
	\end{eqnarray}
	where $\beta$ is independent of $h$. We assume that there exists a finite element space  $\tilde{\bY}_h\subset\tilde{\bY}$. The space of 
	discretely divergence free functions is defined as 
	\begin{align*}
		\bV_h:=\{\bv_h\in \bX_h:(\nabla\cdot \bv_h,q_h)=0,\hspace{2mm}\forall q_h\in Q_h\}.
	\end{align*}
	For simplicity of our analysis, we will use Scott-Vogelius (SV) finite element pair $(\bX_h, Q_h)=((P_k)^d, P_{k-1}^{disc})$,  which satisfies 
	the \textit{inf-sup} condition when the mesh is created as a barycenter refinement of a regular mesh, and the polynomial degree $k\ge d$  
	\cite{arnold1992quadratic,Z05}. Our analysis can be extended without difficulty to any inf-sup stable element choice, {\color{black} however, 
		there will be additional terms that appear in the convergence analysis if non-divergence-free elements are chosen.}

	We have the following approximation properties in $(\bX_h,Q_h)$: \cite{BS08}
	\begin{align}
		\inf_{\bv_h\in \bX_h}\|\bu-\bv_h\|&\leq Ch^{k+1}|\bu|_{k+1},\hspace{2mm}\bu\in \bH^{k+1}(\cD),\label{AppPro1}\\
		\inf_{\bv_h\in \bX_h}\|\grad (\bu-\bv_h)\|&\leq Ch^{k}|\bu|_{k+1},\hspace{5mm}\bu\in \bH^{k+1}(\cD),\label{AppPro2}\\
		\inf_{q_h\in Q_h}\|p-q_h\|&\leq Ch^k|p|_k,\hspace{10mm}p\in H^k(\cD),
	\end{align}
	where $|\cdot|_r$ denotes the $H^r$ or $\bH^r$ seminorm.
	
	We will assume the mesh is sufficiently regular for the inverse inequality to hold.The following lemma for the discrete Gr\"onwall inequality was 
	given in \cite{HR90}.
	\begin{lemma}\label{dgl}
		Let $\Delta t$, $\mathcal{E}$, $a_n$, $b_n$, $c_n$, $d_n$ be non-negative numbers for $n=1,\cdots, M$ such that
		$$a_M+\Delta t \sum_{n=1}^Mb_n\leq \Delta t\sum_{n=1}^{M-1}{d_na_n}+\Delta 
		t\sum_{n=1}^Mc_n+\mathcal{E}\hspace{3mm}\mbox{for}\hspace{2mm}M\in\mathbb{N},$$
		then for all $\Delta t> 0,$
		$$a_M+\Delta t\sum_{n=1}^Mb_n\leq \mbox{exp}\left(\Delta t\sum_{n=1}^{M-1}d_n\right)\lp\Delta 
		t\sum_{n=1}^Mc_n+\mathcal{E}\rp\hspace{2mm}\mbox{for}\hspace{2mm}M\in\mathbb{N}.$$
	\end{lemma}
	
	\section{\large Coupled backward-Euler fully discrete time-stepping scheme}\label{fully-discrete-scheme}
	In this section, we consider a first-order accurate, fully discrete, backward-Euler time-stepping ``Coupled-SMHD'' algorithm as the benchmark UQ scheme for SMHD flow problems which was proposed and analyzed in \cite{mohebujjaman2022efficient}. We state its stability and convergence theorems and give a rigorous proof of the adopted small data assumption which will be used in the analysis of the proposed SCM-SPP-SMHD scheme in Section \ref{SPP-FEM-scheme}. The Coupled-SMHD scheme is based on Els{\"{a}}sser variables formulation, although, its momentum and induction-like equations are decoupled, the velocity and pressure, magnetic field, and magnetic pressure-like variables are still coupled. The unconditionally stable first-order temporally accurate Coupled-SMHD scheme is proven to be optimally accurate in 2D, and sub-optimally accurate in 3D, which is because of the ensemble eddy-viscosity term present in the scheme that leads to the use of the discrete inverse inequality. The Coupled-SMHD scheme was successfully implemented into a benchmark regularized lid-driven cavity and flow over rectangular step problems with high random Reynolds number with random low magnetic diffusivity parameter. The Coupled-SMHD algorithm proposed in \cite{mohebujjaman2022efficient} is efficient in terms of saving a huge computational time and computer memory as it allows to use of block linear solver at each time-step and is able to provide solutions of all realizations using a single block linear solve. However, for the Coupled-SMHD, we are required to solve the saddle-point system at each time since the velocity- and pressure-like variables are coupled together. To decouple the velocity- and pressure-like variables, and to make the Coupled-SMHD scheme more efficient, we combine the stochastic collocation method with a recently proposed grad-div stabilization penalty projection method for N-S flow problems in \cite{linke2017connection}. The Coupled-SMHD is presented in Algorithm \ref{Algn0}.
	\begin{algorithm}[H]\label{Algn0}
		\caption{Coupled-SMHD ensemble scheme \cite{mohebujjaman2022efficient}} Given time-step $\Delta t>0$, end time $T>0$, initial conditions 
		$\bv_j^0, \bw_j^0\in \bV_h$ and $\bif_{1,j}, \bif_{2,j}\in$ $ L^\infty\left( 0,T;\bH^{-1}(\cD)\right)$ for $j=1,2,\cdots\hspace{-0.35mm},N_{sc}$. 
		Set $M=T/\Delta t$ and for $n=1,2,\cdots\hspace{-0.35mm},M-1$, compute:
		Find $(\bv_{j,h}^{n+1},q_{j,h}^{n+1})\in \bX_h\times Q_h$ satisfying, for all $(\bchi_{h},\zeta_h)\in \bX_h\times Q_h$:\vspace{-2mm}
		\begin{align}
			\Bigg(&\frac{\bv_{j,h}^{n+1}-\bv_{j,h}^n}{\Delta t},\; \bchi_{h}\Bigg)+b\left(<\bw_h>^n, 
			\bv_{j,h}^{n+1},\bchi_{h}\right)+\left(\frac{\Bar{\nu}+\Bar{\nu}_{m}}{2}\nabla \bv_{j,h}^{n+1},\nabla 
			\bchi_{h}\right)\nonumber\\&+\left(2\nu_T(\bw^{'}_{h},t^n)\nabla \bv_{j,h}^{n+1},\nabla \bchi_{h}\right)-(q_{j,h}^{n+1},\nabla\cdot\bchi_{h})= 
			\left(\bif_{1,j}(t^{n+1}),\bchi_{h}\right)-b\left(\bw_{j,h}^{'n}, \bv_{j,h}^n,\bchi_{h}\right)\nonumber\\&-\left(\frac{\nu_j-\nu_{m,j}}{2}\nabla 
			\bw_{j,h}^n,\nabla \bchi_{h}\right)-\left(\frac{\nu_j^{'}+\nu_{m,j}^{'}}{2}\nabla \bv_{j,h}^{n},\nabla\bchi_{h}\right),\label{weak1}\\
			&(\nabla\cdot\bv_{j,h}^{n+1},\zeta_h)=0.
		\end{align}
		Find $(\bw_{j,h}^{n+1},r_{j,h}^{n+1})\in \bX_h\times Q_h$ satisfying, for all $(\bl_{h},\eta_h)\in \bX_h\times Q_h$:\vspace{-2mm}
		\begin{align}
			\Bigg(&\frac{\bw_{j,h}^{n+1}-\bw_{j,h}^n}{\Delta t},\bl_{h}\Bigg)+b\left(<\bv_h>^n, 
			\bw_{j,h}^{n+1},\bl_{h}\right)+\left(\frac{\Bar{\nu}+\Bar{\nu}_{m}}{2}\nabla \bw_{j,h}^{n+1},\nabla 
			\bl_{h}\right)\nonumber\\+&\left(2\nu_T(\bv^{'}_{h},t^n)\nabla \bw_{j,h}^{n+1},\nabla \bl_{h}\right)-(r_{j,h}^{n+1},\nabla\cdot\bchi_{h})= 
			\left(\bif_{2,j}(t^{n+1}),\bl_{h}\right)-b\left(\bv_{j,h}^{'n}, \bw_{j,h}^n,\bl_{h}\right)\nonumber\\-&\left(\frac{\nu_j-\nu_{m,j}}{2}\nabla 
			\bv_{j,h}^n,\nabla \bl_{h}\right)-\left(\frac{\nu_j^{'}+\nu_{m,j}^{'}}{2}\nabla \bw_{j,h}^{n},\nabla \bl_{h}\right),\label{weak2}\\
			&(\nabla\cdot\bw_{j,h}^{n+1},\eta_h)=0.\label{weak*}
		\end{align}
	\end{algorithm} 
	\noindent For simplicity of our analysis, we define $\Bar{\nu}_{\min}:=\min\limits_{\bx\in\cD}\Bar{\nu}(\bx)$,  $\Bar{\nu}_{m,\min}:=\min\limits_{\bx\in\cD}\Bar{\nu}_m(\bx)$, and $\alpha_{\min}:=\min\limits_{1\le j\le J}\alpha_j$,  where 
	$\alpha_j:=\Bar{\nu}_{\min}+\Bar{\nu}_{m,\min}-\|\nu_j-\nu_{m,j}\|_\infty-\|\nu_j^{'}+\nu_{m,j}^{'}\|_\infty$, for $j=1,2,\cdots\hspace{-0.35mm}, 
	N_{sc}$, and 
	state the stability and convergence theorems of the Algorithm \ref{Algn0}.
	
	\begin{theorem} Suppose $\textbf{f}_{1,j},\textbf{f}_{2,j}\in L^2\left(0,T;\bH^{-1}(\cD)\right)$, and $\bv_{j,h}^0$, $\bw_{j,h}^0 \in 
		\bH^1(\cD)$, for $j=1,2,\cdots\hspace{-0.35mm}, N_{sc}$ then the solutions to the Algorithm \ref{Algn0} are stable: For any $\Delta t>0$, if 
		$\alpha_j>0$, and $\mu>\frac{1}{2}$
		\begin{align}   \|\bv_{j,h}^M\|^2+\|\bw_{j,h}^M\|^2+\frac{\alpha_{\min}\Delta t}{2}\sum_{n=1}^{M}\Big(\|\nabla \bv_{j,h}^n\|^2+\|\nabla 
			\bw_{j,h}^n\|^2\Big)\le C(data).
		\end{align}\label{stability-theorem}
	\end{theorem}
	
	\begin{proof}
		See Theorem 4.1 in  \cite{mohebujjaman2022efficient}.
	\end{proof}
	\begin{theorem}
		Assume $\left(\bv_j, \bw_j, q_j, r_j\right)$ satisfying \eqref{els1}-\eqref{els3} with regularity assumptions $\bv_j, \bw_j\in 
		L^\infty(0,T;\bH^{k+1}(\cD))$, $\bv_{j,t}, \bw_{j,t}\in L^\infty(0,T;\bH^{2}(\cD))$, $\bv_{j,tt}, \bw_{j,tt}\in L^\infty(0,T;\bL^{2}(\cD))$ for 
		$k\ge 2$ and  $j=1,2,\cdots\hspace{-0.35mm},N_{sc}$, then the solution $(\bv_{j,h},\bw_{j,h})$ to the Algorithm \ref{Algn0} converges to the true 
		solution: For any $\Delta t>0$, if $\alpha_j>0$, and $\mu>\frac12$, then one has
		\begin{align}
			\|\bv_j(T)&-\bv_{j,h}^M\|^2+\|\bw_j(T)-\bw_{j,h}^M\|^2+\frac{\alpha_{\min}\Delta t}{2}\sum_{n=1}^{M}\Big(\|\nabla 
			\left(\bv_j(t^n)-\bv_{j,h}^n\right)\|^2+\|\nabla \left(\bw_j(t^n)-\bw_{j,h}^n\right)\|^2\Big)\nonumber\\&\le  C\big(\Delta 
			t^2+h^{2k}+h^{2k}\Delta t^2+h^{2-d}\Delta t^2+ h^{2k-1}\Delta t\big).\label{error-bound-regular-alg}
		\end{align}
	\end{theorem}
	\begin{proof}
		See the Theorem 5.1, equation (5.18) in \cite{mohebujjaman2022efficient}.
	\end{proof}

	\begin{lemma}\label{lemma1}
		Assume the true solution $\bv_j,\bw_j\in L^\infty(0,T;\bH^2(\cD))$. Then, there exists a constant $C_*$ which is independent of $h$, $\Delta t$, 
		and $\gamma$ such that for sufficiently small $h$ and $\Delta t$, the solutions of the Algorithm \ref{Algn0} satisfies
		\begin{align*}
			\max_{1\le n\le M}\Big(\|\nabla \bv_{j,h}^n\|_{L^3}+\|\nabla \bw_{j,h}^n\|_{L^3}+\|\bv_{j,h}^n\|_\infty+\|\bw_{j,h}^n\|_\infty\Big)&\le 
			C_*,\hspace{2mm}\forall j=1,2,\cdots,N_{sc}.
		\end{align*}
	\end{lemma}
	
	\begin{proof}
		Using triangle inequality, we write
		\begin{align}
			\|\nabla \bv_{j,h}^n\|_{L^3}+\|\bv_{j,h}^n\|_\infty&\le 
			\|\nabla(\bv_{j,h}^n-\bv_j(t^n))\|_{L^3}+\|\nabla\bv_j(t^n)\|_{L^3}+\|\bv_{j,h}^n-\bv_j(t^n)\|_\infty+\|\bv_j(t^n)\|_\infty.\label{lemma-trianlge}
		\end{align}
		Apply Sobolev embedding theorem on the first two terms, and Agmon’s \cite{Robinson2016Three-Dimensional} inequality on the last two terms in the 
		right-hand-side of \eqref{lemma-trianlge}, to provide
		\begin{align}
			\|\nabla \bv_{j,h}^n\|_{L^3}+\|\bv_{j,h}^n\|_\infty\le 
			C\|\nabla(\bv_{j,h}^n-\bv_j(t^n))\|^{\frac12}\|\nabla^2(\bv_{j,h}^n-\bv_j(t^n))\|^{\frac12}+C\|\bv_j(t^n)\|_{H^1}^\frac12\|\bv_j(t^n)\|_{H^2}^\frac12.
		\end{align}
		Apply the regularity assumption of the true solution and
		discrete inverse inequality, to obtain
		\begin{align}
			\|\nabla \bv_{j,h}^n\|_{L^3}+\|\bv_{j,h}^n\|_\infty&\le Ch^{-\frac12}\|\nabla(\bv_{j,h}^n-\bv_j(t^n))\|+C.
		\end{align}
		Consider the $(P_k,P_{k-1})$ element for the pair $(\bv_{j,h},q_{j,h})$, and use the error bounds in \eqref{error-bound-regular-alg}, gives 
		\begin{align*}
			\|\nabla \bv_{j,h}^n\|_{L^3}+\|\bv_{j,h}^n\|_\infty&\le Ch^{-\frac12}\left(\Delta t^\frac12+\frac{h^k}{\Delta t^\frac{1}{2}}+h^k\Delta 
			t^\frac12+h^{1-\frac{d}{2}}\Delta t^\frac12+h^{k-\frac12}\right)+C.
		\end{align*}
		Choose $\Delta t$ so that
		\begin{align*}
			\frac{\Delta t^\frac12}{h^\frac12}\le \frac{1}{C},\;
			\frac{h^{k-\frac12}}{\Delta t^\frac12}\le\frac{1}{C},\;
			h^{k-\frac12}\Delta t^\frac12\le\frac{1}{C},\; h^{\frac{1-d}{2}}\Delta t^\frac12\le \frac{1}{C}\;\text{and}\;
			h^{k-1}\le\frac{1}{C},
		\end{align*}
		which gives
		\begin{align*}
			\|\nabla \bv_{j,h}^n\|_{L^3}+\|\bv_{j,h}^n\|_\infty\le 4+C,
		\end{align*} with the time-step restrictions  $O(h^{2k-1})\le\Delta t\le O(h^{d-1})$. Similarly, we can show \begin{align*}
			\|\nabla \bw_{j,h}^n\|_{L^3}+\|\bw_{j,h}^n\|_\infty\le 4+C.
		\end{align*}
		Therefore, $C_*:=8+C$ completes the proof.
	\end{proof}
	
	\section{Efficient SCM-SPP-SMHD scheme for UQ of SMHD flow problems}\label{SPP-FEM-scheme}
	In this section, we present the proposed efficient, fully discrete, and stable penalty-projection-based decoupled time-stepping algorithm
	SCM-SPP-SMHD,  which combines the SCM for the UQ of SMHD flow problems. We state and prove the unconditional stability theorem and provide an error analysis which shows that as  $\gamma\rightarrow\infty$, the outcomes of the SCM-SPP-SMHD scheme converge to the outcomes of the Coupled-SMHD in Algorithm \ref{Algn0}. 
	
	The SCM-SPP-SMHD computes the solution in four steps. The scheme is designed in a technical way so that at each of these steps, the system matrix remains the same for all realizations, which makes it efficient in saving a huge computational time and computer memory. We describe the  scheme below:
	\begin{algorithm}[H]\label{SPP-FEM}
		\caption{Fully discrete, decoupled grad-div Stabilized Penalty-projection Finite Element Method (SCM-SPP-SMHD)}
		Given time-step $\Delta t>0$, end time $T>0$, stabilization parameter $\gamma>0$, initial conditions $\bv_j^0, \bw_j^0\in \bV_h$ and 
		$\bif_{1,j}, \bif_{2,j}\in L^2\left( 0,T;\bH^{-1}(\cD)\right)$ for $j=1,2,\cdots, N_{sc}$. Set $M=T/\Delta t$ and for $n=1,2,\cdots, M-1$, 
		compute:\\
		Step 1: Find $\bhv_{j,h}^{n+1}\in \bX_h$ satisfying for all $\bchi_{h}\in \bX_h$:\vspace{-2mm}
		\begin{align}
			\Bigg(&\frac{\bhv_{j,h}^{n+1}-\btv_{j,h}^n}{\Delta t}, \bchi_{h}\Bigg)+b\left(<\bhw_h>^n, 
			\bhv_{j,h}^{n+1},\bchi_{h}\right)+\left(\frac{\Bar{\nu}+\Bar{\nu}_{m}}{2}\nabla \bhv_{j,h}^{n+1},\nabla 
			\bchi_{h}\right)\nonumber\\&+\gamma\lp\nabla\cdot\bhv_{j,h}^{n+1},\nabla\cdot\bchi_{h}\rp+\left(2\nu_T(\hat{\bw}^{'}_{h},t^n)\nabla 
			\bhv_{j,h}^{n+1},\nabla \bchi_{h}\right)= \left(\bif_{1,j}(t^{n+1}),\bchi_{h}\right)\nonumber\\&-b\left(\bhw_{j,h}^{'n}, 
			\bhv_{j,h}^n,\bchi_{h}\right)-\left(\frac{\nu_j-\nu_{m,j}}{2}\nabla \bhw_{j,h}^n,\nabla 
			\bchi_{h}\right)-\left(\frac{\nu_j^{'}+\nu_{m,j}^{'}}{2}\nabla \bhv_{j,h}^{n},\nabla\bchi_{h}\right).\label{weaknew1}
		\end{align}
		Step 2: Find $(\btv_{j,h}^{n+1},\hq_{j,h}^{n+1})\in \tilde{\bY}_h\times Q_h$ satisfying for all $(\bxi_{h},\zeta_{h})\in\tilde{\bY}_h\times 
		Q_h$:\vspace{-2mm}
		\begin{align} \Bigg(\frac{\btv_{j,h}^{n+1}-\bhv_{j,h}^{n+1}}{\Delta t}, \bxi_{h}\Bigg)-(\hq_{j,h}^{n+1},\nabla\cdot\bxi_{h}) 
			=&0,\label{weaknew2}\\(\nabla\cdot\btv_{j,h}^{n+1}, \zeta_{h})=&0.\label{weaknew3}
		\end{align}
		Step 3: Find $\bhw_{j,h}^{n+1}\in \bX_h$ satisfying for all $\bl_{h}\in \bX_h$:\vspace{-2mm}
		\begin{align}
			\Bigg(&\frac{\bhw_{j,h}^{n+1}-\btw_{j,h}^n}{\Delta t},\bl_{h}\Bigg)+b\left(<\bhv_h>^n, 
			\bhw_{j,h}^{n+1},\bl_{h}\right)+\left(\frac{\Bar{\nu}+\Bar{\nu}_{m}}{2}\nabla \bhw_{j,h}^{n+1},\nabla 
			\bl_{h}\right)\nonumber\\&+\gamma\lp\nabla\cdot\bhw_{j,h}^{n+1},\nabla\cdot\bl_{h}\rp+\left(2\nu_T(\bhv^{'}_{h},t^n)\nabla 
			\bhw_{j,h}^{n+1},\nabla \bl_{h}\right)= \left(\bif_{2,j}(t^{n+1}),\bl_{h}\right)\nonumber\\&-b\left(\bhv_{j,h}^{'n}, 
			\bhw_{j,h}^n,\bl_{h}\right)-\left(\frac{\nu_j-\nu_{m,j}}{2}\nabla \bhv_{j,h}^n,\nabla 
			\bl_{h}\right)-\left(\frac{\nu_j^{'}+\nu_{m,j}^{'}}{2}\nabla \bhw_{j,h}^{n},\nabla \bl_{h}\right).\label{weaknew4}
		\end{align}
		Step 4: Find $(\btw_{j,h}^{n+1},\hr_{j,h}^{n+1})\in \tilde{\bY}_h\times Q_h$ satisfying for all $(\bs_{h},\eta_{h})\in\tilde{\bY}_h\times 
		Q_h$:\vspace{-2mm}
		\begin{align}
			\Bigg(\frac{\btw_{j,h}^{n+1}-\bhw_{j,h}^{n+1}}{\Delta t}, \bs_{h}\Bigg)-(\hr_{j,h}^{n+1},\nabla\cdot\bs_{h}) 
			=&0,\label{weaknew5}\\(\nabla\cdot\btw_{j,h}^{n+1}, \eta_{h})=&0.\label{weaknew6}\vspace{-2mm}
		\end{align}\vspace{-2mm}
	\end{algorithm}
	Since $\bX_h\subset\tilde{\bY}_h$, we can choose $\bxi_{h}=\bchi_{h}$ in \eqref{weaknew2}, $\bs_{h}=\bl_{h}$ in \eqref{weaknew5} and combine them 
	with equations \eqref{weaknew1} and \eqref{weaknew4}, respectively, to get\vspace{-2mm}
	\begin{align}
		\Bigg(&\frac{\bhv_{j,h}^{n+1}-\bhv_{j,h}^n}{\Delta t}, \bchi_{h}\Bigg)+b\Big(<\bhw_h>^n, 
		\bhv_{j,h}^{n+1},\bchi_{h}\Big)+\left(\frac{\Bar{\nu}+\Bar{\nu}_{m}}{2}\nabla \bhv_{j,h}^{n+1},\nabla 
		\bchi_{h}\right)\nonumber\\&+\gamma\Big(\nabla\cdot\bhv_{j,h}^{n+1},\nabla\cdot\bchi_{h}\Big)+\Big(2\nu_T(\hat{\bw}^{'}_{h},t^n)\nabla 
		\bhv_{j,h}^{n+1},\nabla \bchi_{h}\Big)-\Big(\hq_{j,h}^{n},\nabla\cdot\bchi_{h}\Big)= 
		\Big(\bif_{1,j}(t^{n+1}),\bchi_{h}\Big)\nonumber\\&-b\Big(\bhw_{j,h}^{'n}, \bhv_{j,h}^n,\bchi_{h}\Big)-\left(\frac{\nu_j-\nu_{m,j}}{2}\nabla 
		\bhw_{j,h}^n,\nabla \bchi_{h}\right)-\left(\frac{\nu_j^{'}+\nu_{m,j}^{'}}{2}\nabla \bhv_{j,h}^{n},\nabla\bchi_{h}\right),\label{hatweak1}
	\end{align}
	and
	\begin{align}
		\Bigg(&\frac{\bhw_{j,h}^{n+1}-\bhw_{j,h}^n}{\Delta t},\bl_{h}\Bigg)+b\Big(<\bhv_h>^n, 
		\bhw_{j,h}^{n+1},\bl_{h}\Big)+\left(\frac{\Bar{\nu}+\Bar{\nu}_{m}}{2}\nabla \bhw_{j,h}^{n+1},\nabla 
		\bl_{h}\right)\nonumber\\&+\gamma\Big(\nabla\cdot\bhw_{j,h}^{n+1},\nabla\cdot\bl_{h}\Big)+\Big(2\nu_T(\bhv^{'}_{h},t^n)\nabla 
		\bhw_{j,h}^{n+1},\nabla \bl_{h}\Big)-\Big(\hr_{j,h}^{n},\nabla\cdot\bl_{h}\Big)= 
		\Big(\bif_{2,j}(t^{n+1}),\bl_{h}\Big)\nonumber\\&-b\Big(\bhv_{j,h}^{'n}, \bhw_{j,h}^n,\bl_{h}\Big)-\left(\frac{\nu_j-\nu_{m,j}}{2}\nabla 
		\bhv_{j,h}^n,\nabla \bl_{h}\right)-\left(\frac{\nu_j^{'}+\nu_{m,j}^{'}}{2}\nabla \bhw_{j,h}^{n},\nabla \bl_{h}\right).\label{hatweak2}
	\end{align}
	
	\subsection{Stability Analysis}\label{stability-analysis}
	We now prove stability and well-posedness for the Algorithm \ref{SPP-FEM}. 
	
	\begin{lemma}(Unconditional Stability)
		Let $\big(\bhv_{j,h}^{n+1},\hq_{j,h}^{n+1},\bhw_{j,h}^{n+1},\hr_{j,h}^{n+1}\big)$ be the solution of Algorithm \ref{SPP-FEM} and 
		$\textbf{f}_{1,j},\textbf{f}_{2,j}\in L^2\left(0,T;\bH^{-1}(\cD)\right)$, and $\bv_{j,h}^0$, $\bw_{j,h}^0 \in \bH^1(\cD)$ for $j=1,2,\cdots, 
		N_{sc}$. Then for all $\Delta t>0$, if $\alpha_j>0$, and $\mu>\frac{C}{2\Delta t\alpha_j}$, we have the following stability bound:
		\begin{align} \|\bhv_{j,h}^{M}\|^2&+\|\bhw_{j,h}^{M}\|^2+\frac{\Bar{\nu}_{\min}+\Bar{\nu}_{m,\min}}{2}\Delta t\Big(\|\nabla 
			\bhv_{j,h}^{M}\|^2+\|\nabla \bhw_{j,h}^{M}\|^2\Big)+2\gamma\Delta 
			t\sum_{n=0}^{M-1}\Big(\|\nabla\cdot\bhv_{j,h}^{n+1}\|^2+\|\nabla\cdot\bhw_{j,h}^{n+1}\|^2\Big)\nonumber\\&\le 
			\|\bhv_{j,h}^0\|^2+\|\bhw_{j,h}^0\|^2+\frac{\Bar{\nu}_{\min}+\Bar{\nu}_{m,\min}}{2}\Delta t\Big(\|\nabla \bhv_{j,h}^{0}\|^2+\|\nabla 
			\bhw_{j,h}^{0}\|^2\Big)\nonumber\\&+\frac{2\Delta 
				t}{\alpha_j}\sum_{n=0}^{M-1}\Big(\|\bif_{1,j}(t^{n+1})\|_{-1}^2+\|\bif_{2,j}(t^{n+1})\|_{-1}^2\Big).
		\end{align}
	\end{lemma}
	\begin{proof}
		Taking $\bchi_{h}=\bhv_{j,h}^{n+1}$ in \eqref{weaknew1} and $\bl_{h}=\bhw_{j,h}^{n+1}$ in \eqref{weaknew4}, to obtain
		\begin{align}
			\Bigg(\frac{\bhv_{j,h}^{n+1}-\btv_{j,h}^n}{\Delta t}, \bhv_{j,h}^{n+1}\Bigg)+\frac{1}{2}\|(\Bar{\nu}+\Bar{\nu}_{m})^\frac12\nabla 
			\bhv_{j,h}^{n+1}\|^2+\gamma\|\nabla\cdot\bhv_{j,h}^{n+1}\|^2\nonumber\\+\left(2\nu_T(\hat{\bw}^{'}_{h},t^n)\nabla \bhv_{j,h}^{n+1},\nabla 
			\bhv_{j,h}^{n+1}\right)=\Big(\bif_{1,j}(t^{n+1}),\bhv_{j,h}^{n+1}\Big)\nonumber\\-\Big(\bhw_{j,h}^{'n}\cdot\nabla 
			\bhv_{j,h}^n,\bhv_{j,h}^{n+1}\Big)-\Big(\frac{\nu_j-\nu_{m,j}}{2}\nabla \bhw_{j,h}^n,\nabla 
			\bhv_{j,h}^{n+1}\Big)-\Big(\frac{\nu_j^{'}+\nu_{m,j}^{'}}{2}\nabla \bhv_{j,h}^{n},\nabla\bhv_{j,h}^{n+1}\Big),
		\end{align}
		and
		\begin{align}
			\Bigg(\frac{\bhw_{j,h}^{n+1}-\btw_{j,h}^n}{\Delta t},\bhw_{j,h}^{n+1}\Bigg)+\frac{1}{2}\|(\Bar{\nu}+\Bar{\nu}_{m})^\frac12\nabla 
			\bhw_{j,h}^{n+1}\|^2+\gamma\|\nabla\cdot\bhw_{j,h}^{n+1}\|^2\nonumber\\+\Big(2\nu_T(\hat{\bv}^{'}_{h},t^n)\nabla \bhw_{j,h}^{n+1},\nabla 
			\bhw_{j,h}^{n+1}\Big)=\Big(\bif_{2,j}(t^{n+1}),\bhw_{j,h}^{n+1}\Big)\nonumber\\-\Big(\bhv_{j,h}^{'n}\cdot\nabla 
			\bhw_{j,h}^n,\bhw_{j,h}^{n+1}\Big)-\Big(\frac{\nu_j-\nu_{m,j}}{2}\nabla \bhv_{j,h}^n,\nabla 
			\bhw_{j,h}^{n+1}\Big)-\Big(\frac{\nu_j^{'}+\nu_{m,j}^{'}}{2}\nabla \bhw_{j,h}^{n},\nabla\bhw_{j,h}^{n+1}\Big).
		\end{align}
		Using polarization identity, \eqref{eddy-viscosity-term} and
		$(2\nu_T(\bhw^{'}_{h},t^n)\nabla \bhv_{j,h}^{n+1},\nabla \bhv_{j,h}^{n+1})=2\mu\Delta t\|l^{n}_{\bhw,h}\nabla \bhv_{j,h}^{n+1}\|^2$, we get 
		\begin{align}
			\frac{1}{2\Delta 
				t}\Big(\|\bhv_{j,h}^{n+1}\|^2-\|\btv_{j,h}^n\|^2+\|\bhv_{j,h}^{n+1}-\btv_{j,h}^n\|^2\Big)+\frac{1}{2}\|(\Bar{\nu}+\Bar{\nu}_{m})^\frac12\nabla 
			\bhv_{j,h}^{n+1}\|^2+\gamma\|\nabla\cdot\bhv_{j,h}^{n+1}\|^2\nonumber\\+2\mu\Delta t\|l^{n}_{\bhw,h}\nabla 
			\bhv_{j,h}^{n+1}\|^2=\Big(\bif_{1,j}(t^{n+1}),\bhv_{j,h}^{n+1}\Big)-b\Big(\bhw_{j,h}^{'n}, 
			\bhv_{j,h}^n,\bhv_{j,h}^{n+1}\Big)\nonumber\\-\Big(\frac{\nu_j-\nu_{m,j}}{2}\nabla \bhw_{j,h}^n,\nabla 
			\bhv_{j,h}^{n+1}\Big)-\Big(\frac{\nu_j^{'}+\nu_{m,j}^{'}}{2}\nabla \bhv_{j,h}^{n},\nabla\bhv_{j,h}^{n+1}\Big),\label{pol11}
		\end{align}
		and
		\begin{align}
			\frac{1}{2\Delta 
				t}\Big(\|\bhw_{j,h}^{n+1}\|^2-\|\btw_{j,h}^n\|^2+\|\bhw_{j,h}^{n+1}-\btw_{j,h}^n\|^2\Big)+\frac{1}{2}\|(\Bar{\nu}+\Bar{\nu}_{m})^\frac12\nabla 
			\bhw_{j,h}^{n+1}\|^2+\gamma\|\nabla\cdot\bhw_{j,h}^{n+1}\|^2\nonumber\\+2\mu\Delta t\|l^{n}_{\bhv,h}\nabla 
			\bhw_{j,h}^{n+1}\|^2=\Big(\bif_{2,j}(t^{n+1}),\bhw_{j,h}^{n+1}\Big)-b\Big(\bhv_{j,h}^{'n}, 
			\bhw_{j,h}^n,\bhw_{j,h}^{n+1}\Big)\nonumber\\-\Big(\frac{\nu_j-\nu_{m,j}}{2}\nabla \bhv_{j,h}^n,\nabla 
			\bhw_{j,h}^{n+1}\Big)-\Big(\frac{\nu_j^{'}+\nu_{m,j}^{'}}{2}\nabla \bhw_{j,h}^{n},\nabla\bhw_{j,h}^{n+1}\Big).\label{pol22}
		\end{align}
		Adding \eqref{pol11} and \eqref{pol22}, using $\|\ba\cdot\nabla\bb\|\le\sqrt{2}\||\ba|\nabla \bb\|$, Cauchy-Schwarz's, Poincar\'e inequality, \eqref{eddy-viscosity-term}, and Young's 
		inequality in
		\begin{align*}
			\left(\bhw_{j,h}^{'n}\cdot\nabla \bhv_{j,h}^n,\bhv_{j,h}^{n+1}\right)&=-\left(\bhw_{j,h}^{'n}\cdot\nabla \bhv_{j,h}^{n+1}, 
			\bhv_{j,h}^n\right)\\&\le\|\bhw_{j,h}^{'n}\cdot\nabla \bhv_{j,h}^{n+1}\|\hspace{1mm}\|\bhv_{j,h}^n\|\\&\le\sqrt{2}\||\bhw_{j,h}^{'n}|\nabla 
			\bhv_{j,h}^{n+1}\|\hspace{1mm}\|\bhv_{j,h}^n\|\nonumber\\&\le C\|l^{n}_{\bhw,h}\nabla 
			\bhv_{j,h}^{n+1}\|\hspace{1mm}\|\nabla\bhv_{j,h}^n\|\\&\le\frac{C}{\alpha_j}\|l^{n}_{\bhw,h}\nabla 
			\bhv_{j,h}^{n+1}\|^2+\frac{\alpha_j}{4}\|\nabla\bhv_{j,h}^n\|^2,
		\end{align*}
		and then applying the Cauchy-Schwarz inequality to the forcing term and H\"older's inequality  to the last two terms in the right-hand-side, 
		reduces to
		\begin{align}
			\frac{1}{2\Delta 
				t}\Big(\|\bhv_{j,h}^{n+1}\|^2-\|\btv_{j,h}^n\|^2+\|\bhv_{j,h}^{n+1}-\btv_{j,h}^n\|^2+\|\bhw_{j,h}^{n+1}\|^2-\|\btw_{j,h}^n\|^2+\|\bhw_{j,h}^{n+1}-\btw_{j,h}^n\|^2\Big)\nonumber\\+\frac{\Bar{\nu}_{\min}+\Bar{\nu}_{m,\min}}{2}\Big(\|\nabla 
			\bhv_{j,h}^{n+1}\|^2+\|\nabla 
			\bhw_{j,h}^{n+1}\|^2\Big)+\gamma\Big(\|\nabla\cdot\bhv_{j,h}^{n+1}\|^2+\|\nabla\cdot\bhw_{j,h}^{n+1}\|^2\Big)\nonumber\\+\left(2\mu\Delta 
			t-\frac{C}{\alpha_j}\right)\Big(\|l^{n}_{\bhw,h}\nabla \bhv_{j,h}^{n+1}\|^2+\|l^{n}_{\bhv,h}\nabla \bhw_{j,h}^{n+1}\|^2\Big)\nonumber\\\le 
			\frac{\alpha_j}{4}\|\nabla\bhv_{j,h}^n\|^2+\frac{\alpha_j}{4}\|\nabla\bhw_{j,h}^n\|^2+\|\bif_{1,j}(t^{n+1})\|_{-1}\|\nabla\bhv_{j,h}^{n+1}\|+\|\bif_{2,j}(t^{n+1})\|_{-1}\|\nabla\bhw_{j,h}^{n+1}\|\nonumber\\+\frac{\|\nu_j-\nu_{m,j}\|_\infty}{2}\Big(\|\nabla 
			\bhw_{j,h}^n\|\|\nabla \bhv_{j,h}^{n+1}\|+\|\nabla \bhv_{j,h}^n\|\|\nabla 
			\bhw_{j,h}^{n+1}\|\Big)\nonumber\\+\frac{\|\nu_j^{'}+\nu_{m,j}^{'}\|_{\infty}}{2}\Big(\|\nabla 
			\bhv_{j,h}^{n}\|\|\nabla\bhv_{j,h}^{n+1}\|+\|\nabla \bhw_{j,h}^{n}\|\|\nabla\bhw_{j,h}^{n+1}\|\Big).
		\end{align}
		Using Young's inequality and reducing, we have
		\begin{align}
			\frac{1}{2\Delta 
				t}\Big(\|\bhv_{j,h}^{n+1}\|^2-\|\btv_{j,h}^n\|^2+\|\bhv_{j,h}^{n+1}-\btv_{j,h}^n\|^2+\|\bhw_{j,h}^{n+1}\|^2-\|\btw_{j,h}^n\|^2+\|\bhw_{j,h}^{n+1}-\btw_{j,h}^n\|^2\Big)\nonumber\\+\frac{\Bar{\nu}_{\min}+\Bar{\nu}_{m,\min}}{4}\Big(\|\nabla 
			\bhv_{j,h}^{n+1}\|^2+\|\nabla 
			\bhw_{j,h}^{n+1}\|^2\Big)+\gamma\Big(\|\nabla\cdot\bhv_{j,h}^{n+1}\|^2+\|\nabla\cdot\bhw_{j,h}^{n+1}\|^2\Big)\nonumber\\+\left(2\mu\Delta 
			t-\frac{C}{\alpha_j}\right)\Big(\|l^{n}_{\bhw,h}\nabla \bhv_{j,h}^{n+1}\|^2+\|l^{n}_{\bhv,h}\nabla \bhw_{j,h}^{n+1}\|^2\Big)\nonumber\\\le 
			\frac{1}{\alpha_j}\Big(\|\bif_{1,j}(t^{n+1})\|_{-1}^2+\|\bif_{2,j}(t^{n+1})\|_{-1}^2\Big)+\frac{\Bar{\nu}_{\min}+\Bar{\nu}_{m,\min}}{4}\Big(\|\nabla 
			\bhv_{j,h}^{n}\|^2+\|\nabla \bhw_{j,h}^{n}\|^2\|\Big).
		\end{align}
		Assuming $\mu>\frac{C}{2\Delta t\alpha_j}$, and dropping non-negative terms from the left-hand-side, this reduces to
		\begin{align}
			\frac{1}{2\Delta 
				t}\Big(\|\bhv_{j,h}^{n+1}\|^2-\|\btv_{j,h}^n\|^2+\|\bhw_{j,h}^{n+1}\|^2-\|\btw_{j,h}^n\|^2\Big)\nonumber\\+\frac{\Bar{\nu}_{\min}+\Bar{\nu}_{m,\min}}{4}\Big(\|\nabla 
			\bhv_{j,h}^{n+1}\|^2-\|\nabla \bhv_{j,h}^{n}\|^2+\|\nabla \bhw_{j,h}^{n+1}\|^2-\|\nabla 
			\bhw_{j,h}^{n}\|^2\Big)\nonumber\\+\gamma\Big(\|\nabla\cdot\bhv_{j,h}^{n+1}\|^2+\|\nabla\cdot\bhw_{j,h}^{n+1}\|^2\Big)\le 
			\frac{1}{\alpha_j}\Big(\|\bif_{1,j}(t^{n+1})\|_{-1}^2+\|\bif_{2,j}(t^{n+1})\|_{-1}^2\Big).\label{proj1}
		\end{align}
		Now choose $\bxi_{h}=\btv_{j,h}^{n+1}$ in \eqref{weaknew2}, $\zeta_{h}=\hq_{j,h}^{n+1}$ in \eqref{weaknew3} and  $\bs_{h}=\btw_{j,h}^{n+1}$ in 
		\eqref{weaknew5}, $\eta_{h}=\hr_{j,h}^{n+1}$ in \eqref{weaknew6}.
		Then apply Cauchy-Schwarz and Young’s inequalities, to obtain
		\begin{align*}    \|\btv_{j,h}^{n+1}\|^2\le\|\bhv_{j,h}^{n+1}\|^2,~\text{and}~~
			\|\btw_{j,h}^{n+1}\|^2\le\|\bhw_{j,h}^{n+1}\|^2,
		\end{align*}
		for all $n=0,1,2,\cdots,M-1$. Plugging these estimates into \eqref{proj1}, results in 
		\begin{align}
			\frac{1}{2\Delta 
				t}\Big(\|\bhv_{j,h}^{n+1}\|^2-\|\bhv_{j,h}^n\|^2+\|\bhw_{j,h}^{n+1}\|^2-\|\bhw_{j,h}^n\|^2\Big)\nonumber\\+\frac{\Bar{\nu}_{\min}+\Bar{\nu}_{m,\min}}{4}\Big(\|\nabla 
			\bhv_{j,h}^{n+1}\|^2-\|\nabla \bhv_{j,h}^{n}\|^2+\|\nabla \bhw_{j,h}^{n+1}\|^2-\|\nabla 
			\bhw_{j,h}^{n}\|^2\Big)\nonumber\\+\gamma\Big(\|\nabla\cdot\bhv_{j,h}^{n+1}\|^2+\|\nabla\cdot\bhw_{j,h}^{n+1}\|^2\Big)\le 
			\frac{1}{\alpha_j}\Big(\|\bif_{1,j}(t^{n+1})\|_{-1}^2+\|\bif_{2,j}(t^{n+1})\|_{-1}^2\Big).
		\end{align}
		Multiplying both sides by $2\Delta t$ and summing over the time steps, completes the proof.
	\end{proof}
	
	We now prove the Algorithm \ref{SPP-FEM} converges to Algorithm \ref{Algn0} as $\gamma\rightarrow\infty$. Thus, we need to define the space 
	$\bR_h:=\bV_h^\perp\subset\bX_h$ to be the orthogonal complement of $\bV_h$ with respect to the $\bH^1(\cD)$ norm.
	
	\begin{lemma}\label{CR-lemma}
		Let the finite element pair $(\bX_h,Q_h)\subset(\bX,Q)$ satisfy the \textit{inf-sup condition} \eqref{infsup} and the divergence-free property, 
		i.e., $\nabla\cdot\bX_h\subset Q_h$. Then there exists a constant $C_R$ independent of $h$ such that $$\|\nabla\bv_h\|\le 
		C_R\|\nabla\cdot\bv_h\|,\hspace{3mm}\forall\bv_h\in \bR_h.$$
	\end{lemma}
	\begin{proof}
		See \cite{GR86, linke2017connection}
	\end{proof}
	
	\begin{assumption}\label{assump1}
		We assume there exists a constant $C_*$ which is independent of $h$, and $\Delta t$, such that for sufficiently small $h$ for a fixed mesh and 
		fixed $\Delta t$ as $\gamma\rightarrow\infty$, the solution of the Algorithm \ref{SPP-FEM} satisfies
		\begin{align}
			\max_{1\le n\le M}\Big\{\|\bhv_{j,h}^n\|_\infty,\|\bhw_{j,h}^n\|_\infty\Big\}&\le C_*,\hspace{2mm}\forall j=1,2,\cdots,{N_{sc}}.
		\end{align}
	\end{assumption}
	The Assumption \ref{assump1} is proved later in Lemma \ref{uniform-boundedness-lemma-proof}. The idea of utilizing the Assumption \ref{assump1} in the following convergence anaysis is taken from the finite element analysis of reaction-diffusion equation in \cite{mohebujjaman2024decoupled}.
	
	\begin{theorem}[Convergence]\label{gamma-convergence}
		Let $(\bv_{j,h}^{n+1},
		\bw_{j,h}^{n+1},q_{j,h}^{n+1})$, and $(\bhv_{j,h}^{n+1},
		\bhw_{j,h}^{n+1},\hq_{j,h}^{n+1})$ for $j=1,2,\cdots,N_{sc}$, are the solutions to the Algorithm \ref{Algn0}, and Algorithm \ref{SPP-FEM}, 
		respectively, for $n=0,1,\cdots,M-1$. We then have
		\begin{align}
			\Delta 
			t\sum_{n=1}^{M}\Big(\|\nabla\hspace{-1mm}<\hspace{-1mm}\bv_h\hspace{-1mm}>^n-\nabla\hspace{-1mm}<\hspace{-1mm}\bhv_h\hspace{-1mm}>^n\|^2+\|\nabla\hspace{-1mm}<\hspace{-1mm}\bw_h\hspace{-1mm}>^n-\nabla\hspace{-1mm}<\hspace{-1mm}\bhw_h\hspace{-1mm}>^n\|^2\Big)\nonumber\\\le 
			\frac{CC_R^2}{\gamma^2} \left(\frac{1}{\alpha_{\min}^3\Delta t}+\frac{1}{\alpha_{\min}\Delta t}+\frac{\Delta t}{\alpha_{\min}}+1\right) exp 
			\left(\frac{CC_*^2}{\alpha_{\min}}+\frac{C\Delta t}{h^3\alpha_{\min}}\right)\nonumber\\\times \lp  \Delta 
			t\sum_{n=0}^{M-1}\sum_{j=1}^J\Big(\|q_{j,h}^{n+1}-\hq_{j,h}^n\|^2+\|r_{j,h}^{n+1}-\hr_{j,h}^n\|^2\Big)\rp.
		\end{align}
	\end{theorem}
	
	\begin{remark}
		The above theorem states the first order convergence of the penalty-projection algorithm to the Algorithm \ref{Algn0} as 
		$\gamma\rightarrow\infty$ for a fixed mesh and time-step size.
	\end{remark}
	
	\begin{proof}
		Denote $\be_{j}^{n+1}:=\bv_{j,h}^{n+1}-\bhv_{j,h}^{n+1}$, and $\bep_{j}^{n+1}:=\bw_{j,h}^{n+1}-\bhw_{j,h}^{n+1}$ and use the following 
		$H^1$-orthogonal decomposition of the errors:
		$$\be_{j}^{n+1}:=\be_{j,0}^{n+1}+\be_{j,\bR}^{n+1}, \;\text{and}\;\bep_{j}^{n+1}:=\bep_{j,0}^{n+1}+\bep_{j,\bR}^{n+1},$$
		with $\be_{j,0}^{n+1},\bep_{j,0}^{n+1}\in\bV_h$, and $\be_{j,\bR}^{n+1},\bep_{j,\bR}^{n+1}\in\bR_h$, for $n=0,1,\cdots,M-1$.\\
		\textbf{Step 1:} Estimate of $\be_{j,\bR}^{n+1}$,  and $\bep_{j,\bR}^{n+1}$: Subtracting the equation \eqref{weak1} from \eqref{hatweak1} and 
		\eqref{weak2} from \eqref{hatweak2} produces
		\begin{align}
			\frac{1}{\Delta t}\Big(\be_j^{n+1}&-\be_j^n,\bchi_{h}\Big)+\left(\frac{\Bar{\nu}+\Bar{\nu}_{m}}{2}\nabla \be_j^{n+1},\nabla 
			\bchi_{h}\right)+\gamma\Big(\nabla\cdot\be_{j,\bR}^{n+1},\nabla\cdot\bchi_{h}\Big)+b\Big(\hspace{-1.5mm}\lab\bhw_h\rab^n,\be^{n+1}_j,\bchi_{h}\Big)\nonumber\\&+b\Big(\hspace{-1.5mm}\lab\bep\rab^n,\bv_{j,h}^{n+1},\bchi_{h}\Big)-\Big(q_{j,h}^{n+1}-\hq_{j,h}^n,\nabla\cdot\bchi_{h}\Big)+2\mu\Delta 
			t\Big((l_{\bhw,h}^n)^2\nabla\be_j^{n+1},\nabla\bchi_{h}\Big)\nonumber\\&+2\mu\Delta 
			t\Big(\big\{(l_{\bw,h}^n)^2-(l_{\bhw,h}^n)^2\big\}\nabla\bv_{j,h}^{n+1},\nabla\bchi_{h}\Big)=-b\Big(\bhw^{'n}_{j,h},\be^n_j,\bchi_{h}\Big)-b\Big(\bep^{'n}_{j},\bv_{j,h}^n,\bchi_{h}\Big)\nonumber\\&-\Big(\frac{\nu_j-\nu_{m,j}}{2}\nabla 
			\bep_{j}^n,\nabla \bchi_{h}\Big)-\Big(\frac{\nu_j^{'}+\nu_{m,j}^{'}}{2}\nabla \be_{j}^{n},\nabla\bchi_{h}\Big),\label{step-1-eqn-1}
		\end{align}
		and
		\begin{align}
			\frac{1}{\Delta t}\Big(\bep_j^{n+1}&-\bep_j^n,\bl_{h}\Big)+\left(\frac{\Bar{\nu}+\Bar{\nu}_{m}}{2}\nabla \bep_j^{n+1},\nabla 
			\bl_{h}\right)+\gamma\Big(\nabla\cdot\bep_{j,\bR}^{n+1},\nabla\cdot\bl_{h}\Big)+b\Big(\hspace{-1.5mm}\lab\bhv_h\rab^n,\bep^{n+1}_j,\bl_{h}\Big)\nonumber\\&+b\Big(\hspace{-1.5mm}\lab\be\rab^n,\bw_{j,h}^{n+1},\bl_{h}\Big)-\Big(r_{j,h}^{n+1}-\hr_{j,h}^n,\nabla\cdot\bl_{h}\Big)+2\mu\Delta 
			t\Big((l_{\bhv,h}^n)^2\nabla\bep_j^{n+1},\nabla\bl_{h}\Big)\nonumber\\&+2\mu\Delta 
			t\Big(\big\{(l_{\bv,h}^n)^2-(l_{\bhv,h}^n)^2\big\}\nabla\bw_{j,h}^{n+1},\nabla\bl_{h}\Big)=-b\Big(\bhv^{'n}_{j,h},\bep^n_j,\bl_{h}\Big)-b\Big(\be^{'n}_{j},\bw_{j,h}^n,\bl_{h}\Big)\nonumber\\&-\Big(\frac{\nu_j-\nu_{m,j}}{2}\nabla 
			\be_{j}^n,\nabla \bl_{h}\Big)-\Big(\frac{\nu_j^{'}+\nu_{m,j}^{'}}{2}\nabla \bep_{j}^{n},\nabla\bl_{h}\Big). \label{step-1-eqn-2}
		\end{align}
		Take $\bchi_{h}=\be_j^{n+1}$ in \eqref{step-1-eqn-1}, and $\bl_{h}=\bep_j^{n+1}$ in \eqref{step-1-eqn-2}, which yield
		\begin{align*}
			b\Big(\hspace{-1.5mm}\lab\bhv_h\rab^n,\bep^{n+1}_j,\bl_{h}\Big)=0,\hspace{2mm}\text{and}\hspace{2mm}
			b\Big(\hspace{-1.5mm}\lab\bhw_h\rab^n,\be^{n+1}_j,\bchi_{h}\Big)=0,
		\end{align*}
		and use polarization identity to get
		\begin{align}
			\frac{1}{2\Delta t}\Big(\|\be_j^{n+1}\|^2-\|\be_j^n\|^2+\|\be_j^{n+1}-\be_j^n\|^2\Big)+\frac{1}{2}\|(\Bar{\nu}+\Bar{\nu}_{m})^\frac12\nabla 
			\be_j^{n+1}\|^2+\gamma\|\nabla\cdot\be_{j,\bR}^{n+1}\|^2\nonumber\\+b\Big(\hspace{-1.5mm}\lab\bep\rab^n,\bv_{j,h}^{n+1},\be_j^{n+1}\Big)-\Big(q_{j,h}^{n+1}-\hq_{j,h}^n,\nabla\cdot\be_{j,\bR}^{n+1}\Big)+2\mu\Delta 
			t\|l_{\bhw,h}^n\nabla\be_j^{n+1}\|^2\nonumber\\+2\mu\Delta 
			t\Big(\big\{(l_{\bw,h}^n)^2-(l_{\bhw,h}^n)^2\big\}\nabla\bv_{j,h}^{n+1},\nabla\be_j^{n+1}\Big)=-b\Big(\bhw^{'n}_{j,h},\be^n_j,\be_j^{n+1}\Big)\nonumber\\-b\Big(\bep^{'n}_{j},\bv_{j,h}^n,\be_j^{n+1}\Big)-\Big(\frac{\nu_j-\nu_{m,j}}{2}\nabla 
			\bep_{j}^n,\nabla \be_j^{n+1}\Big)-\Big(\frac{\nu_j^{'}+\nu_{m,j}^{'}}{2}\nabla \be_{j}^{n},\nabla\be_j^{n+1}\Big),\label{pol-1}
		\end{align}
		and
		\begin{align}
			\frac{1}{2\Delta 
				t}\Big(\|\bep_j^{n+1}\|^2-\|\bep_j^n\|^2+\|\bep_j^{n+1}-\bep_j^n\|^2\Big)+\frac{1}{2}\|(\Bar{\nu}+\Bar{\nu}_{m})^\frac12\nabla 
			\bep_j^{n+1}\|^2+\gamma\|\nabla\cdot\bep_{j,\bR}^{n+1}\|^2\nonumber\\
			+b\Big(\hspace{-1.5mm}\lab\be\rab^n,\bw_{j,h}^{n+1},\bep_j^{n+1}\Big)-\Big(r_{j,h}^{n+1}-\hr_{j,h}^n,\nabla\cdot\bep_{j,\bR}^{n+1}\Big)+2\mu\Delta 
			t\|l_{\bhv,h}^n\nabla\bep_j^{n+1}\|^2\nonumber\\+2\mu\Delta t\Big(\big\{(l_{\bv,h}^n)^2-(l_{\bhv, 
				h}^n)^2\big\}\nabla\bw_{j,h}^{n+1},\nabla\bep_j^{n+1}\Big)=-b\Big(\bhv^{'n}_{j,h},\bep^n_j,\bep_j^{n+1}\Big)\nonumber\\-b\Big(\be^{'n}_{j},\bw_{j,h}^n,\bep_j^{n+1}\Big)-\Big(\frac{\nu_j-\nu_{m,j}}{2}\nabla 
			\be_{j}^n,\nabla \bep_j^{n+1}\Big)-\Big(\frac{\nu_j^{'}+\nu_{m,j}^{'}}{2}\nabla \bep_{j}^{n},\nabla\bep_j^{n+1}\Big).\label{pol-2}
		\end{align}
		Now, we find the bound of the terms in \eqref{pol-1} first.  Rearranging and applying Cauchy-Schwarz inequality, \eqref{eddy-viscosity-term}, and Young’s inequality in the following 
		nonlinear term yields
		\begin{align*}
			-b\Big(\bhw^{'n}_{j,h},\be^n_j,\be_j^{n+1}\Big)&= 
			-b\Big(\bhw^{'n}_{j,h},\be_j^{n+1},\be_j^{n+1}-\be^n_j\Big)\\&\le\|\bhw^{'n}_{j,h}\cdot\nabla\be_j^{n+1}\|\|\be_j^{n+1}-\be^n_j\|\\&\le\sqrt{2}\|l_{\bhw,h}^n\nabla\be_j^{n+1}\|\|\be_j^{n+1}-\be^n_j\|\\&\le 
			2\Delta t\|l_{\bhw,h}^n\nabla\be_j^{n+1}\|^2+\frac{1}{4\Delta t}\|\be_j^{n+1}-\be^n_j\|^2.
		\end{align*}
		Applying H\"older's and Young’s inequalities, we have
		\begin{align*}
			\Big|\Big(\frac{\nu_j-\nu_{m,j}}{2}\nabla \bep_{j}^n,\nabla \be_j^{n+1}\Big)\Big|&\le\frac{\|\nu_j-\nu_{m,j}\|_\infty}{4}\Big(\|\nabla 
			\bep_{j}^n\|^2+\|\nabla \be_j^{n+1}\|^2\Big),\\
			\Big|\Big(\frac{\nu_j^{'}+\nu_{m,j}^{'}}{2}\nabla \be_{j}^{n},\nabla\be_j^{n+1}\Big)\Big|&\le 
			\frac{\|\nu_j^{'}+\nu_{m,j}^{'}\|_\infty}{4}\Big(\|\nabla \be_{j}^{n}\|^2+\|\nabla \be_j^{n+1}\|^2\Big).
		\end{align*}
		Applying Cauchy-Schwarz and Young’s inequalities, we have
		\begin{align*}    
			\Big|\Big(q_{j,h}^{n+1}-\hq_{j,h}^n,\nabla\cdot\be_{j,\bR}^{n+1}\Big)\Big|&\le\frac{1}{2\gamma}\|q_{j,h}^{n+1}-\hq_{j,h}^n\|^2+\frac{\gamma}{2}\|\nabla\cdot\be_{j,\bR}^{n+1}\|^2.
		\end{align*}
		Using H\"older's inequality, estimate in Lemma \ref{lemma1},  Sobolev embedding theorem, Poincar\'e, and Young's inequalities provides
		\begin{align*}
			\Big|b\Big(\hspace{-1.5mm}\lab\bep\rab^n,\bv_{j,h}^{n+1},\be_j^{n+1}\Big)\Big|&\le  
			\|\hspace{-1.mm}\lab\bep\rab^n\hspace{-1.mm}\|\|\nabla\bv_{j,h}^{n+1}\|_{L^3}\|\be_j^{n+1}\|_{L^6}\\&\le 
			CC_*\|\hspace{-1.mm}\lab\bep\rab^n\hspace{-1.mm}\|\|\nabla\be_j^{n+1}\|\\
			&\le \frac{\alpha_j}{12}\|\nabla\be_j^{n+1}\|^2+\frac{CC_*^2}{\alpha_j}\|\hspace{-1.mm}\lab\bep\rab^n\hspace{-1.mm}\|^2,\\
			\Big|b\Big(\bep^{'n}_{j},\bv_{j,h}^n,\be_j^{n+1}\Big)\Big|&\le \|\bep^{'n}_{j}\|\|\nabla\bv_{j,h}^n\|_{L^3}\|\be_j^{n+1}\|_{L^6}\\&\le CC_* 
			\|\bep^{'n}_{j}\|\|\nabla\be_j^{n+1}\|\\&\le \frac{\alpha_j}{12}\|\nabla\be_j^{n+1}\|^2+\frac{CC_*^2}{\alpha_j}\|\bep^{'n}_{j}\|^2.
		\end{align*}
		For the third non-linear term, we apply H\"older’s and triangle inequalities, the stability estimate of Algorithm \ref{Algn0}, uniform 
		boundedness in Lemma \ref{lemma1} and in Assumption \ref{assump1}, Agmon’s \cite{Robinson2016Three-Dimensional}, discrete inverse, and Young's 
		inequalities,  to get
		\begin{align}
			2\mu\Delta t\Big(\big\{(l_{\bw,h}^n)^2-&(l_{\bhw,h}^n)^2\big\}\nabla\bv_{j,h}^{n+1},\nabla\be_{j}^{n+1}\Big)\nonumber\\
			&\le 2\mu\Delta t\|(l_{\bw,h}^n)^2-(l_{\bhw,h}^n)^2\|_{\infty}\|\nabla\bv_{j,h}^{n+1}\|\|\nabla\be_j^{n+1}\|\nonumber\\
			&=2\mu\Delta t 
			\|\sum_{i=1}^{N_{sc}}\left(|\bw_{i,h}^{'n}|^2-|\bhw_{i,h}^{'n}|^2\right)\|_\infty\|\nabla\bv_{j,h}^{n+1}\|\|\nabla\be_j^{n+1}\|\nonumber\\
			&\le 2\mu\Delta t 
			\sum_{i=1}^{N_{sc}}\|(\bw_{i,h}^{'n}-\bhw_{i,h}^{'n})\cdot(\bw_{i,h}^{'n}+\bhw_{i,h}^{'n})\|_\infty\|\nabla\bv_{j,h}^{n+1}\|\|\nabla\be_j^{n+1}\|\nonumber\\&\le 
			2\mu\Delta t 
			\sum_{i=1}^{N_{sc}}\|\bw_{i,h}^{'n}-\bhw_{i,h}^{'n}\|_\infty\|\bw_{i,h}^{'n}+\bhw_{i,h}^{'n}\|_\infty\|\nabla\bv_{j,h}^{n+1}\|\|\nabla\be_j^{n+1}\|\nonumber\\
			&\le C\Delta t^{\frac{1}{2}} 
			\sum_{i=1}^{N_{sc}}\|\bep_i^{'n}\|_\infty\left(\|\bw_{i,h}^{'n}\|_\infty+\|\bhw_{i,h}^{'n}\|_\infty\right)\|\nabla\be_j^{n+1}\|\nonumber\\
			&\le C\Delta t^{\frac{1}{2}}\sum_{i=1}^{N_{sc}}\|\bep_i^n\|_\infty\|\nabla\be_j^{n+1}\|\le C\Delta 
			t^\frac{1}{2}h^{-\frac{3}{2}}\sum_{i=1}^{N_{sc}}\|\bep_i^n\|\|\nabla\be_j^{n+1}\|\nonumber\\
			&\le\frac{\alpha_j}{12}\|\nabla\be_j^{n+1}\|^2+\frac{C\Delta 
				t}{h^3\alpha_j}\sum_{i=1}^{N_{sc}}\|\bep_i^n\|^2.\label{mixing-length-difference}
		\end{align}
		Using the above estimates in \eqref{pol-1}, choosing $\mu>\max\{ 1,\frac{1}{2\Delta t}\}$, dropping non-negative terms and reducing produces
		\begin{align}
			\frac{1}{2\Delta t}\Big(\|\be_j^{n+1}\|^2-\|\be_j^n\|^2\Big)+\frac{\Bar{\nu}_{\min}+\Bar{\nu}_{m,\min}}{4}\|\nabla 
			\be_j^{n+1}\|^2+\frac{\gamma}{2}\|\nabla\cdot\be_{j,\bR}^{n+1}\|^2\le \frac{\|\nu_j-\nu_{m,j}\|_\infty}{4}\|\nabla 
			\bep_{j}^n\|^2\nonumber\\+\frac{\|\nu_j^{'}+\nu_{m,j}^{'}\|_\infty}{4}\|\nabla 
			\be_{j}^{n}\|^2+\frac{1}{2\gamma}\|q_{j,h}^{n+1}-\hq_{j,h}^n\|^2+\frac{CC_*^2}{\alpha_j}\Big(\|\hspace{-1.mm}\lab\bep\rab^n\hspace{-1.mm}\|^2+\|\bep^{'n}_{j}\|^2\Big)+\frac{C\Delta 
				t}{h^3\alpha_j}\sum_{i=1}^{N_{sc}}\|\bep_i^n\|^2.\label{upperbound1}
		\end{align}
		Now, apply similar estimates to the right-hand-side terms of \eqref{pol-2}, to produce
		\begin{align}
			\frac{1}{2\Delta t}\Big(\|\bep_j^{n+1}\|^2-\|\bep_j^n\|^2\Big)+\frac{\Bar{\nu}_{\min}+\Bar{\nu}_{m,\min}}{4}\|\nabla 
			\bep_j^{n+1}\|^2+\frac{\gamma}{2}\|\nabla\cdot\bep_{j,\bR}^{n+1}\|^2\le \frac{\|\nu_j-\nu_{m,j}\|_\infty}{4}\|\nabla 
			\be_{j}^n\|^2\nonumber\\+\frac{\|\nu_j^{'}+\nu_{m,j}^{'}\|_\infty}{4}\|\nabla 
			\bep_{j}^{n}\|^2+\frac{1}{2\gamma}\|r_{j,h}^{n+1}-\hr_{j,h}^n\|^2+\frac{CC_*^2}{\alpha_j}\Big(\|\hspace{-1.mm}\lab\be\rab^n\hspace{-1.mm}\|^2+\|\be^{'n}_{j}\|^2\Big)+\frac{C\Delta 
				t}{h^3\alpha_j}\sum_{i=1}^{N_{sc}}\|\be_i^n\|^2.\label{upperbound2}
		\end{align}
		Add \eqref{upperbound1} and \eqref{upperbound2}, and rearrange
		\begin{align}
			&\frac{1}{2\Delta 
				t}\Big(\|\be_j^{n+1}\|^2-\|\be_j^n\|^2+\|\bep_j^{n+1}\|^2-\|\bep_j^n\|^2\Big)+\frac{\Bar{\nu}_{\min}+\Bar{\nu}_{m,\min}}{4}\Big(\|\nabla 
			\be_j^{n+1}\|^2-\|\nabla \be_j^{n}\|^2+\|\nabla \bep_j^{n+1}\|^2-\|\nabla \bep_j^{n}\|^2\Big)\nonumber\\&+\frac{\alpha_j}{4}\Big(\|\nabla 
			\be_j^{n}\|^2+\|\nabla \bep_j^{n}\|^2\Big)+\frac{\gamma}{2}\left(\|\nabla\cdot\be_{j,\bR}^{n+1}\|^2+\|\nabla\cdot\bep_{j,\bR}^{n+1}\|^2\right)\le 
			\frac{1}{2\gamma}\left(\|q_{j,h}^{n+1}-\hq_{j,h}^n\|^2+\|r_{j,h}^{n+1}-\hr_{j,h}^n\|^2\right)\nonumber\\&+\frac{CC_*^2}{\alpha_j}\Big(\|\hspace{-1.mm}\lab\be\rab^n\hspace{-1.mm}\|^2+\|\be^{'n}_{j}\|^2+\|\hspace{-1.mm}\lab\bep\rab^n\hspace{-1.mm}\|^2+\|\bep^{'n}_{j}\|^2\Big)+\frac{C\Delta 
				t}{h^3\alpha_j}\sum_{i=1}^{N_{sc}}\Big(\|\be_i^n\|^2+\|\bep_i^n\|^2\Big).
		\end{align}
		Now, multiply both sides by $2\Delta t$, and sum over the time steps $n=0,1,\cdots,M-1$, to get
		\begin{align}
			\|\be_j^{M}\|^2+\|\bep_j^{M}\|^2+\frac{\Bar{\nu}_{\min}+\Bar{\nu}_{m,\min}}{2}\Delta 
			t\Big(\|\nabla\be_j^{M}\|^2+\|\nabla\bep_j^{M}\|^2\Big)+\frac{\alpha_j}{2}\Delta 
			t\sum_{n=0}^{M-1}\Big(\|\nabla\be_j^{n}\|^2+\|\nabla\bep_j^{n}\|^2\Big)\nonumber\\+\Delta 
			t\sum_{n=0}^{M-1}\gamma\Big(\|\nabla\cdot\be_{j,\bR}^{n+1}\|^2+\|\nabla\cdot\bep_{j,\bR}^{n+1}\|^2\Big)\le\frac{\Delta 
				t}{\gamma}\sum_{n=0}^{M-1}\Big(\|q_{j,h}^{n+1}-\hq_{j,h}^n\|^2+\|r_{j,h}^{n+1}-\hr_{j,h}^n\|^2\Big)\nonumber\\+\frac{CC_*^2}{\alpha_j}\Delta 
			t\sum_{n=0}^{M-1}\Big(\|\hspace{-1mm}\lab\be\rab^n\hspace{-1mm}\|^2+\|\be^{'n}_{j}\|^2+\|\hspace{-1mm}\lab\bep\rab^n\hspace{-1mm}\|^2+\|\bep^{'n}_{j}\|^2\Big)+\frac{C\Delta 
				t^2}{h^3\alpha_j}\sum_{n=1}^{M-1}\sum_{i=1}^{N_{sc}}\Big(\|\be_i^n\|^2+\|\bep_i^n\|^2\Big).
		\end{align} 
		Using triangle, Cauchy-Schwarz, and Young's inequalities, to get
		\begin{align}
			\|\be_j^{M}\|^2+\|\bep_j^{M}\|^2+\frac{\alpha_j\Delta t}{2}\sum_{n=1}^{M}\Big(\|\nabla\be_j^{n}\|^2+\|\nabla\bep_j^{n}\|^2\Big)+\Delta 
			t\sum_{n=1}^{M}\gamma\Big(\|\nabla\cdot\be_{j,\bR}^{n}\|^2+\|\nabla\cdot\bep_{j,\bR}^{n}\|^2\Big)\nonumber\\\le\frac{\Delta 
				t}{\gamma}\sum_{n=0}^{M-1}\Big(\|q_{j,h}^{n+1}-\hq_{j,h}^n\|^2+\|r_{j,h}^{n+1}-\hr_{j,h}^n\|^2\Big)+\left(\frac{CC_*^2}{\alpha_j}\Delta 
			t+\frac{C\Delta t^2}{h^3\alpha_j}\right)\sum_{n=1}^{M-1}\sum_{j=1}^{N_{sc}}\Big(\|\be_j^n\|^2+\|\bep_j^n\|^2\Big).
		\end{align} 
		Summing over $j=1,2,\cdots, {N_{sc}}$, we have
		\begin{align}   \sum_{j=1}^{N_{sc}}\|\be_j^{M}\|^2+\sum_{j=1}^{N_{sc}}\|\bep_j^{M}\|^2+\frac{\alpha_{\min}\Delta 
				t}{2}\sum_{n=1}^{M}\sum_{j=1}^{N_{sc}}\Big(\|\nabla\be_j^{n}\|^2+\|\nabla\bep_j^{n}\|^2\Big)\nonumber\\+\gamma\Delta 
			t\sum_{n=1}^{M}\sum_{j=1}^{N_{sc}}\Big(\|\nabla\cdot\be_{j,\bR}^{n}\|^2+\|\nabla\cdot\bep_{j,\bR}^{n}\|^2\Big)\le\frac{\Delta 
				t}{\gamma}\sum_{n=0}^{M-1}\sum_{j=1}^{N_{sc}}\Big(\|q_{j,h}^{n+1}-\hq_{j,h}^n\|^2+\|r_{j,h}^{n+1}-\hr_{j,h}^n\|^2\Big)\nonumber\\+\Delta 
			t\sum_{n=1}^{M-1}\left(\frac{CC_*^2}{\alpha_{\min}}+\frac{C\Delta 
				t}{h^3\alpha_{\min}}\right)\sum_{j=1}^{N_{sc}}\Big(\|\be_j^n\|^2+\|\bep_j^n\|^2\Big).
		\end{align} 
		Apply discrete Gr\"onwall inequality given in Lemma \ref{dgl}, to get 
		\begin{align} \sum_{j=1}^{N_{sc}}&\|\be_j^{M}\|^2+\sum_{j=1}^{N_{sc}}\|\bep_j^{M}\|^2+\frac{\alpha_{\min}\Delta 
				t}{2}\sum_{n=1}^{M}\sum_{j=1}^{N_{sc}}\Big(\|\nabla\be_j^{n}\|^2+\|\nabla\bep_j^{n}\|^2\Big)+\gamma\Delta 
			t\sum_{n=1}^{M}\sum_{j=1}^{N_{sc}}\Big(\|\nabla\cdot\be_{j,\bR}^{n}\|^2+\|\nabla\cdot\bep_{j,\bR}^{n}\|^2\Big)\nonumber\\&\le \frac{1}{\gamma} 
			exp\lp CT\left(\frac{C_*^2}{\alpha_{\min}}+\frac{\Delta t}{h^3\alpha_{\min}}\right)\rp\lp  \Delta 
			t\sum_{n=0}^{M-1}\sum_{j=1}^{N_{sc}}\Big(\|q_{j,h}^{n+1}-\hq_{j,h}^n\|^2+\|r_{j,h}^{n+1}-\hr_{j,h}^n\|^2\Big)\rp.\label{after-gronwall}
		\end{align}
		Using Lemma \ref{CR-lemma} with \eqref{after-gronwall} yields the following bound
		\begin{align}
			\Delta t\sum_{n=1}^{M}\sum_{j=1}^{N_{sc}}\Big(\|\nabla\be_{j,\bR}^{n}\|^2+\|\nabla\bep_{j,\bR}^{n}\|^2\Big)\le C_R^2\Delta 
			t\sum_{n=1}^{M}\sum_{j=1}^{N_{sc}}\Big(\|\nabla\cdot\be_{j,\bR}^{n}\|^2+\|\nabla\cdot\bep_{j,\bR}^{n}\|^2\Big)\nonumber\\\le\frac{C_R^2}{\gamma^2} 
			exp\lp \frac{CC_*^2}{\alpha_{\min}}+\frac{C\Delta t}{h^3\alpha_{\min}}\rp\lp  \Delta 
			t\sum_{n=0}^{M-1}\sum_{j=1}^{N_{sc}}\Big(\|q_{j,h}^{n+1}-\hq_{j,h}^n\|^2+\|r_{j,h}^{n+1}-\hr_{j,h}^n\|^2\Big)\rp.\label{step-1-bound}
		\end{align}
		\textbf{Step 2:} Estimate of $\be_{j,0}^{n}$, and $\bep_{j,0}^{n}$: To find a bound on $\Delta 
		t\sum\limits_{n=1}^{M}\sum\limits_{j=1}^{N_{sc}}\Big(\|\nabla\be_{j,0}^{n}\|^2+\|\nabla\bep_{j,0}^{n}\|^2\Big),$ take $\bchi_{h}=\be_{j,0}^{n+1}$ 
		in \eqref{step-1-eqn-1}, and $\bl_{h}=\bep_{j,0}^{n+1}$ in \eqref{step-1-eqn-2}, which yield
		\begin{align}
			&\frac{1}{\Delta t}\Big(\be_j^{n+1}-\be_j^n,\be_{j,0}^{n+1}\Big)+\frac{1}{2}\|(\Bar{\nu}+\Bar{\nu}_{m})^\frac12\nabla 
			\be_{j,0}^{n+1}\|^2=-b\Big(\hspace{-1.5mm}\lab\bhw_h\rab^n,\be^{n+1}_{j,\bR},\be_{j,0}^{n+1}\Big)-b\Big(\hspace{-1.5mm}\lab\bep\rab^n,\bv_{j,h}^{n+1},\be_{j,0}^{n+1}\Big)\nonumber\\&-2\mu\Delta 
			t\Big((l_{\bhw,h}^n)^2\nabla\be_j^{n+1},\nabla\be_{j,0}^{n+1}\Big)-2\mu\Delta 
			t\Big(\big\{(l_{\bw,h}^n)^2-(l_{\bhw,h}^n)^2\big\}\nabla\bv_{j,h}^{n+1},\nabla\be_{j,0}^{n+1}\Big)-b\Big(\bhw^{'n}_{j,h},\be^n_j,\be_{j,0}^{n+1}\Big)\nonumber\\&-b\Big(\bep^{'n}_{j},\bv_{j,h}^n,\be_{j,0}^{n+1}\Big)-\Big(\frac{\nu_j-\nu_{m,j}}{2}\nabla 
			\bep_{j,0}^n,\nabla \be_{j,0}^{n+1}\Big)-\Big(\frac{\nu_j^{'}+\nu_{m,j}^{'}}{2}\nabla 
			\be_{j,0}^{n},\nabla\be_{j,0}^{n+1}\Big),\label{step-2-eqn-1}
		\end{align}
		and
		\begin{align}
			\frac{1}{\Delta t}&\Big(\bep_j^{n+1}-\bep_j^n,\bep_{j,0}^{n+1}\Big)+\frac{1}{2}\|(\Bar{\nu}+\Bar{\nu}_{m})^\frac12\nabla 
			\bep_{j,0}^{n+1}\|^2=-b\Big(\hspace{-1.5mm}\lab\bhv_h\rab^n,\bep^{n+1}_{j,\bR},\bep_{j,0}^{n+1}\Big)-b\Big(\hspace{-1.5mm}\lab\be\rab^n,\bw_{j,h}^{n+1},\bep_{j,0}^{n+1}\Big)\nonumber\\&-2\mu\Delta 
			t\Big((l_{\bhv,h}^n)^2\nabla\bep_j^{n+1},\nabla\bep_{j,0}^{n+1}\Big)-2\mu\Delta 
			t\Big(\big\{(l_{\bv,h}^n)^2-(l_{\bhv,h}^n)^2\big\}\nabla\bw_{j,h}^{n+1},\nabla\bep_{j,0}^{n+1}\Big)-b\Big(\bhv^{'n}_{j,h},\bep^n_j,\bep_{j,0}^{n+1}\Big)\nonumber\\&-b\Big(\be^{'n}_{j},\bw_{j,h}^n,\bep_{j,0}^{n+1}\Big)-\Big(\frac{\nu_j-\nu_{m,j}}{2}\nabla 
			\be_{j,0}^n,\nabla \bep_{j,0}^{n+1}\Big)-\Big(\frac{\nu_j^{'}+\nu_{m,j}^{'}}{2}\nabla \bep_{j,0}^{n},\nabla\bep_{j,0}^{n+1}\Big). 
			\label{step-2-eqn-2}
		\end{align}
		Apply the non-linear bound given in \eqref{nonlinearbound}, and H\"older's inequality for the first, and second non-linear terms of 
		\eqref{step-2-eqn-1}, respectively, to obtain
		\begin{align}
			&\frac{1}{\Delta t}\Big(\be_j^{n+1}-\be_j^n,\be_{j,0}^{n+1}\Big)+\frac{\Bar{\nu}_{\min}+\Bar{\nu}_{m,\min}}{2}\|\nabla 
			\be_{j,0}^{n+1}\|^2+2\mu\Delta t\|l_{\bhw,h}^n\nabla\be_{j,0}^{n+1}\|^2\nonumber\\&\le 
			C\|\nabla\hspace{-1mm}\lab\bhw_h\rab^n\hspace{-1mm}\|\|\nabla\be^{n+1}_{j,\bR}\|\|\nabla\be_{j,0}^{n+1}\|+C\|\hspace{-1.mm}\lab\bep\rab^n\hspace{-1.mm}\|\|\nabla\bv_{j,h}^{n+1}\|_{L^3}\|\nabla\be_{j,0}^{n+1}\|\nonumber\\&-2\mu\Delta 
			t\Big((l_{\bhw,h}^n)^2\nabla\be_{j,\bR}^{n+1},\nabla\be_{j,0}^{n+1}\Big)-2\mu\Delta 
			t\Big(\big\{(l_{\bw,h}^n)^2-(l_{\bhw,h}^n)^2\big\}\nabla\bv_{j,h}^{n+1},\nabla\be_{j,0}^{n+1}\Big)-b\Big(\bhw^{'n}_{j,h},\be^n_j,\be_{j,0}^{n+1}\Big)\nonumber\\&-b\Big(\bep^{'n}_{j},\bv_{j,h}^n,\be_{j,0}^{n+1}\Big)-\Big(\frac{\nu_j-\nu_{m,j}}{2}\nabla 
			\bep_{j,0}^n,\nabla \be_{j,0}^{n+1}\Big)-\Big(\frac{\nu_j^{'}+\nu_{m,j}^{'}}{2}\nabla \be_{j,0}^{n},\nabla\be_{j,0}^{n+1}\Big).
		\end{align}

		Using Cauchy-Schwarz inequality, \eqref{eddy-viscosity-term}, and Young's inequalities yields
		\begin{align*}    
			-b\Big(\bhw^{'n}_{j,h},\be^n_j,\be_{j,0}^{n+1}\Big)&=b\Big(\bhw^{'n}_{j,h},\be_{j,0}^{n+1},\be^n_j\Big)\\&\le\|\bhw^{'n}_{j,h}\cdot\nabla\be_{j,0}^{n+1}\|\|\be^n_j\|\\&\le\sqrt{2}\||\bhw^{'n}_{j,h}|\nabla\be_{j,0}^{n+1}\|\|\be^n_j\|\\&\le\sqrt{2}\|l_{\bhw,h}^n\nabla\be_{j,0}^{n+1}\|\|\be^n_j\|\\&\le\|l_{\bhw,h}^n\nabla\be_{j,0}^{n+1}\|^2+\frac{1}{2}\|\be^n_j\|^2.
		\end{align*}
		Using the above bound, triangle inequality, stability estimate, Assumption \ref{assump1}, and finally rearranging, we have
		\begin{align}
			\frac{1}{\Delta t}\Big(\be_j^{n+1}-\be_j^n,\be_{j,0}^{n+1}\Big)+\frac{\Bar{\nu}_{\min}+\Bar{\nu}_{m,\min}}{2}\|\nabla 
			\be_{j,0}^{n+1}\|^2+\left(2\mu\Delta t-1\right)\|l_{\bhw,h}^n\nabla\be_{j,0}^{n+1}\|^2\nonumber\\\le 
			\frac{C}{(\Bar{\nu}_{\min}+\Bar{\nu}_{m,\min})^{\frac12}\Delta 
				t^{\frac12}}\|\nabla\be^{n+1}_{j,\bR}\|\|\nabla\be_{j,0}^{n+1}\|+CC_*\|\hspace{-1.mm}\lab\bep\rab^n\hspace{-1.mm}\|\|\nabla\be_{j,0}^{n+1}\|\nonumber\\+2\mu\Delta 
			t\Big|\Big((l_{\bhw,h}^n)^2\nabla\be_{j,\bR}^{n+1},\nabla\be_{j,0}^{n+1}\Big)\Big|+2\mu\Delta 
			t\Big|\Big(\big\{(l_{\bw,h}^n)^2-(l_{\bhw,h}^n)^2\big\}\nabla\bv_{j,h}^{n+1},\nabla\be_{j,0}^{n+1}\Big)\Big|\nonumber\\+\frac{1}{2}\|\be^n_{j}\|^2+\Big|b\Big(\bep^{'n}_{j},\bv_{j,h}^n,\be_{j,0}^{n+1}\Big)\Big|+\Big|\Big(\frac{\nu_j-\nu_{m,j}}{2}\nabla 
			\bep_{j,0}^n,\nabla \be_{j,0}^{n+1}\Big)\Big|+\Big|\Big(\frac{\nu_j^{'}+\nu_{m,j}^{'}}{2}\nabla 
			\be_{j,0}^{n},\nabla\be_{j,0}^{n+1}\Big)\Big|.\label{before-time-derivative-estimate}
		\end{align}
		To evaluate the discrete time-derivative term, we use polarization identity, Cauchy-Schwarz, Young's, and Poincar\'e's inequalities 
		\begin{align*}
			\frac{1}{\Delta t}\Big(\be_j^{n+1}-\be_j^n,\be_{j,0}^{n+1}\Big)&=\frac{1}{\Delta 
				t}\Big(\be_j^{n+1}-\be_j^n,\be_{j}^{n+1}-\be_{j,\bR}^{n+1}\Big)\\&=\frac{1}{2\Delta 
				t}\Big(\|\be_j^{n+1}-\be_j^n\|^2+\|\be_j^{n+1}\|^2-\|\be_j^n\|^2\Big)-\frac{1}{\Delta t}\Big(\be_j^{n+1}-\be_j^n,\be_{j,\bR}^{n+1}\Big)\\&\ge 
			\frac{1}{2\Delta t}\Big(\|\be_j^{n+1}\|^2-\|\be_j^n\|^2\Big)-\frac{C}{\Delta t}\|\nabla\be_{j,\bR}^{n+1}\|^2.
		\end{align*}
		Plugging the above estimate into \eqref{before-time-derivative-estimate} and using H\"older's, and Young's inequalities, yields
		\begin{align}
			\frac{1}{2\Delta t}\Big(\|\be_j^{n+1}\|^2-\|\be_j^n\|^2\Big)+\frac{\Bar{\nu}_{\min}+\Bar{\nu}_{m,\min}}{2}\|\nabla 
			\be_{j,0}^{n+1}\|^2+\left(2\mu\Delta t-1\right)\|l_{\bhw,h}^n\nabla\be_{j,0}^{n+1}\|^2\nonumber\\\le \frac{C}{\Delta 
				t}\left(\frac{1}{\alpha_j^2}+1\right)\|\nabla\be^{n+1}_{j,\bR}\|^2+\frac{CC_*^2}{\alpha_j}\|\hspace{-1.mm}\lab\bep\rab^n\hspace{-1.mm}\|^2+2\mu\Delta 
			t\Big|\Big((l_{\bhw,h}^n)^2\nabla\be_{j,\bR}^{n+1},\nabla\be_{j,0}^{n+1}\Big)\Big|\nonumber\\+2\mu\Delta 
			t\Big|\Big(\big\{(l_{\bw,h}^n)^2-(l_{\bhw,h}^n)^2\big\}\nabla\bv_{j,h}^{n+1},\nabla\be_{j,0}^{n+1}\Big)\Big|+\frac{\alpha_j}{12}\|\nabla\be_{j,0}^{n+1}\|^2+\frac{1}{2}\|\be^n_{j}\|^2+\Big|b\Big(\bep^{'n}_{j},\bv_{j,h}^n,\be_{j,0}^{n+1}\Big)\Big|\nonumber\\+\frac{\|\nu_j-\nu_{m,j}\|_\infty}{4}\Big(\|\nabla 
			\bep_{j,0}^n\|^2+\|\nabla \be_{j,0}^{n+1}\|^2\Big)+\frac{\|\nu_j^{'}+\nu_{m,j}^{'}\|_\infty}{4}\Big(\|\nabla 
			\be_{j,0}^{n}\|^2+\|\nabla\be_{j,0}^{n+1}\|^2\Big).\label{before-non-linear-bounds}
		\end{align}
		We now find the bounds for the non-linear terms. For the first non-linear term, we use Cauchy-Schwarz, and  Young's inequalities, uniform 
		boundedness in Assumption \ref{assump1}, and the stability estimate, to obtain
		\begin{align*}
			2\mu\Delta t\Big|\Big((l_{\bhw,h}^n)^2\nabla\be_{j,\bR}^{n+1},\nabla\be_{j,0}^{n+1}\Big)\Big|&\le 2\mu\Delta 
			t\|l_{\bhw,h}^n\nabla\be_{j,\bR}^{n+1}\|\|l_{\bhw,h}^n\nabla\be_{j,0}^{n+1}\|\\
			&\le\mu\Delta t\|l_{\bhw,h}^n\nabla\be_{j,\bR}^{n+1}\|^2+\mu\Delta t\|l_{\bhw,h}^n\nabla\be_{j,0}^{n+1}\|^2\\
			&\le\mu\Delta t\|l_{\bhw,h}^n\|_\infty^2\|\nabla\be_{j,\bR}^{n+1}\|^2+\mu\Delta t\|l_{\bhw,h}^n\nabla\be_{j,0}^{n+1}\|^2\\
			&\le C\Delta t\|\nabla\be_{j,\bR}^{n+1}\|^2+\mu\Delta t\|l_{\bhw,h}^n\nabla\be_{j,0}^{n+1}\|^2.
		\end{align*}
		For the second non-linear term, we follow the same treatment as in \eqref{mixing-length-difference}, and get
		\begin{align*}
			2\mu\Delta 
			t\Big(\big\{(l_{\bw,h}^n)^2-(l_{\bhw,h}^n)^2\big\}\nabla\bv_{j,h}^{n+1},\nabla\be_{j,0}^{n+1}\Big)\le\frac{\alpha_j}{12}\|\nabla\be_{j,0}^{n+1}\|^2+\frac{C\Delta 
				t}{h^3\alpha_j}\sum_{i=1}^{N_{sc}}\|\bep_i^n\|^2.
		\end{align*}
		For the last non-linear term,  we apply H\"older's inequality, estimate in Lemma \ref{lemma1}, Sobolev embedding theorem, Poincar\'e  and Young's 
		inequalities to get
		\begin{align*}
			|b\Big(\bep^{'n}_{j},\bv_{j,h}^n,\be_{j,0}^{n+1}\Big)|&\le \|\bep^{'n}_{j}\|\|\nabla\bv_{j,h}^n\|_{L^3}\|\|\be_{j,0}^{n+1}\|_{L^6}\\&\le 
			CC_*\|\bep^{'n}_{j}\|\|\nabla\be_{j,0}^{n+1}\|\\&\le\frac{\alpha_j}{12}\|\nabla\be_{j,0}^{n+1}\|^2+\frac{CC_*^2}{\alpha_j}\|\bep^{'n}_{j}\|^2.
		\end{align*}
		Use the above estimates, assume $\mu>\max\{1,\frac{1}{2\Delta t}\}$ to drop non-negative terms from left-hand-side, and reducing, the equation 
		\eqref{before-non-linear-bounds} becomes
		\begin{align}
			\frac{1}{2\Delta t}\Big(\|\be_j^{n+1}\|^2-\|\be_j^n\|^2\Big)+\frac{\Bar{\nu}_{\min}+\Bar{\nu}_{m,\min}}{4}\|\nabla 
			\be_{j,0}^{n+1}\|^2\nonumber\\\le \frac{C}{\Delta t}\left(\frac{1}{\alpha_j^2}+1+\Delta 
			t^2\right)\|\nabla\be^{n+1}_{j,\bR}\|^2+\frac{CC_*^2}{\alpha_j}\Big(\|\hspace{-1.mm}\lab\bep\rab^n\hspace{-1.mm}\|^2+\|\bep^{'n}_{j}\|^2\Big)\nonumber\\+\frac{C\Delta 
				t}{h^3\alpha_j}\sum_{i=1}^{N_{sc}}\|\bep_i^n\|^2+\frac{1}{2}\|\be^n_{j}\|^2+\frac{\|\nu_j-\nu_{m,j}\|_\infty}{4}\|\nabla 
			\bep_{j,0}^n\|^2+\frac{\|\nu_j^{'}+\nu_{m,j}^{'}\|_\infty}{4}\|\nabla \be_{j,0}^{n}\|^2.\label{e-equ-bounded}
		\end{align}
		Apply similar techniques to \eqref{step-2-eqn-2}, yields
		\begin{align}
			\frac{1}{2\Delta t}\Big(\|\bep_j^{n+1}\|^2-\|\bep_j^n\|^2\Big)+\frac{\Bar{\nu}_{\min}+\Bar{\nu}_{m,\min}}{4}\|\nabla 
			\bep_{j,0}^{n+1}\|^2\nonumber\\\le \frac{C}{\Delta t}\left(\frac{1}{\alpha_j^2}+1+\Delta 
			t^2\right)\|\nabla\bep^{n+1}_{j,\bR}\|^2+\frac{CC_*^2}{\alpha_j}\Big(\|\hspace{-1.mm}\lab\be\rab^n\hspace{-1.mm}\|^2+\|\be^{'n}_{j}\|^2\Big)\nonumber\\+\frac{C\Delta 
				t}{h^3\alpha_j}\sum_{i=1}^{N_{sc}}\|\be_i^n\|^2+\frac{1}{2}\|\bep^n_{j}\|^2+\frac{\|\nu_j-\nu_{m,j}\|_\infty}{4}\|\nabla 
			\be_{j,0}^n\|^2+\frac{\|\nu_j^{'}+\nu_{m,j}^{'}\|_\infty}{4}\|\nabla \bep_{j,0}^{n}\|^2.\label{ep-equ-bounded}
		\end{align}
		Add equations \eqref{e-equ-bounded}, and \eqref{ep-equ-bounded}, and use triangle inequality, to get 
		\begin{align}
			\frac{1}{2\Delta 
				t}\Big(\|\be_j^{n+1}\|^2-\|\be_j^n\|^2+\|\bep_j^{n+1}\|^2-\|\bep_j^n\|^2\Big)+\frac{\Bar{\nu}_{\min}+\Bar{\nu}_{m,\min}}{4}\Big(\|\nabla 
			\be_{j,0}^{n+1}\|^2+\|\nabla \bep_{j,0}^{n+1}\|^2\Big)\nonumber\\\le \frac{C}{\Delta t}\left(\frac{1}{\alpha_j^2}+1+\Delta 
			t^2\right)\lp\|\nabla\be^{n+1}_{j,\bR}\|^2+\|\nabla\bep^{n+1}_{j,\bR}\|^2\rp+\left(\frac{CC_*^2}{\alpha_j}+\frac{C\Delta 
				t}{h^3\alpha_j}+\frac{1}{2}\right)\sum_{j=1}^{N_{sc}}\Big(\|\be_j^n\|^2+\|\bep_j^n\|^2\Big)\nonumber\\+\frac{\|\nu_j-\nu_{m,j}\|_\infty+\|\nu_j^{'}+\nu_{m,j}^{'}\|_\infty}{4}\Big(\|\nabla 
			\be_{j,0}^{n}\|^2+\|\nabla \bep_{j,0}^{n}\|^2\Big).
		\end{align}
		Rearranging
		\begin{align}
			\frac{1}{2\Delta 
				t}\Big(\|\be_j^{n+1}\|^2-\|\be_j^n\|^2+\|\bep_j^{n+1}\|^2-\|\bep_j^n\|^2\Big)\nonumber\\+\frac{\Bar{\nu}_{\min}+\Bar{\nu}_{m,\min}}{4}\Big(\|\nabla 
			\be_{j,0}^{n+1}\|^2-\|\nabla \be_{j,0}^{n}\|^2+\|\nabla \bep_{j,0}^{n+1}\|^2-\|\nabla 
			\bep_{j,0}^{n}\|^2\Big)\nonumber\\+\frac{\alpha_j}{4}\Big(\|\nabla \be_{j,0}^{n}\|^2+\|\nabla \bep_{j,0}^{n}\|^2\Big)\le \frac{C}{\Delta 
				t}\left(\frac{1}{\alpha_j^2}+1+\Delta 
			t^2\right)\lp\|\nabla\be^{n+1}_{j,\bR}\|^2+\|\nabla\bep^{n+1}_{j,\bR}\|^2\rp\nonumber\\+\left(\frac{CC_*^2}{\alpha_{j}}+\frac{C\Delta 
				t}{h^3\alpha_j}+\frac{1}{2}\right)\sum_{j=1}^{N_{sc}}\Big(\|\be_j^n\|^2+\|\bep_j^n\|^2\Big).
		\end{align}
		Multiply both sides by $2\Delta t$, and summing over the time-step $n=0,1,\cdots,M-1$, results in
		\begin{align} \|\be_j^M\|^2+\|\bep_j^M\|^2&+\frac{\Bar{\nu}_{\min}+\Bar{\nu}_{m,\min}}{2}\Delta 
			t\Big(\|\nabla\be_{j,0}^M\|^2+\|\nabla\bep_{j,0}^M\|^2\Big)+\frac{\alpha_j}{2}\Delta t\sum_{n=1}^{M-1}\Big(\|\nabla \be_{j,0}^{n}\|^2+\|\nabla 
			\bep_{j,0}^{n}\|^2\Big)\nonumber\\&\le C\left(\frac{1}{\alpha_j^2}+1+\Delta 
			t^2\right)\sum_{n=1}^M\Big(\|\nabla\be^{n}_{j,\bR}\|^2+\|\nabla\bep^{n}_{j,\bR}\|^2\Big)\nonumber\\&+\Delta 
			t\left(\frac{CC_*^2}{\alpha_j}+\frac{C\Delta t}{h^3\alpha_j}+1\right)\sum_{n=1}^{M-1}\sum_{j=1}^{N_{sc}}\Big(\|\be_j^n\|^2+\|\bep_j^n\|^2\Big).
		\end{align}
		Now, simplifying, and summing over $j=1,2,\cdots,{N_{sc}}$, we have
		\begin{align}
			\sum_{j=1}^{N_{sc}}\Big(\|\be_j^M\|^2+\|\bep_j^M\|^2\Big)+\Delta t\sum_{n=1}^{M}\frac{\alpha_{\min}}{2}\sum_{j=1}^{N_{sc}}\Big(\|\nabla 
			\be_{j,0}^{n}\|^2+\|\nabla \bep_{j,0}^{n}\|^2\Big)\nonumber\\\le  \sum_{n=1}^MC\left(\frac{1}{\alpha_{\min}^2}+1+\Delta 
			t^2\right)\sum_{j=1}^{N_{sc}}\Big(\|\nabla\be^{n}_{j,\bR}\|^2+\|\nabla\bep^{n}_{j,\bR}\|^2\Big)\nonumber\\+\Delta 
			t\sum_{n=1}^{M-1}\left(\frac{CC_*^2}{\alpha_{\min}}+\frac{C\Delta 
				t}{h^3\alpha_{\min}}+{N_{sc}}\right)\sum_{j=1}^{N_{sc}}\Big(\|\be_j^n\|^2+\|\bep_j^n\|^2\Big).
		\end{align}
		Apply the version of the discrete Gr\"onwall inequality given in Lemma \ref{dgl}
		\begin{align}
			\sum_{j=1}^{N_{sc}}\Big(\|\be_j^M\|^2+\|\bep_j^M\|^2\Big)+\frac{\alpha_{\min}}{2}\Delta t\sum_{n=1}^{M}\sum_{j=1}^{N_{sc}}\Big(\|\nabla 
			\be_{j,0}^{n}\|^2+\|\nabla \bep_{j,0}^{n}\|^2\Big)\nonumber\\\le exp \left(\frac{CC_*^2}{\alpha_{\min}}+\frac{C\Delta 
				t}{h^3\alpha_{\min}}+{N_{sc}}*T\right)\Bigg[ C\left(\frac{1}{\alpha_{\min}^2}+1+\Delta 
			t^2\right)\sum_{n=1}^M\sum_{j=1}^{N_{sc}}\Big(\|\nabla\be^{n}_{j,\bR}\|^2+\|\nabla\bep^{n}_{j,\bR}\|^2\Big)\Bigg],\label{step-2-after-gronwell}
		\end{align}
		and use the estimate  \eqref{step-1-bound} in \eqref{step-2-after-gronwell}, to get
		\begin{align}
			\Delta t\sum_{n=1}^{M}\sum_{j=1}^{N_{sc}}\Big(\|\nabla \be_{j,0}^{n}\|^2+\|\nabla \bep_{j,0}^{n}\|^2\Big)\le 
			\frac{CC_R^2}{\gamma^2\alpha_{\min}} exp \left(\frac{CC_*^2}{\alpha_{\min}}+\frac{C\Delta t}{h^3\alpha_{\min}}\right)\nonumber\\\times\Bigg[ 
			\left(\frac{1}{\alpha_{\min}^2}+1+\Delta 
			t^2\right)\sum_{n=0}^{M-1}\sum_{j=1}^{N_{sc}}\Big(\|q_{j,h}^{n+1}-\hq_{j,h}^n\|^2+\|r_{j,h}^{n+1}-\hr_{j,h}^n\|^2\Big)\Bigg].
		\end{align}
		Using triangle and Young's inequalities
		\begin{align}
			\Delta t\sum_{n=1}^{M}\left(\|\nabla \hspace{-1.mm}\lab\be_{0}\rab^{n}\|^2+\|\nabla\hspace{-1mm} 
			\lab\bep_{0}\rab^{n}\|^2\right)\le\frac{2\Delta t}{{N^2_{sc}}}\sum_{n=1}^{M}\sum_{j=1}^{N_{sc}}\left(\|\nabla \be_{j,0}^{n}\|^2+\|\nabla 
			\bep_{j,0}^{n}\|^2\right)\nonumber\\\le \frac{CC_R^2}{\gamma^2\alpha_{\min}} exp \left(\frac{CC_*^2}{\alpha_{\min}}+\frac{C\Delta 
				t}{h^3\alpha_{\min}}\right)\nonumber\\\times\Bigg[ \left(\frac{1}{\alpha_{\min}^2}+1+\Delta 
			t^2\right)\sum_{n=0}^{M-1}\sum_{j=1}^{N_{sc}}\Big(\|q_{j,h}^{n+1}-\hq_{j,h}^n\|^2+\|r_{j,h}^{n+1}-\hr_{j,h}^n\|^2\Big)\Bigg],\label{step-2-bound-ensemble}
		\end{align}
		and
		\begin{align}
			\Delta 
			t\sum_{n=1}^{M}\Big(\|\nabla\hspace{-1mm}\lab\be_{\bR}\rab^{n}\|^2+\|\nabla\hspace{-1mm}\lab\bep_{\bR}\rab^{n}\|^2\Big)\le\frac{2\Delta 
				t}{{N_{sc}^2}}\sum_{n=1}^{M}\sum_{j=1}^{N_{sc}}\Big(\|\nabla\be_{j,\bR}^{n}\|^2+\|\nabla\bep_{j,\bR}^{n}\|^2\Big)\nonumber\\\le\frac{CC_R^2}{\gamma^2} 
			exp\lp \frac{CC_*^2}{\alpha_{\min}}+\frac{C\Delta t}{h^3\alpha_{\min}}\rp\lp  \Delta 
			t\sum_{n=0}^{M-1}\sum_{j=1}^{N_{sc}}\Big(\|q_{j,h}^{n+1}-\hq_{j,h}^n\|^2+\|r_{j,h}^{n+1}-\hr_{j,h}^n\|^2\Big)\rp.
		\end{align}
		Finally, apply triangle and Young's inequalities on 
		$$\|\nabla\hspace{-1mm}\lab\bv_h\rab^n-\nabla\hspace{-1mm}\lab\bhv_h\rab^n\|^2+\|\nabla\hspace{-1mm}\lab\bw_h\rab^n-\nabla\hspace{-1mm}\lab\bhw_h\rab^n\|^2$$ 
		to obtain the desire result.
	\end{proof}

	We prove the following Lemma by strong mathematical induction.
	\begin{lemma}\label{uniform-boundedness-lemma-proof}
		If $\gamma\rightarrow\infty$ then there exists a constant $C_*$ which is independent of $h$, and $\Delta t$, such that for sufficiently small $h$ 
		for a fixed mesh and fixed $\Delta t$, the solution of the Algorithm \ref{SPP-FEM} satisfies
		\begin{align}
			\max_{0\le n\le M}\Big\{\|\bhv_{j,h}^n\|_\infty,\|\bhw_{j,h}^n\|_\infty\Big\}&\le C_*,\hspace{2mm}\forall j=1,2,\cdots,{N_{sc}}.
		\end{align}
	\end{lemma}
	\begin{proof}
		Basic step: $\bhv_{j,h}^0=I_h(\bv_j^{true}(0,\bx)),$ where $I_h$ is an appropriate interpolation operator. Because of the regularity assumption 
		of $\bv_j^{true}(0,\bx)$, we have $\|\bhv_{j,h}^0\|_\infty\le C_*$, for some constant $C_*>0$.\\
		Inductive step: Assume for some $L\in\mathbb{N}$ and $L<M$, $\|\bhv_{j,h}^n\|_\infty\le C_*$ holds true for $n=0,1,\cdots,L$. Then, using 
		triangle inequality and Lemma \ref{lemma1}, we have
		\begin{align*}
			\|\bhv_{j,h}^{L+1}\|_{\infty}\le\|\bhv_{j,h}^{L+1}-\bv_{j,h}^{L+1}\|_{\infty}+C_*.
		\end{align*}
		Using Agmon’s inequality \cite{Robinson2016Three-Dimensional}, and discrete inverse inequality, yields
		\begin{align}
			\|\bhv_{j,h}^{L+1}\|_{\infty}\le Ch^{-\frac32}\|\bhv_{j,h}^{L+1}-\bv_{j,h}^{L+1}\|+C_*.
		\end{align}
		Next, using equation \eqref{after-gronwall}
		\begin{align}
			&\|\bhv_{j,h}^{L+1}\|_{\infty}\le C_*\nonumber\\&+\frac{C}{h^{\frac32}\gamma^{\frac12} } exp\lp 
			CT\left(\frac{C_*^2}{\alpha_{\min}}+\frac{\Delta t}{h^3\alpha_{\min}}\right)\rp\lp  \Delta 
			t\sum_{n=0}^{L}\sum_{j=1}^{N_{sc}}\Big(\|q_{j,h}^{n+1}-\hq_{j,h}^n\|^2+\|\lambda_{j,h}^{n+1}-\hlam_{j,h}^n\|^2\Big)\rp^{\frac12}.
		\end{align}
		For a fixed mesh, and timestep size, as $\gamma\rightarrow \infty$, yields $\|\bhv_{j,h}^{L+1}\|_{\infty}\le C_*$. Hence, by the principle of 
		strong mathematical induction, $\|\bhv_{j,h}^{n}\|_{\infty}\le C_*$ holds true for $0\le n\le M$.
		
		Similarly, we can prove the uniform boundedness of $\bhw^n_{j,h}$.
	\end{proof}
	\section{SCMs}\label{scm}
	As SCMs, in this work, we consider sparse grid algorithm \cite{FTW2008}, where for a given time $t_n$ and a set of sample points $\{\by^j\}_{j=0}^{N_{sc}}\subset\bGamma$, we approximate the exact solution of \eqref{weak_formulation_final-start}-\eqref{weak_formulation_final-end} by solving a discrete scheme. Then, for a basis $\{\phi_l\}_{l=1}^{N_p}$ of dimension $N_p$ for the space $L_\rho^2(\bGamma)$, a discrete approximation is constructed with coefficients $c_l(t_n,\bx)$ of the form \begin{align*}
		\bu_h^{sc}(t_n,\bx,\by)=\sum_{l=1}^{N_{p}}c_l(t_n,\bx)\phi_l(\by),
	\end{align*}
	which is essentially an interpolant. In the sparse grid algorithm, we consider Leja and Clenshaw--Curtis points as the interpolation points that come with the associated weights $\{w^j\}_{j=1}^{N_{sc}}$. SCMs were developed for the UQ of the Quantity of Interest (QoI), $\psi$, which can be lift, drag, and energy. SCMs provide statistical information about QoI, that is, $$\mathbb{E}[\psi(\bu)]=\int_{\Gamma} \psi(\bu,\by)\rho(\by)dy\approx\sum_{j=1}^{N_{sc}}w^j\psi(\bu,\by^j).$$
	\section{Numerical Experiments}
	\label{numerical-exp}
	To test the proposed Algorithm \ref{SPP-FEM} (SCM-SPP-SMHD method) and the associated theory, in this section, we present the results of 
	numerical experiments. For MHD simulations, it is crucial to enforce the discrete solenoidal constraint $\nabla\cdot\bB_{j,h}=0$ strongly, 
	otherwise, it can produce large errors in the solution \cite{case2010high, linke2009collision}. Moreover, to have the divergence-free condition 
	of the magnetic field at all times, the initial magnetic field must need to be zero. This is because the curl of the electric field is equal to 
	and opposite of the time derivative of the magnetic flux density. Thus, it is popular to use pointwise divergence-free elements such as 
	Scott-Vogelius (SV) elements on barycenter refined regular
	triangular meshes to enforce the divergence constraints \cite{AKMR15,HMR17, kuberry2012numerical,MR17,mohebujjaman2022efficient}. However, using 
	SV elements require higher degrees of freedom (dof) which is quite demanding. Throughout this numerical section, we will use $(P_2,P_1^{disc})$ 
	SV element in the Coupled-SMHD method for the velocity-pressure and magnetic flux density-magnetic pressure variables and their outcomes will be 
	considered as the benchmark solutions. Also, in the SCM-SPP-SMHD method, we will use $(P_2, P_1)$ Taylor-hood (TH) element (which is weakly 
	divergence-free and requires less dof) with large a $\gamma$. Both methods will be employed on a barycenter refined triangular mesh.
	
	In the first experiment, we verify the predicted convergence rates given in Theorem \ref{gamma-convergence} as $\gamma$ varies and compute the 
	spatial and temporal convergence rates with manufactured solutions. We implement the scheme on a channel flow over a step problem and a 
	regularized lid-driven cavity problem in the second and third experiments, respectively. Finally, we examine the sparse grid algorithm as SCM in the lid-driven cavity problem.

	\subsection{Convergence rate verification}
	
	We will begin this experiment with $\bx=(x_1,x_2)$ and the following manufactured analytical functions,\[
	{\bv}=\left(\begin{array}{c} \cos x_2+(1+e^t)\sin x_2 \\ \sin x_1+(1+e^t)\cos x_1 \end{array} \right), \
	{\bw}=\left(\begin{array}{c} \cos x_2-(1+e^t)\sin x_2 \\ \sin x_1-(1+e^t)\cos x_1 \end{array} \right), \ q =\sin(x_1+x_2)(1+e^t), \ r=0.
	\]
	
	Clearly, $\nabla\cdot\bv=\nabla\cdot\bw=0$. Next, introducing a perturbation parameter $\epsilon$ we introduce noise in the above analytical 
	functions as below to create manufactured solutions
	\begin{align}
		\bv_j(t,\bx):=\left(1+k_j\epsilon\right)\bv,\hspace{1mm} \bw_j(t,\bx):=\left(1+k_j\epsilon\right)\bw,\hspace{1mm} 
		q_j:=(1+k_j\epsilon)q,\hspace{1mm}\text{and}\hspace{1mm}r_j:=0,\label{manufactured-solution}
	\end{align}
	where $k_j:=\frac{(-1)^{j+1}4\lceil j/2\rceil}{N_{sc}}$, and $j=1,2,\cdots, N_{sc}$, where $N_{sc}=20$. We consider the kinematic viscosity $\nu$ 
	and magnetic diffusivity $\nu_m$ are  continuous random variables with uniform distribution. In this experiment, we consider $\nu\sim 
	\mathcal{U}(0.0009, 0.0011)$ with $E[\nu]=0.001$,  $\nu\sim \mathcal{U}(0.009, 0.011)$ with $E[\nu]=0.01$, and $\nu_m\sim\mathcal{U}(0.0009, 
	0.0011)$ with $E[\nu_m]=0.001$. For each of the cases, we collect a i.i.d sample size of 20 which leads us to have two two-dimensional random 
	samples. For a fixed $j$ together with pair $(\nu_j,\nu_{m,j})$, and the analytical solution in \eqref{manufactured-solution}, we compute the 
	forcing functions as below:
	\begin{eqnarray*}
		\bif_{1,j}=\bv_{j,t}+\bw_j\cdot\nabla \bv_j-\frac{\nu_j+\nu_{m,j}}{2}\Delta \bv_j-\frac{\nu_j-\nu_{m,j}}{2}\Delta \bw_j+\nabla q_j,\\
		\bif_{2,j}=\bw_{j,t}+\bv_j\cdot\nabla \bw_j-\frac{\nu_j+\nu_{m,j}}{2}\Delta \bw_j-\frac{\nu_j-\nu_{m,j}}{2}\Delta \bv_j+\nabla r_j.
	\end{eqnarray*}
	We consider a domain $\cD=(0,1)^2$, and $\bv_{j,h}^0=\bv_j(0,\bx)$ and $\bw_{j,h}^0=\bw_j(0,\bx)$ as the initial conditions for both algorithms. 
	The boundary conditions for the Algorithm \ref{Algn0} are considered as $\bv_{j,h}|_{\partial\cD}=\bv_j$ in Step 1, and 
	$\bw_{j,h}|_{\partial\cD}=\bw_j$ in Step 2, whereas the boundary conditions for the Algorithm \ref{SPP-FEM} are considered as 
	$\bhv_{j,h}|_{\partial\cD}=\bv_j$ in Step 1, $\btv\cdot\bnh|_{\partial\cD}=0$ in Step 2, $\bhw_{j,h}|_{\partial\cD}=\bw_j$ in Step 3, and 
	$\btw\cdot\bnh|_{\partial\cD}=0$ in Step 4. We compute the ensemble average solutions 
	$(<\hspace{-1mm}\bv_h\hspace{-1mm}>^n,<\hspace{-1mm}\bw_h\hspace{-1mm}>^n)$, and 
	$(<\hspace{-1mm}\bhv_{h,\gamma}\hspace{-1mm}>^n,<\hspace{-1mm}\bhw_{h,\gamma}\hspace{-1mm}>^n)$ at $t=t^n$ using the Algorithm \ref{Algn0}, and 
	the penalty-projection based Algorithm \ref{SPP-FEM}, respectively, and compare them by computing the difference between the two algorithms.
	\subsubsection{Convergence with $\gamma$ varies} To observe the convergence of the SCM-SPP-SMHD to the Coupled-SMHD scheme, for $\bz=\bv$ or 
	$\bw$, we define $<\hspace{-1mm}\be_{h,\gamma}^{\bz}\hspace{-1mm}>:=<\hspace{-1mm}\bz_{h}\rab-\lab\bhz_{h,\gamma}\hspace{-1mm}>$ and compute
	$$\|\hspace{-1mm}<\hspace{-1mm}\be_{h,\gamma}^{\bz}\hspace{-1mm}>\hspace{-1mm}\|_{2,1}:=\|\hspace{-1mm}<\hspace{-1mm}\bz_{h}\rab-\lab\bhz_{h,\gamma}\hspace{-1mm}>\hspace{-1mm}\|_{L^2(0,T;H^1(\cD)^2)}.$$ 
	\begin{table}[!ht] 
		\begin{center}
			\small\begin{tabular}{|c|c|c|c|c|c|c|c|c|}\hline
				\multicolumn{9}{|c|}{Fixed $T=1.0$, $\Delta t=T/10$, $h=1/32$ with $j=1,2,\cdots\hspace{-0.35mm}, 20$}\\\hline
				$\hspace{-1mm}\epsilon=0.01\hspace{-1mm}$&\multicolumn{4}{|c|}{$\hspace{-1mm}\{(\nu_j,\nu_{m,j})\in[0.009,0.011]\times[0.0009,0.0011]\}\hspace{-1mm}$}&\multicolumn{4}{|c|}{$\hspace{-1mm}\{(\nu_j,\nu_{m,j})\in[0.009,0.011]\times[0.09,0.11]\}\hspace{-1mm}$}\\\hline
				$\gamma$ & $\|\hspace{-1mm}<\hspace{-1mm}\be_{h,\gamma}^{\bv}\hspace{-1mm}>\hspace{-1mm}\|_{2,1}$ & rate   
				&$\|\hspace{-1mm}<\hspace{-1mm}\be_{h,\gamma}^{\bw}\hspace{-1mm}>\hspace{-1mm}\|_{2,1}$  & rate 
				&$\|\hspace{-1mm}<\hspace{-1mm}\be_{h,\gamma}^{\bv}\hspace{-1mm}>\|_{2,1}$ & rate   
				&$\|\hspace{-1mm}<\hspace{-1mm}\be_{h,\gamma}^{\bw}\hspace{-1mm}>\|_{2,1}$  & rate\\ \hline
				1 & 2.5284e-0  & & 1.8966e-0 & & 1.5279e-0  & & 1.4276e-0 & \\ \hline
				10 & 3.7456e-1  & 0.83 & 3.2483e-1 & 0.77 &2.0998e-1  & 0.86 & 1.8832e-1 & 0.88 \\ \hline
				$10^2$ & 3.9501e-2  & 0.98 & 3.5633e-2 & 0.96 &2.1919e-2  & 0.98 & 1.9578e-2 & 0.98 \\ \hline
				$10^3$ & 3.9690e-3  & 1.00 & 3.5951e-3 & 1.00 &2.2018e-3  & 1.00 & 1.9660e-3 & 1.00 \\ \hline
				$10^4$ & 3.9450e-4  & 1.00 & 3.5761e-4 & 1.00 &2.2026e-4  & 1.00 & 1.9668e-4 & 1.00 \\ \hline
				$10^5$ & 3.8554e-5  & 1.01 & 3.6154e-5 & 1.00 &2.2036e-5  & 1.00 & 1.9649e-5 & 1.00 \\ \hline
			\end{tabular}
		\end{center}
		\caption{\footnotesize Convergence of the SCM-SPP-SMHD scheme to the Coupled-SMHD scheme as $\gamma$ 
			increases.}\label{Gamma-convergence-nu-0.01-num-0.001-ep-0.01-wang_projection}
	\end{table}
	In Table \ref{Gamma-convergence-nu-0.01-num-0.001-ep-0.01-wang_projection}, we represent the above errors and convergence rates as $\gamma$ 
	increases with fixed $\epsilon=0.01$, $T=1.0$, $\Delta t=T/10$, $h=1/32$, and two 2D samples 
	$\{(\nu_j,\nu_{m,j})\in[0.009,0.011]\times[0.0009,0.0011]\}$, and $\{(\nu_j,\nu_{m,j})\in[0.009,0.011]\times[0.09,0.11]\}$. We observe a 
	first-order convergence as the $\gamma$ increases, which is in excellent agreement with the Theorem \ref{gamma-convergence}.
	To observe the spatial and temporal errors and their convergence rates, we define 
	$<\hspace{-1mm}\be_{\bz}\hspace{-1mm}>:=<\hspace{-1mm}\bz\rab-\lab\bhz_{h}\hspace{-1mm}>$ for $\bz=\bv$ or $\bw$, which are the difference 
	between the outcomes of the SCM-SPP-SMHD scheme and the true analytical solutions stated above.
	\begin{table}[!ht] 
		\begin{center}
			\small\begin{tabular}{|c|c|c|c|c|c|c|c|c|}\hline
				\multicolumn{9}{|c|}{Temporal convergence (fixed $h=1/64$, $T=1$) with $j=1,2,\cdots\hspace{-0.35mm}, 20$}\\\hline
				$\hspace{-1mm}\epsilon=0.01\hspace{-1mm}$&\multicolumn{4}{|c|}{$\hspace{-1mm}\{(\nu_j,\nu_{m,j})\in[0.0009,0.0011]\times[0.0009,0.0011]\}\hspace{-1mm}$}&\multicolumn{4}{|c|}{\hspace{-1mm}$\{(\nu_j,\nu_{m,j})\in[0.009,0.011]\times[0.0009,0.0011]\}\hspace{-2mm}$}\\\hline
				$\Delta t$ & $\|\hspace{-1mm}<\hspace{-1mm}\be_{\bv}\hspace{-1mm}>\hspace{-1mm}\|_{2,1}$ & rate   
				&$\|\hspace{-1mm}<\hspace{-1mm}\be_{\bw}\hspace{-1mm}>\hspace{-1mm}\|_{2,1}$  & rate & 
				$\|\hspace{-1mm}<\hspace{-1mm}\be_{\bv}\hspace{-1mm}>\hspace{-1mm}\|_{2,1}$ & rate   
				&$\|\hspace{-1mm}<\hspace{-1mm}\be_{\bw}\hspace{-1mm}>\hspace{-1mm}\|_{2,1}$  & rate \\ \hline
				$\frac{T}{4}$ & 9.5972e-1 & & 7.3434e-1 &  & 4.3904e-1  & & 3.2078e-1 & \\ \hline
				$\frac{T}{8}$ & 4.7453e-1 &1.02& 3.7110e-1 & 0.98 & 2.4865e-1  &0.82 & 1.9816e-1 &0.69\\ \hline
				$\frac{T}{16}$ & 2.3095e-1 &1.04& 1.8421e-1 & 1.01 & 1.4170e-1  &0.81 & 1.2082e-1 &0.71\\ \hline
				$\frac{T}{32}$ & 1.1009e-1 &1.07& 9.0537e-2 & 1.02 & 7.7004e-2  &0.88 & 6.8704e-2 &0.81\\ \hline
				$\frac{T}{64}$ & 5.0218e-2 &1.13& 4.4106e-2 & 1.04 & 3.8578e-2  &1.00 & 3.5557e-2 &0.95\\ \hline
			\end{tabular}
		\end{center}
		\caption{\footnotesize Temporal errors and convergence rates of SCM-SPP-SMHD scheme for $\bv$ and $\bw$ with $\epsilon=0.01$, and 
			$\gamma=500$.}\label{temporal-convergence-ep-0.01-gamma500} 
	\end{table}
	
	\begin{table}[!ht] 
		\begin{center}
			\small\begin{tabular}{|c|c|c|c|c|c|c|c|c|}\hline
				\multicolumn{9}{|c|}{Spatial convergence (fixed $T=0.001$, $\Delta t=T/8$) with $j=1,2,\cdots\hspace{-0.35mm}, 20$}\\\hline
				$\hspace{-1mm}\epsilon=0.01\hspace{-1mm}$&\multicolumn{4}{|c|}{$\hspace{-1mm}\{(\nu_j,\nu_{m,j})\in[0.0009,0.0011]\times[0.0009,0.0011]\}\hspace{-1mm}$}&\multicolumn{4}{|c|}{\hspace{-1mm}$\{(\nu_j,\nu_{m,j})\in[0.009,0.011]\times[0.0009,0.0011]\}\hspace{-2mm}$}\\\hline
				$h$ & $\|\hspace{-1mm}<\hspace{-1mm}\be_{\bv}\hspace{-1mm}>\hspace{-1mm}\|_{2,1}$ & rate   
				&$\|\hspace{-1mm}<\hspace{-1mm}\be_{\bw}\hspace{-1mm}>\hspace{-1mm}\|_{2,1}$  & rate & 
				$\|\hspace{-1mm}<\hspace{-1mm}\be_{\bv}\hspace{-1mm}>\hspace{-1mm}\|_{2,1}$ & rate   
				&$\|\hspace{-1mm}<\hspace{-1mm}\be_{\bw}\hspace{-1mm}>\hspace{-1mm}\|_{2,1}$  & rate \\ \hline
				$\frac{1}{4}$ & 1.1422e-4 & & 2.1889e-4 &  & 1.1422e-4  & & 2.1889e-4 & \\ \hline
				$\frac{1}{8}$ & 2.8712e-5 &1.99& 5.4683e-5 & 2.00 & 2.8712e-5  &1.99 & 5.4683e-5 &2.00\\ \hline
				$\frac{1}{16}$ & 7.1882e-6 &2.00& 1.3669e-5 & 2.00 & 7.1882e-6  &2.00 & 1.3669e-5 &2.00\\ \hline
				$\frac{1}{32}$ & 1.8012e-6 &2.00& 3.4189e-6 & 2.00 & 1.8010e-6  &2.00 & 3.4188e-6 &2.00\\ \hline
				$\frac{1}{64}$ & 4.7631e-7 &1.92& 8.6788e-7 & 1.98 & 4.7252e-7  &1.93 & 8.6622e-7 &1.98\\ \hline
			\end{tabular}
		\end{center}
		\caption{\footnotesize Spatial errors and convergence rates of SCM-SPP-SMHD scheme for $\bv$ and $\bw$ with $\epsilon=0.01$, and 
			$\gamma=10^6$.}\label{sp-convergence-ep-0.01-1000000} 
	\end{table}
	To receive the temporal convergence, we use a fixed mesh width of $h=1/64$, end time $T=1$, vary timestep size as $\Delta t=T/4, T/8, T/16, 
	T/32,$ and  $T/64$, on the other hand, to get the spatial convergence, we use a small end time $T=0.001$, a fixed timestep size $\Delta t=T/8$, 
	vary mesh size as $h=1/4,1/8,1/16,1/32,$ and $1/64$. For both cases, we run the simulations using the proposed Algorithm \ref{SPP-FEM} varying 
	the perturbation parameter $\epsilon$ (which introduces noise in the initial, and boundary conditions and forcing functions), and the two 2D 
	random samples $\{(\nu_j,\nu_{m,j})\in[0.0009,0.0011]\times[0.0009,0.0011]\}$ and $\{(\nu_j,\nu_{m,j})\in[0.009,0.011]\times[0.0009,0.0011]\}$. 
	Then, we record the errors, compute the convergence rates, and present them in Tables 
	\ref{temporal-convergence-ep-0.01-gamma500}-\ref{sp-convergence-ep-0.01-1000000}. In Table \ref{temporal-convergence-ep-0.01-gamma500}, we 
	observe the first-order temporal convergence which is the optimal convergence rate of a first-order time-stepping algorithm. In Tables 
	\ref{sp-convergence-ep-0.01-1000000}, we observe a second-order spatial convergence which is also consistent with the theory as we have used 
	$(P_2,P_1)$ element.
	
	\subsection{SMHD channel flow past a unit step: A comparison between SCM-SPP-SMHD and Coupled-SMHD schemes}We now implement the SCM-SPP-SMHD and Coupled-SMHD schemes in a 2D channel of electrically conducting fluid flow past a unit step under the influence of a magnetic field and compare their outcomes. 
	\begin{figure} [ht]
		\centering	{\includegraphics[width=0.49\textwidth,height=0.2\textwidth]{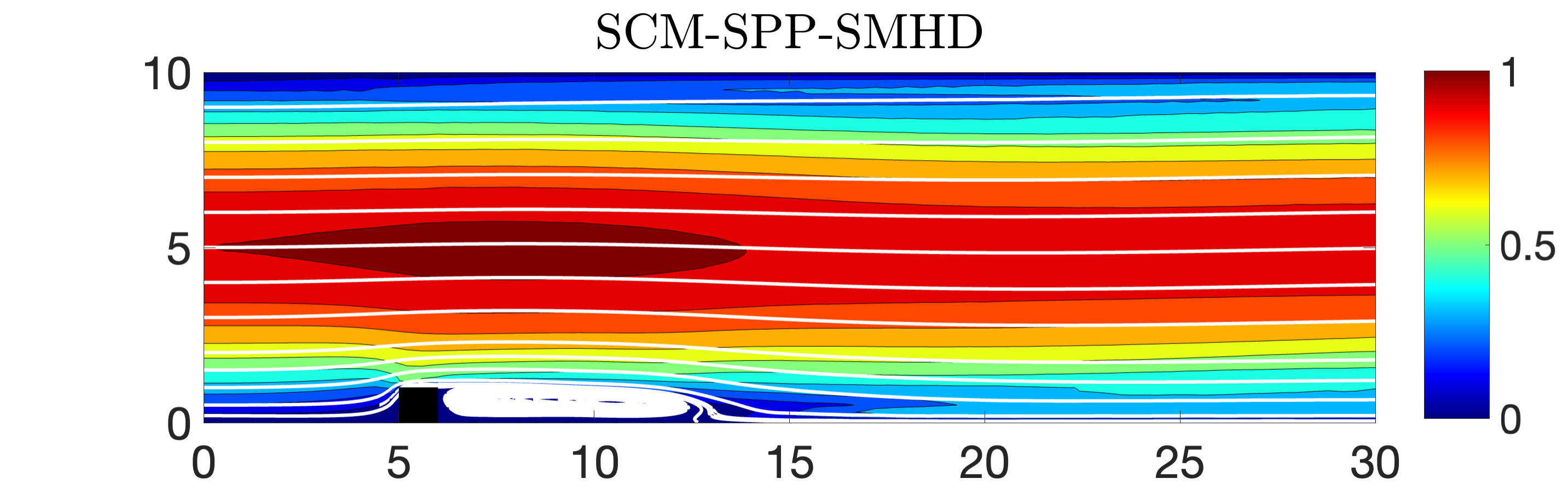}}	{\includegraphics[width=0.49\textwidth,height=0.2\textwidth]{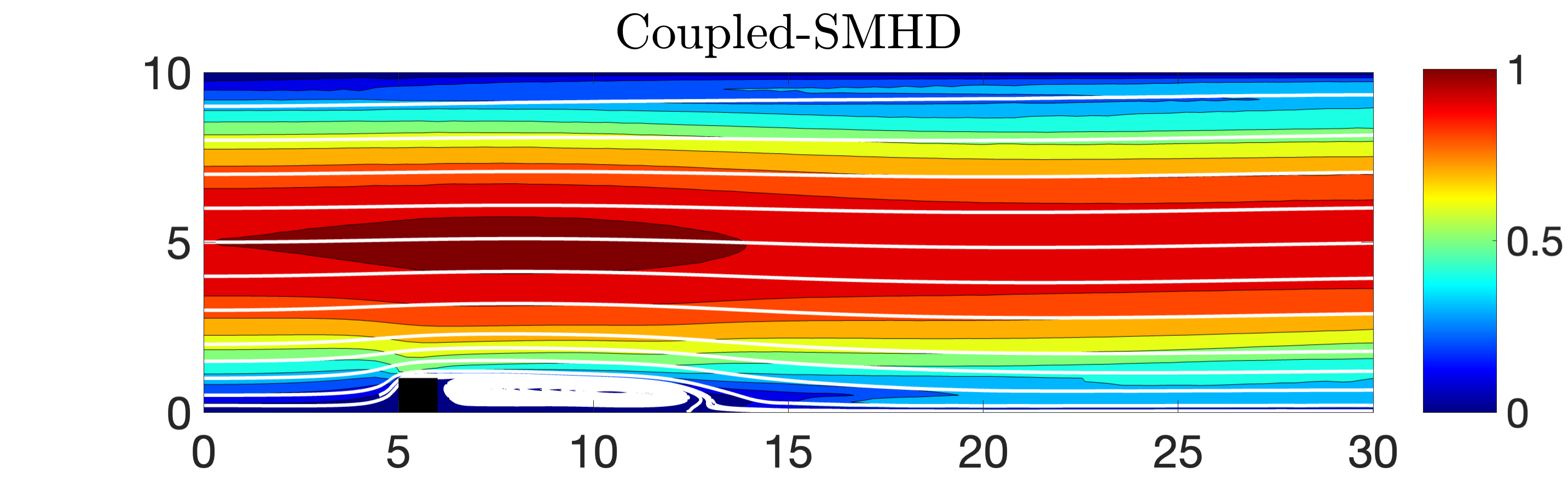}}
		{\includegraphics[width=0.49\textwidth,height=0.2\textwidth]{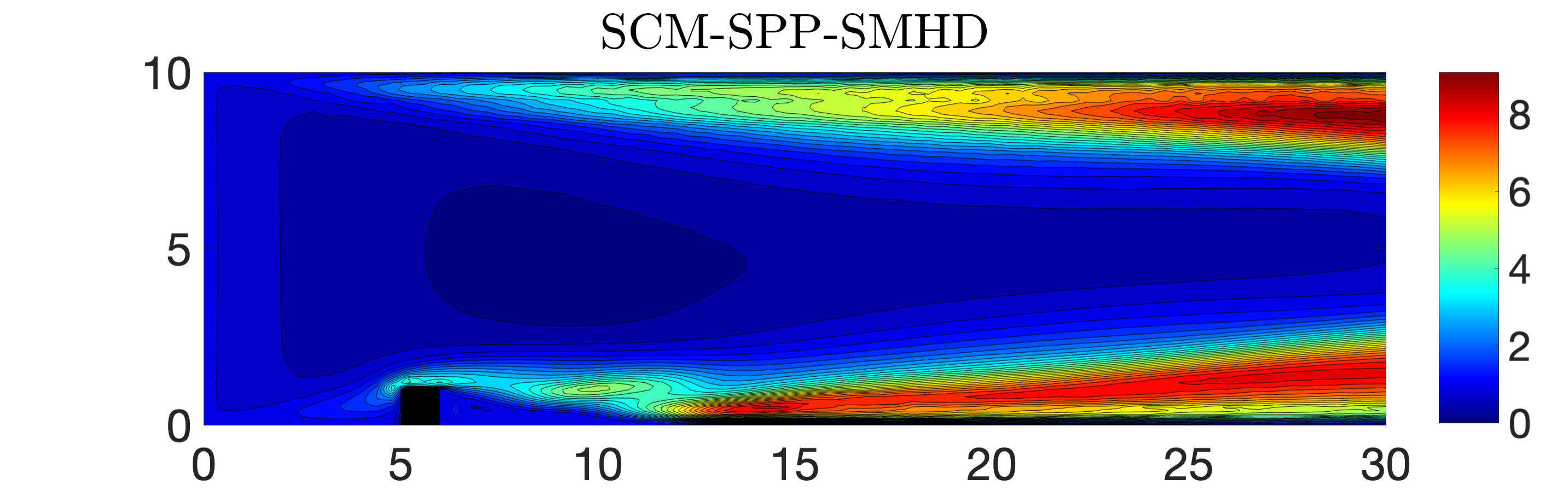}}
		{\includegraphics[width=0.49\textwidth,height=0.2\textwidth]{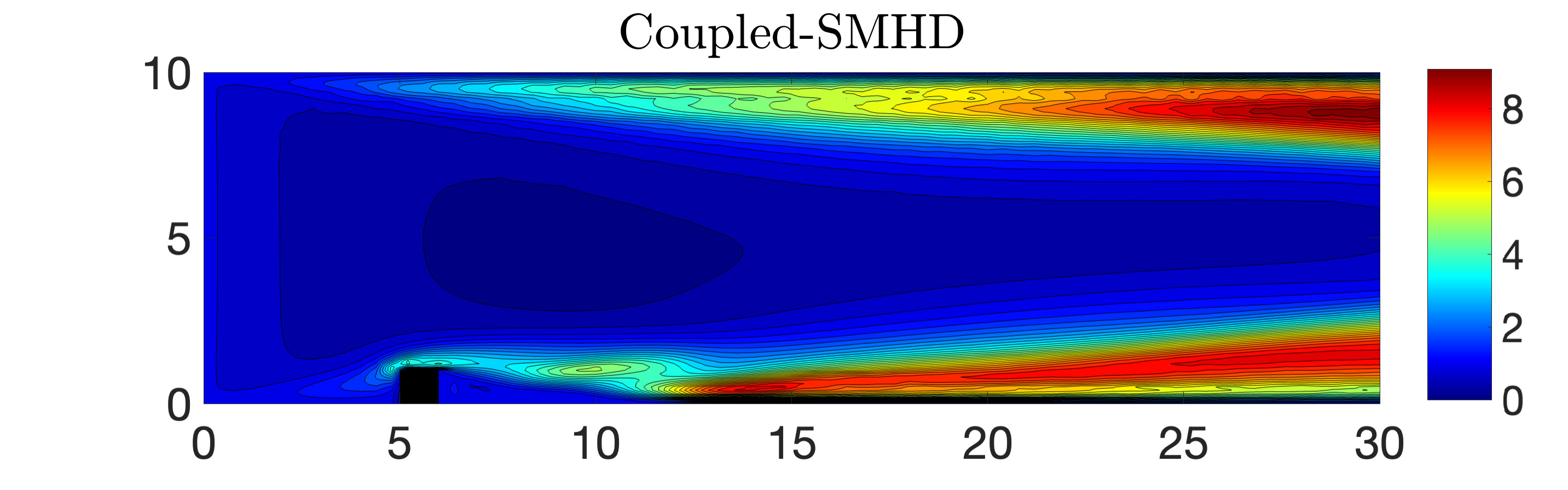}}
		\caption{SCM-SPP-SMHD vs. Coupled-SMHD: Velocity ensemble average solutions shown as streamlines over speed contours (top) and the magnetic field strength (bottom) at $T=40$ with $\epsilon=0.01$, $\mathbb{E}[\nu]=0.001$ and $\mathbb{E}[\nu_{m}]=0.01$.}\label{channel-comparison}
	\end{figure}
	The domain of the flow is a $30\times 10$ rectangular channel over a $1 \times 1$ step on the lower wall which is five units away from the inflow. At the inflow, we set $\bu_j=(1+k_j\epsilon)\begin{pmatrix}
		\frac{x_2(10-x_2)}{25} \\0
	\end{pmatrix}$ and $B_j=\begin{pmatrix}
		0\\1
	\end{pmatrix}$, and the outflow condition uses a channel extension of 10 units, and at the end of the extension, we set outflow velocity and magnetic field equal to their counterpart in the inflow. We consider the initial conditions as \begin{align*}
		\bu_j^0=(1+k_j\epsilon)\begin{pmatrix}\frac{x_2(10-x_2)}{25}\\ 0\end{pmatrix}\hspace{-0.8mm},\hspace{1mm}\text{and}\hspace{1mm} 
		\bB_j^0=\begin{pmatrix}0\\0\end{pmatrix}.
	\end{align*}For the Coupled-SMHD scheme, on the walls, we implement no-slip boundary conditions $\bu_j|_{\Gamma_1}=\begin{pmatrix}0\\0\end{pmatrix}$ and 
	$\bB_j|_{\Gamma_1}=(1+k_j\epsilon)\begin{pmatrix}0\\1\end{pmatrix}$. For the SCM-SPP-SMHD scheme, since the velocity and pressure-like variables appear in different steps, in Step 1, and Step 3, we consider $\bu_j|_{\Gamma_1}=\begin{pmatrix}0\\0\end{pmatrix}$ and $\bB_j|_{\Gamma_1}=(1+k_j\epsilon)\begin{pmatrix}0\\1\end{pmatrix}$ on the walls, and in Step 3, and Step 4, we define the following space for the Els\"asser variables $\bv$, and $\bw$:
	$$\tilde{\bY}_h:=\{\bv_h\in\mathcal{P}_k(\tau_h)^d\cap\bH^1 (\cD)^d:\bv_h\cdot\bnh|_{\Gamma_1}=0\}.$$
	The timestep size $\Delta t=0.05$, $N_{sc}=20$, $\gamma=10^5$, $\bif=\bg=\textbf{0}$ (no-external source), and a constant coupling parameter $s=0.001$ are considered. We consider the mean kinematic viscosity $E[\nu]=0.001$, and mean magnetic diffusivity $E[\nu_{m}]=0.01$ for random samples with distribution $\nu\sim\mathcal{U}(0.0009,0.0011)$, and $\nu_m\sim\mathcal{U}(0.009,0.011)$, respectively. A triangular unstructured mesh of the domain that provides a total of $419058$ dof is considered, where velocity dof $=186134$, magnetic field dof $=186134$, pressure dof 
	$=23395$, and magnetic pressure dof $=23395$. We run the simulations until $T=40$ and plot the speed contour and magnetic field strength in 
	Fig. \ref{channel-comparison} for both SCM-SPP-SMHD and Coupled-SMHD algorithms. We observe a very good agreement between the solutions of the two algorithms which supports our claim in the theory.
	
	\subsection{Variable 5D Random Viscosities with regularized lid-driven cavity problem}
\begin{figure} [ht]
		\centering
		{\includegraphics[width=0.28\textwidth,height=0.2\textwidth]{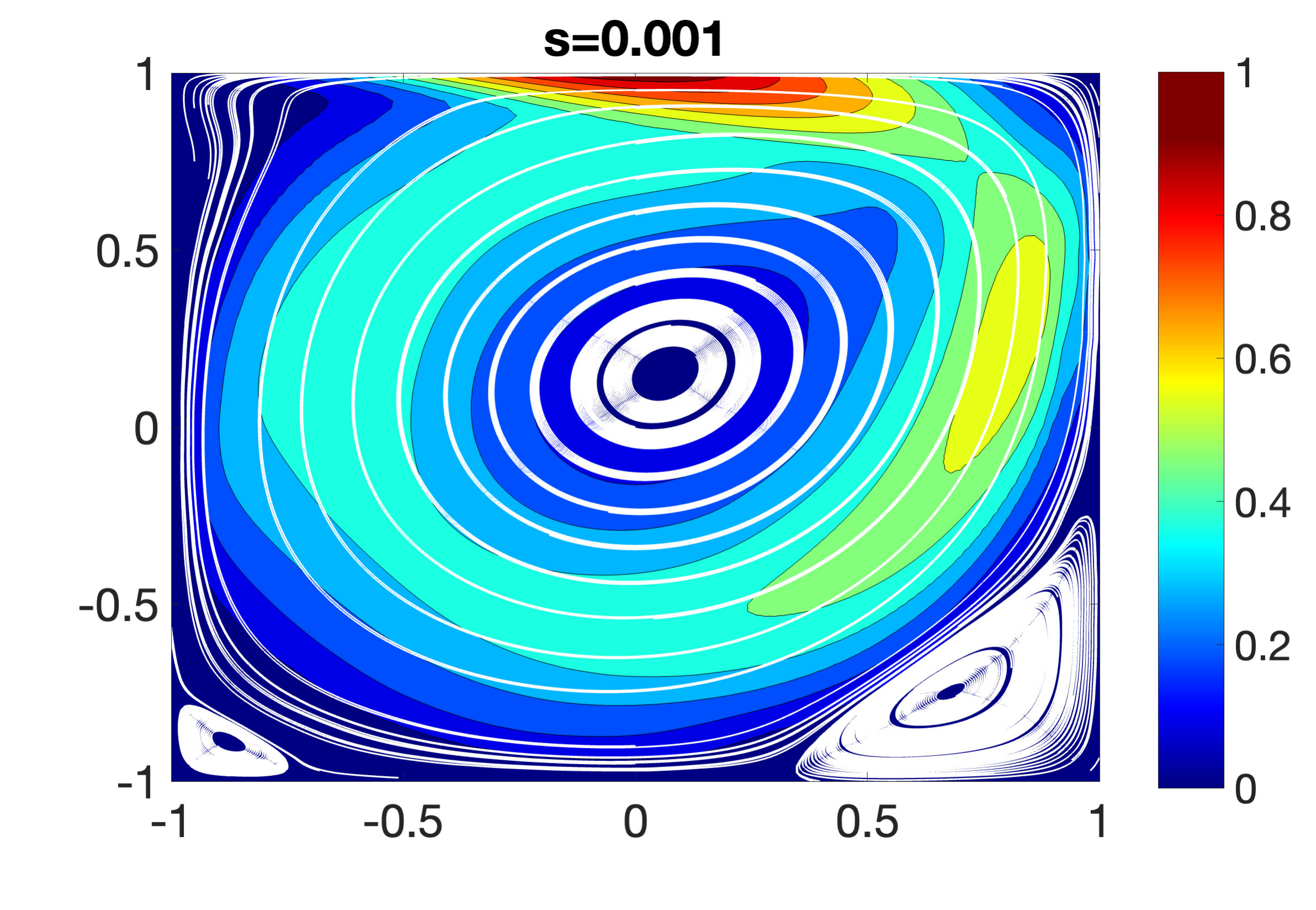}}
		{\includegraphics[width=0.28\textwidth,height=0.2\textwidth]{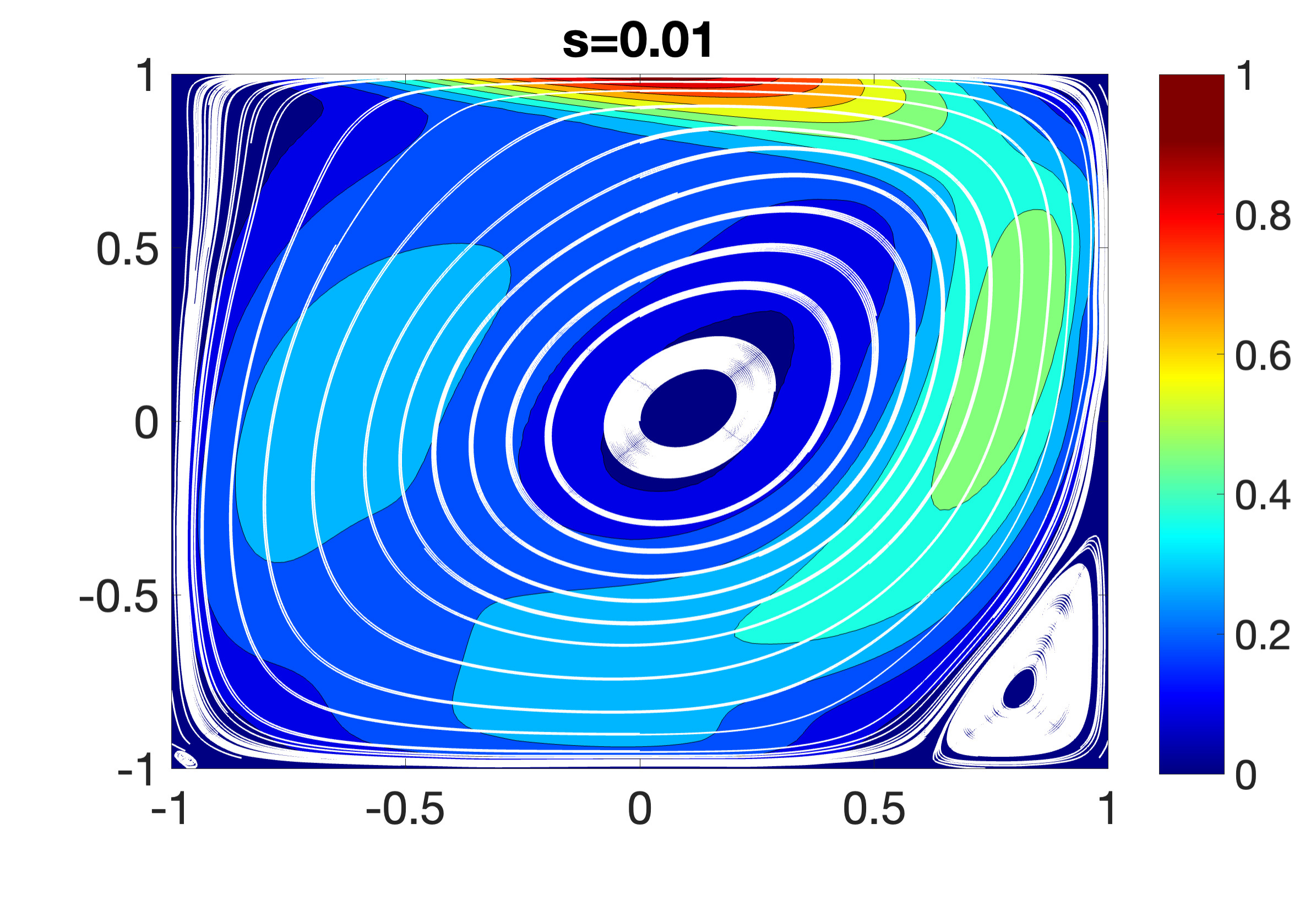}}\\
		{\includegraphics[width=0.28\textwidth,height=0.2\textwidth]{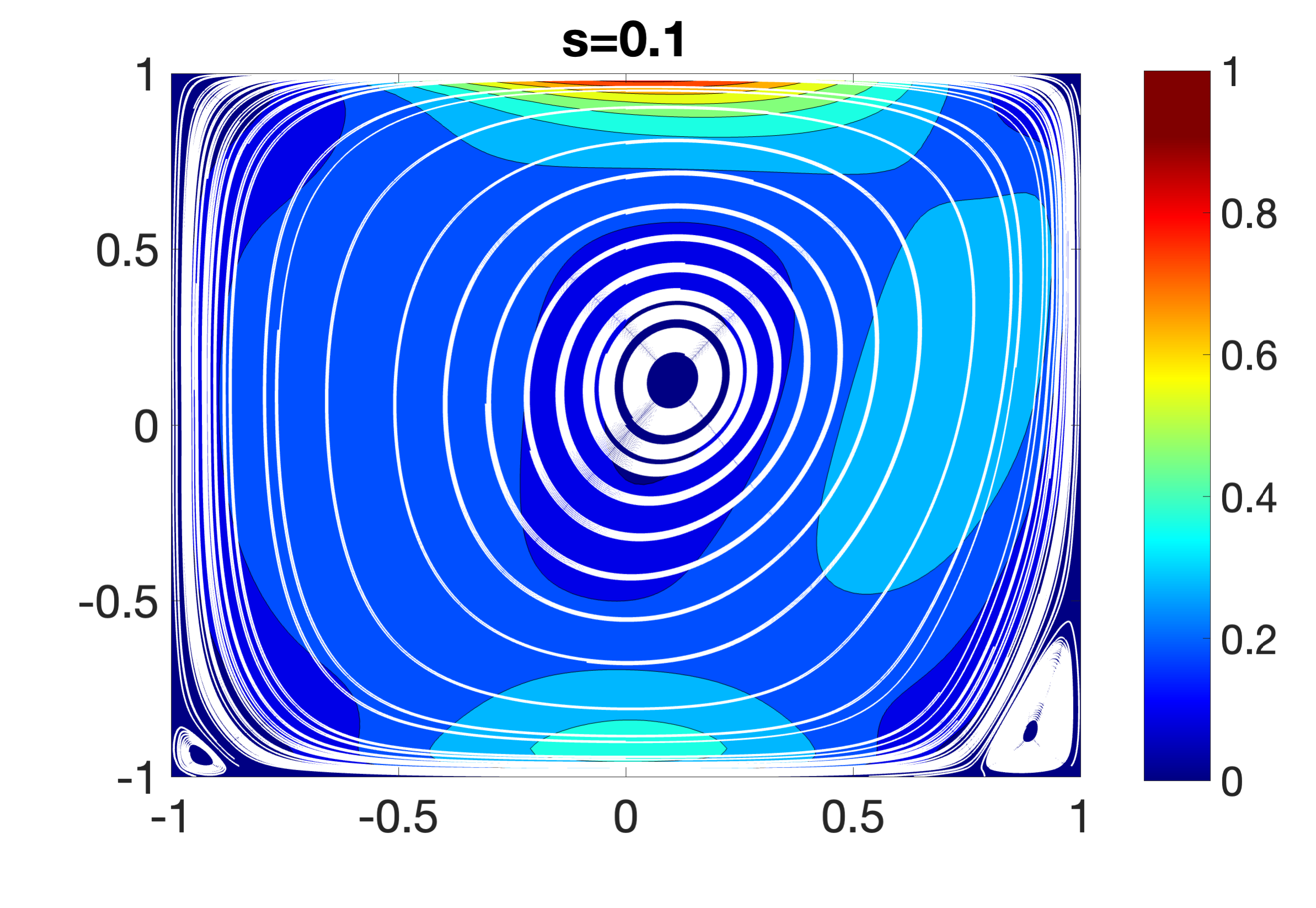}}
		{\includegraphics[width=0.28\textwidth,height=0.2\textwidth]{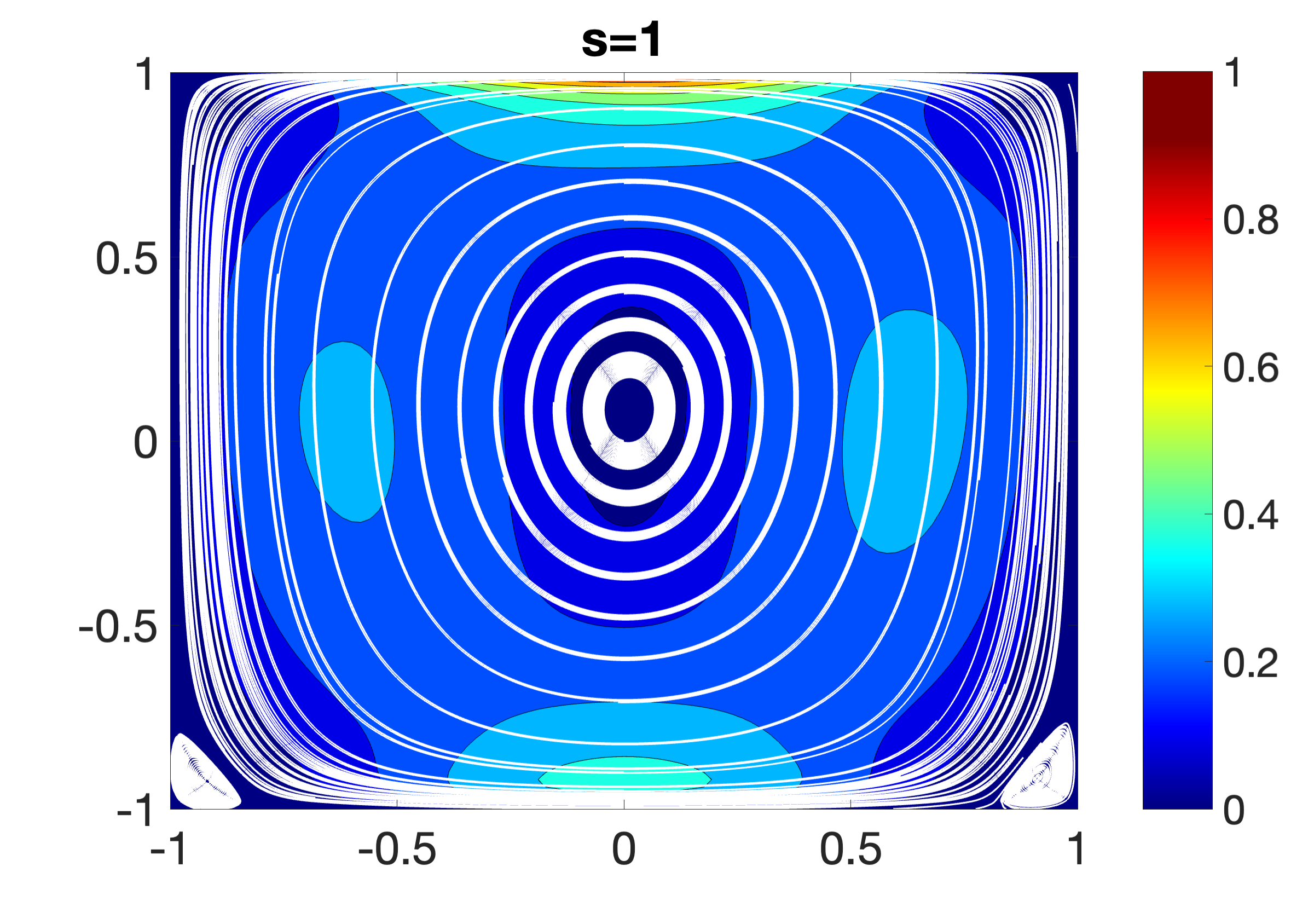}}\vspace{-4mm}
		\caption{Variable 5D random viscosity and magnetic diffusivity in lid-driven problem: Plot of speed contour for varying $s$ at 
			$t=600$ with $\gamma=10000$, $\E[Re]=15000$, and $\E[\nu_m]=0.01$.}	\label{noise-velocity}
	\end{figure}

\begin{figure} [ht]
		\centering		{\includegraphics[width=0.28\textwidth,height=0.2\textwidth]{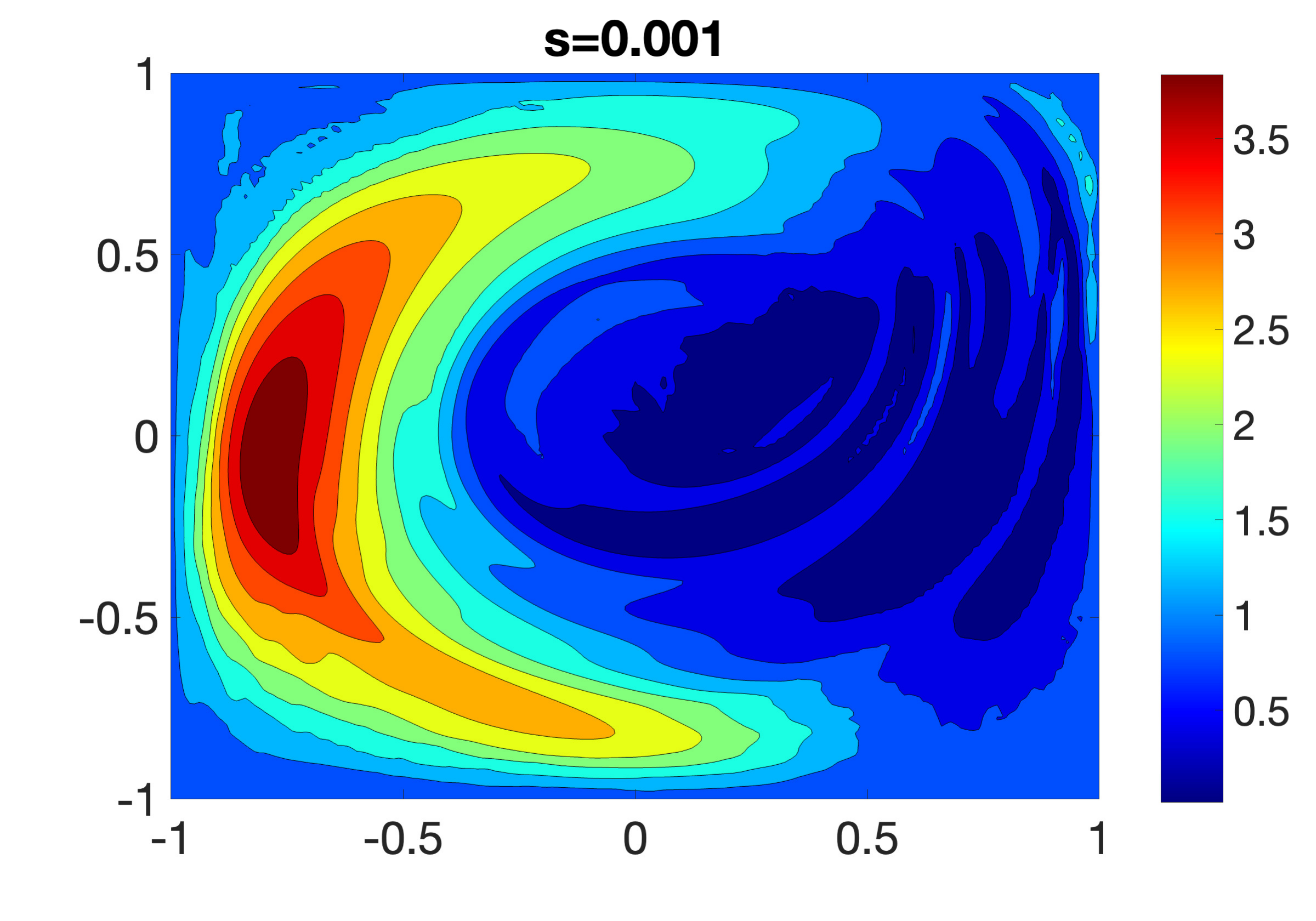}}
		{\includegraphics[width=0.28\textwidth,height=0.2\textwidth]{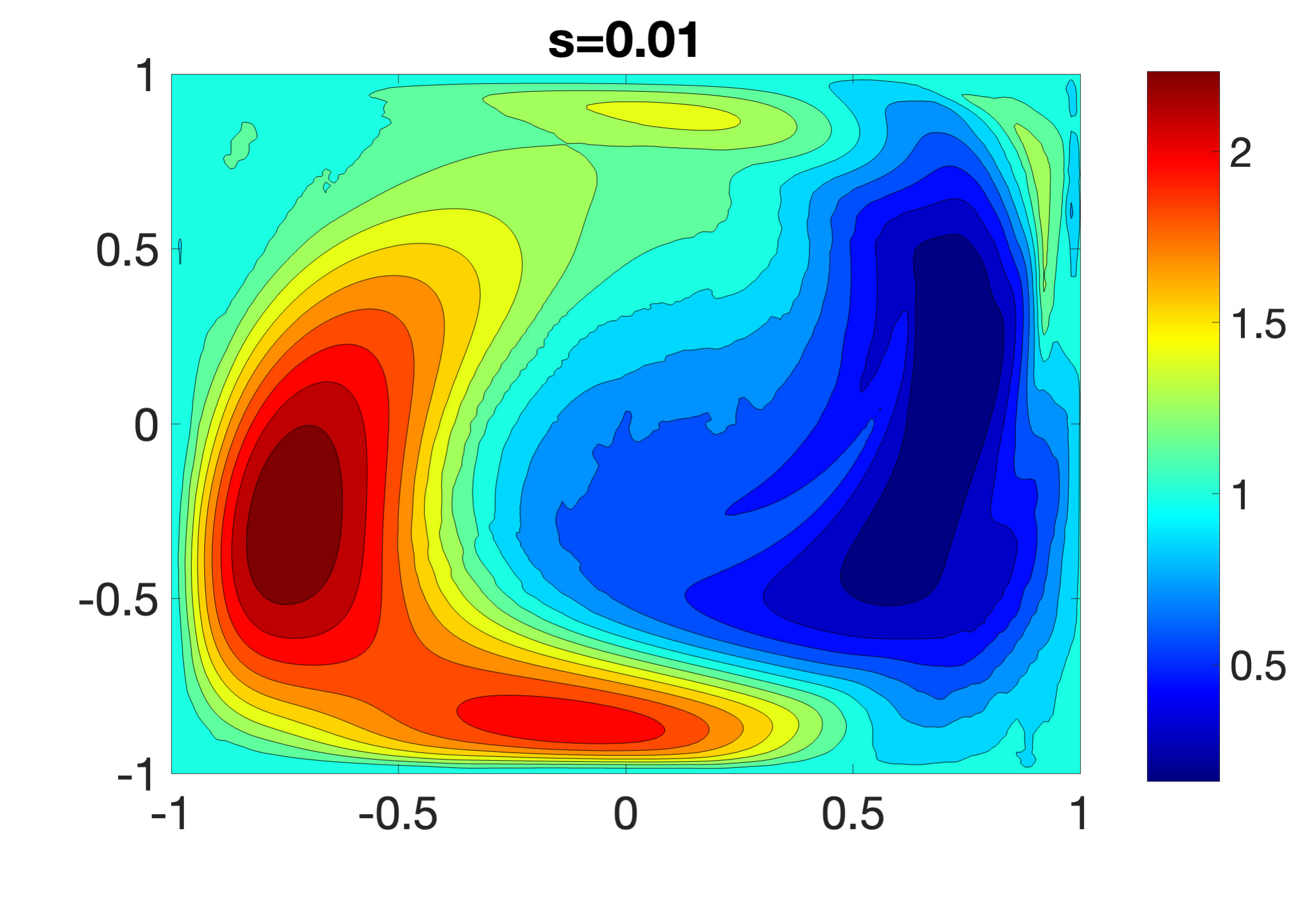}}\\
		{\includegraphics[width=0.28\textwidth,height=0.2\textwidth]{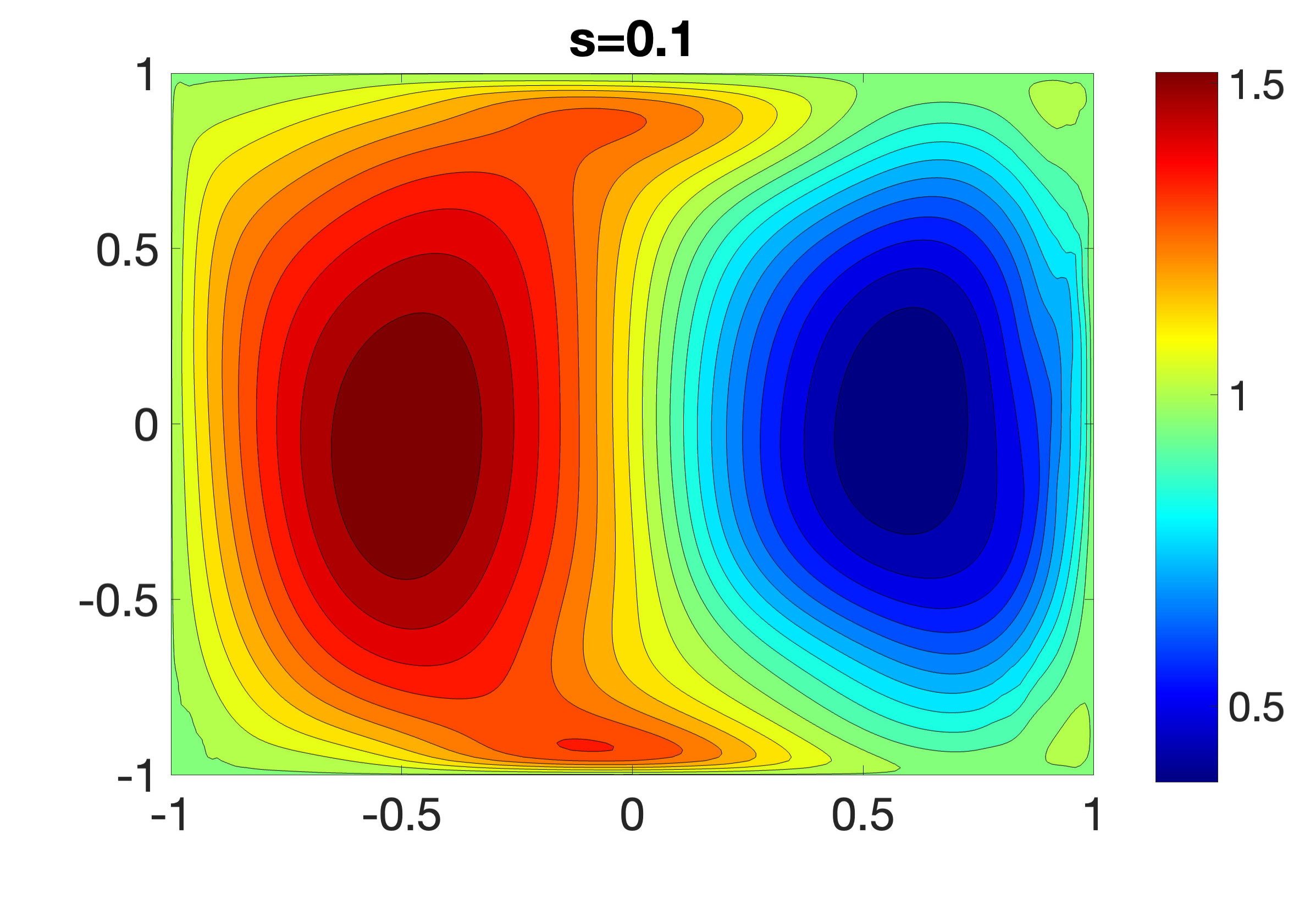}}
		{\includegraphics[width=0.28\textwidth,height=0.2\textwidth]{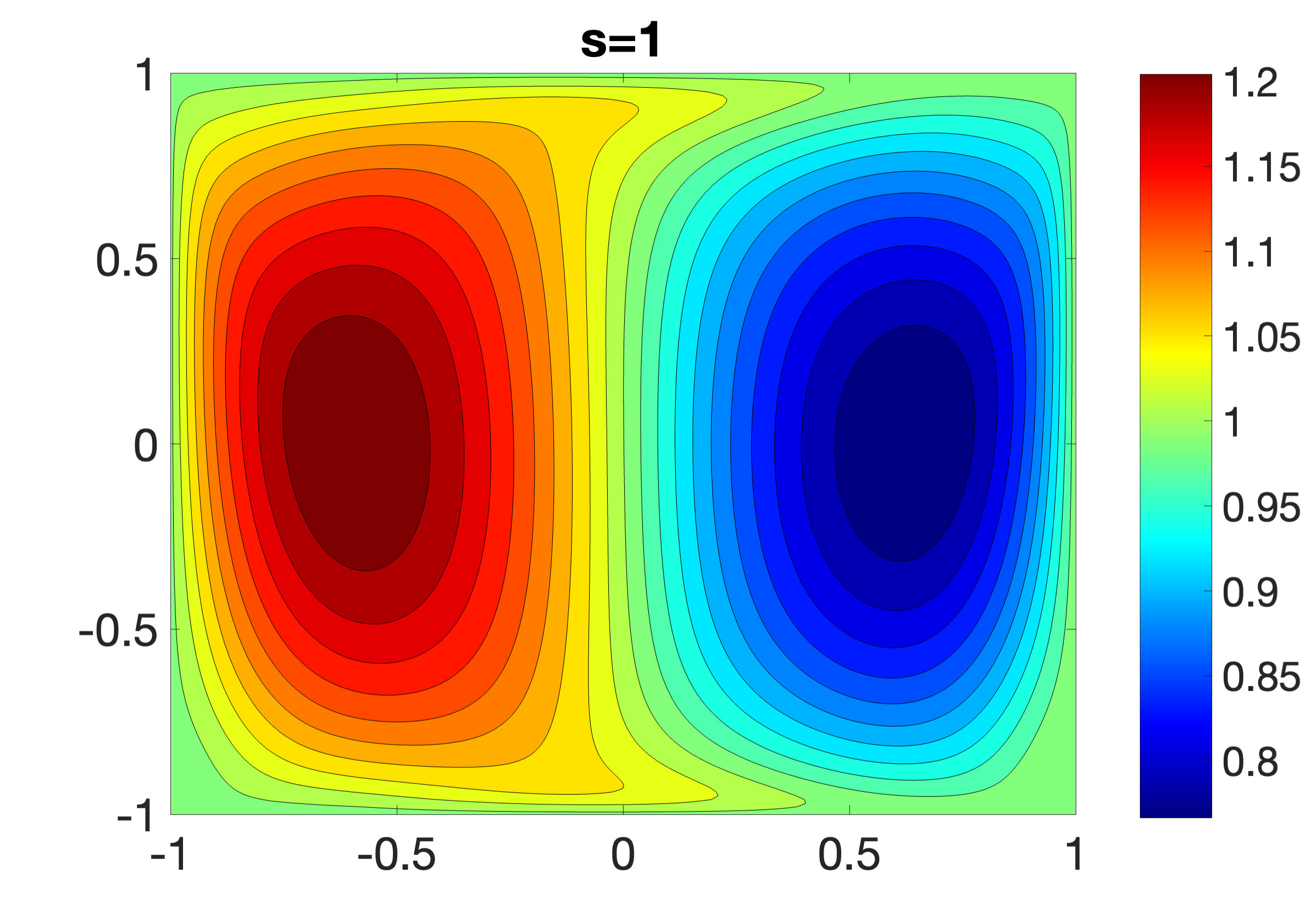}}
		\caption{Variable 5D random viscosity and magnetic diffusivity in lid-driven problem: Plot of magnetic field strength for varying 
			$s$ at $t=600$ with $\gamma=10000$, $\E[Re]=15000$, and $\E[\nu_m]=0.01$.}	\label{noise-magnetic}
	\end{figure}

	We now consider a 2D benchmark regularized lid-driven cavity problem \cite{balajewicz2013low,fick2018stabilized,lee2019study} with a domain 
	$\Omega=(-1,1)^2$. No-slip boundary conditions are applied to all sides except on the top wall (lid) of the cavity where we impose the following 
	boundary condition:
	\begin{align*}
		\bu_j|_{lid}=(1+k_j\epsilon){{(1-x_1^2)^2}\choose{0}}.
	\end{align*}
	On all sides of the cavity, we enforce the following the magnetic field boundary condition:
	\begin{align*}
		\bB_j=(1+k_j\epsilon){{0}\choose{1}}.
	\end{align*}
	The maximum speed of the lid is 1 and the characteristic length is $2$. In this experiment, we consider $\gamma=10000$, $\bif=\bg=0$, and 
	Clenshaw--Curtis sparse grid as the SCM, generated via the software package TASMANIAN \cite{stoyanov2015tasmanian,doecode_6305} with $N_{sc}=11$. 
	The generated computational barycentered refined mesh of the domain provides a total of 729840 degrees of freedom (dof), where velocity dof 
	$=324266$, magnetic field dof $=324266$, pressure dof $=40654$, and magnetic pressure dof $=40654$. The flow initiates from the state of rest in 
	absence of the magnetic flux density. 
	In this section, we consider the equations~\eqref{momentum}-\eqref{magnet-initial} with a random viscosity $\nu({\bx},{\by})$, and magnetic 
	diffusivity $\nu_m(\bx,\by)$, where ${\by}=(y_1,y_2,\cdots,y_d)\in\Gamma\subset\mathbb{R}^d$ is a higher-dimensional random variable, 
	$\mathbb{E}[\nu](\bx)=\frac{2c}{15000}$, and $\mathbb{E}[\nu_m](\bx)=0.01c$ for a suitable $c>0$, 
	$\mathbb{C}ov[\nu]({x},{x^{'}})=\frac{4}{15000^2}exp\left(-\frac{({x}-{x^{'}})^2}{l^2}\right)$, and $l$ is the correlation length. 
	This random field can be represented by the Karhunen-Lo\'eve expansion:
	\begin{align}
		\nu({\bx}, {\by})=\frac{2}{15000}\psi(\bx,\by),\;\;\text{and}\;\;\nu_m(\bx,\by)=\frac{1}{100}\psi(\bx,\by), \label{eq:var-vis}
	\end{align}
	where
	\begin{align*}
		\psi(\bx,\by)=\bigg(c+\left(\frac{\sqrt{\pi}l}{2}\right)^{\frac12}y_{1}(\omega)+\sum_{j=1}^{q}\sqrt{\xi_j}&\bigg(\sin\left(\frac{j\pi 
			x_1}{2}\right)\sin\left(\frac{j\pi x_2}{2}\right)y_{2j}(\omega)\nonumber\\&+\cos\left(\frac{j\pi x_1}{2}\right)\cos\left(\frac{j\pi 
			x_2}{2}\right)y_{2j+1}(\omega)\bigg),
	\end{align*}
	in which the infinite series is truncated up to the first $q$ terms. 
	The uncorrelated random variables $y_{j}$ have zero mean and unit variance, and the eigenvalues are equal to
	$$\sqrt{\xi_j}=(\sqrt{\pi}l)^{\frac12}exp\left(-\frac{(j\pi l)^2}{8}\right).$$
	For our test problem, we consider the random variables $y_{j}(\omega)\in[-\sqrt{3},\sqrt{3}]$, the correlation length $l=0.01$, $d=5$, $c=1$, and $q=2$. We run the simulation with time-step size $\Delta t=5$ until the simulation end time $T=600$ for various values of the coupling parameter $s$ together with the perturbation parameter $\epsilon=0.01$ in the initial and boundary conditions. The Fig.~\ref{noise-velocity}-\ref{noise-magnetic} illustrate the velocity solution as the speed contour, and the magnetic field strength for $s=0.001,\;0.01,\;0.1,$ and $1$ and are the outcomes of the SCM-SPP-SMHD scheme given in Algorithm \ref{SPP-FEM}. As $s$ grows, the impact of the Lorentz force gets stronger in the flow field, which in turn slow down the evolve over time process. This can be observed as the speed and the size of the vorticities get reduced in Fig.~\ref{noise-velocity} while the magnetic field strength realizes a type of reflection symmetry in Fig.~\ref{noise-magnetic}.

	\begin{figure} [ht]
		\centering
		{\includegraphics[width=0.5\textwidth,height=0.30\textwidth]{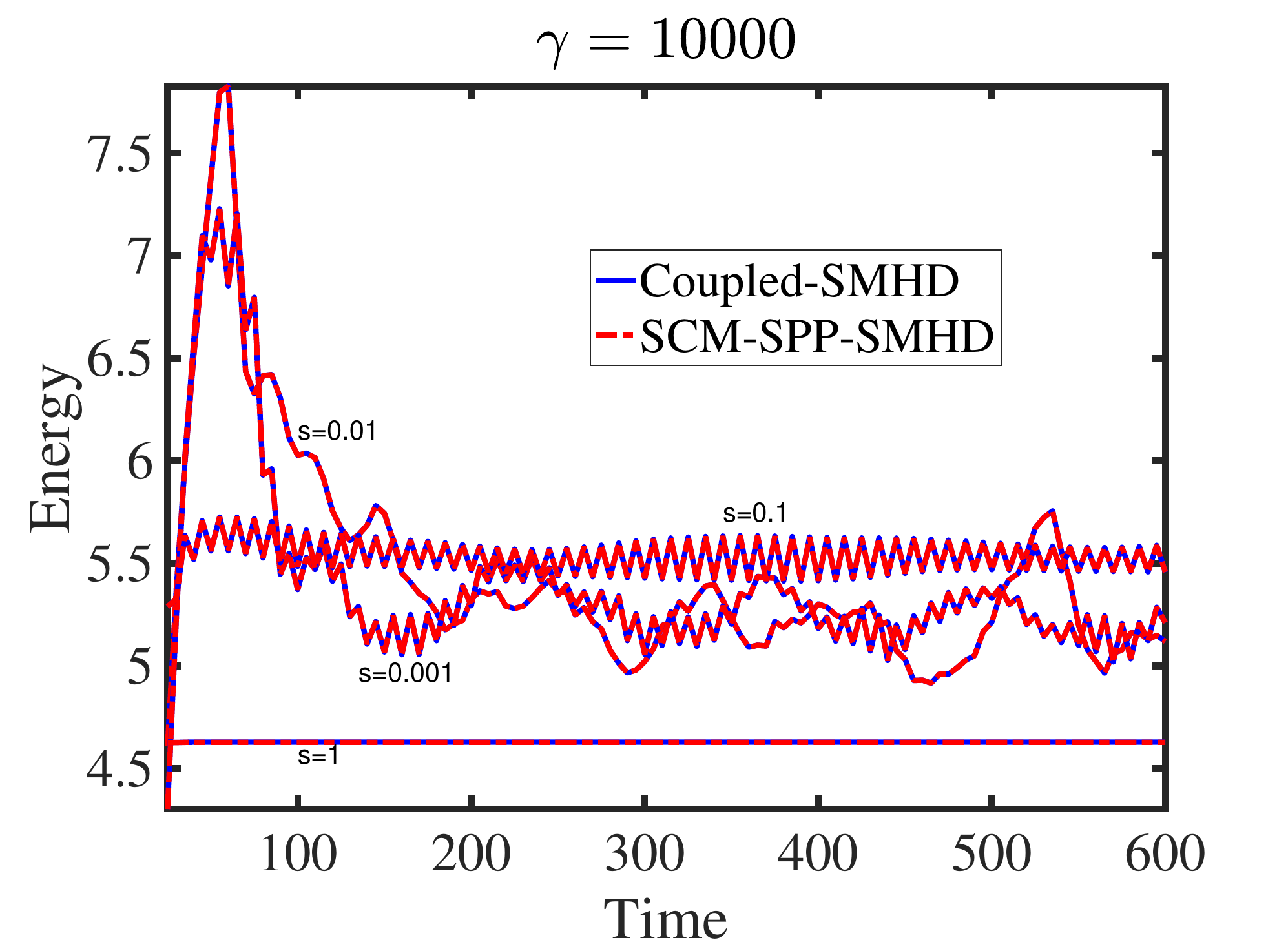}}
		\caption{Variable 5D random viscosity and magnetic diffusivity in lid-driven problem: Plot of energy vs. time for various $s$ 
			with $\gamma=10000$, $\mathbb{E}[Re]=15000$, and $\mathbb{E}[\nu_m]=0.01$.}\label{energy-plot-lid-driven-cavity}
	\end{figure}
	For this experiment, in Fig. \ref{energy-plot-lid-driven-cavity}, we plot the system energy vs. time for various values of $s$ using both the 
	SCM-SPP-SMHD and Coupled-SMHD methods. To compare the two models, we compute the weighted mean energy at time $t=t_n$, which is defined as the weighted average of $\frac12\|\bu(t_n,\bx,\by^j)\|$ for all sample points. We found excellent agreements between the energy plots from the solution of the Coupled-SMHD scheme and 
	the solution of penalty projection based SCM-SPP-SMHD method with $\gamma=10000$, which support the theory.
	
	\section{Conclusion and future works}\label{conclusion-future-works}
	In this paper, we propose, analyze, and test an efficient and accurate grad-div stabilized penalty-projection SCM-SPP-SMHD scheme in conjunction with SCM for solving stochastic MHD flow problems. The intriguing algorithm has several features that make it efficient and accurate: (1) The use of Els\"asser variables formulation allows for a stable decoupling of the coupled PDEs, (2) A discrete Hodge decomposition is used for decoupling further which allows to use two much easier linear solves instead of using a difficult solve of the saddle point problems for each realization at each time-step, (3) The four sub-problems are designed in an elegant way that at each time-step, the system matrix remains common to all realizations but with different right-hand-side vectors, which saves a huge computer memory and assembly time of assembling several global different system matrices; Furthermore, this allows to take the advantage of using block linear solvers, (4) The use of ensemble eddy-viscosity terms provide stability of flows that are not resolved on particular meshes. (5) The large (but optimal) coefficient of the grad-div stabilization parameter provides accuracy of the splitting algorithm equivalent to a coupled scheme, and (6) The sparse grid SCM wrapper helps to use fewer realization.
	
	The SCM-SPP-SMHD algorithm is rigorously proven to be stable and converges to the equivalent coupled Coupled-SMHD ensemble scheme for large grad-div stabilization parameters. The numerical test verifies the first-order convergence of the SCM-SPP-SMHD scheme to the Coupled-SMHD scheme.  The optimal spatial and first-order temporal convergence rates of the scheme are verified with synthetic data for analytical test problems with random noise in the parameters values. We implement the scheme on benchmark channel flow over a step problem and a regularized lid-driven cavity problem with space-dependent 5D random high Reynolds and high magnetic Reynolds numbers. We found the efficient SCM-SPP-SMHD scheme performs well with high grad-div stabilization parameters. This penalty-projection-based efficient algorithm will be an enabling tool for large-scale simulation of complex 3D MHD problems. 
	
	In the future, we will implement this scheme on the 3D Taylor-Green vortex problem, and examine its performance together with various linear solvers and appropriate preconditioners. As a future work, we will propose, analyze, and test first- and second-order accurate time-stepping penalty-projection schemes for the UQ of N-S flow problems following the work in \cite{linke2017connection}. An evolve-filter-relax stabilized Reduced Order (ROM) SCM for the time-dependent MHD flow will be proposed following the recent work in \cite{gunzburger2019evolve}.
	
	\section{Acknowledgement} The National Science Foundation (NSF) is acknowledged for supporting this research through the grant DMS-221327. We also acknowledge the Texas A\&M International University for providing logistic support. The authors also thank Dr. Leo G. Rebholz for sharing his thoughts which greatly improved the manuscript.
	\bibliographystyle{plain}
	\bibliography{Penalty-Pro}
\end{document}